\documentclass[12pt]{amsart}

\usepackage{latexsym}
\usepackage{amssymb}
\usepackage{amsmath}
\usepackage{amscd}
\usepackage{amsthm}
\usepackage[all]{xy}  

\numberwithin{equation}{section}

\newtheoremstyle{fancy1}{10pt}{10pt}{\itshape}{12pt}{\textsc\bgroup}{.\egroup}{8pt}{
}
\newtheoremstyle{fancy2}{10pt}{10pt}{}{12pt}{\itshape}{.}{8pt}{ }

\theoremstyle{fancy1}

\newtheorem{cor}[equation]{Corollary}
\newtheorem{lem}[equation]{Lemma}
\newtheorem{prop}[equation]{Proposition}

\newtheorem*{thethm*}{Theorem A}

\newtheorem{main}{Theorem}
\newtheorem*{main*}{Result}

\newtheorem*{cor*}{Corollary}

\theoremstyle{fancy2}
\newtheorem{defn}[equation]{Definition}
\newtheorem*{def*}{Definition}

\newtheorem*{rem*}{Remark}

\newtheorem*{example*}{Example}

\newtheorem*{examples*}{Examples}

\theoremstyle{remark}
\newtheorem*{case}{Case}

\newcommand{\cref}[1]{Corollary~\ref{#1}}

\newcommand{\tref}[1]{Theorem~\ref{#1}}

\newcommand{\CP}{\mathbb{C\mkern1mu P}}
\newcommand{\HP}{\mathbb{H\mkern1mu P}}
\newcommand{\CaP}{\mathrm{Ca}\mathbb{\mkern1mu P}^2}

\newcommand{\C}{{\mathbb{C}}}
\newcommand{\R}{{\mathbb{R}}}
\newcommand{\Z}{{\mathbb{Z}}}
\renewcommand{\H}{{\mathbb{H}}}

\newcommand{\F}{\ensuremath{\operatorname{\mathsf{F}}}}
\newcommand{\G}{\ensuremath{\operatorname{\mathsf{G}}}}

\newcommand{\SO}{\ensuremath{\operatorname{\mathsf{SO}}}}

\newcommand{\Sp}{\ensuremath{\operatorname{\mathsf{Sp}}}}
\newcommand{\U}{\ensuremath{\operatorname{\mathsf{U}}}}
\newcommand{\SU}{\ensuremath{\operatorname{\mathsf{SU}}}}
\newcommand{\Spin}{\ensuremath{\operatorname{\mathsf{Spin}}}}

\renewcommand{\S}{\ensuremath{\operatorname{\mathsf{S}}}}

\newcommand{\fp}{{\mathfrak{p}}}

\newcommand{\fsu}{{\mathfrak{su}}}

\def\con#1=#2(#3){#1 \equiv #2 \bmod{#3}}

\newcommand{\set}[1]{\left\{#1\right\}}

\newcommand{\tr}{\ensuremath{\operatorname{tr}}}
\newcommand{\diag}{\ensuremath{\operatorname{diag}}}

\newcommand{\rank}{\ensuremath{\operatorname{rk}}}

\renewcommand{\Im}{\ensuremath{\operatorname{Im}}}
\newcommand{\Ad}{\ensuremath{\operatorname{Ad}}}

\newcommand{\proj}{\ensuremath{\operatorname{proj}}}
\newcommand{\Isom}{\ensuremath{\operatorname{Isom}}}
\newcommand{\del}{\partial}

 \DeclareMathOperator{\id}{id}
\DeclareMathOperator{\Diff}{Diff}
\DeclareMathOperator{\Id}{Id}

\newcommand{\wbar}{\overline}  

\renewcommand{\F}{\mathsf{F}}

\begin{document}


\title[Cohomogeneity One Manifolds in Low Dimensions]{Classification of Cohomogeneity One Manifolds in Low Dimensions}

\author{Corey A. Hoelscher}

\begin{abstract}
A cohomogeneity one manifold is a manifold with the action of a compact Lie group, whose quotient is one dimensional. Such manifolds are of interest in Riemannian geometry, in the context of nonnegative sectional curvature, as well as in other areas of geometry and in physics. In this paper we classify compact simply connected cohomogeneity one manifolds in dimensions 5, 6 and 7. We also show that all such manifolds admit metrics of nonnegative sectional curvature, with the possible exception of two families of manifolds.
\end{abstract}

\maketitle

%
%
%
%
%
%

\renewcommand{\thetable}{\Roman{table}}

\section*{Introduction}

Manifolds with non-negative curvature have always played a special role in Riemannian geometry and yet finding new examples of such manifolds is particularly difficult. Recently in \cite{GZ1}, Grove-Ziller constructed a large class of non-negatively curved metrics on certain cohomogeneity one manifolds, that is, manifolds with an action by a compact Lie group whose orbit space is one dimensional. In particular they showed that all principal $S^3$ bundles over $S^4$ can be written as cohomogeneity one manifolds with metrics of non-negative sectional curvature. It was also shown in \cite{GZ2} that every compact cohomogeneity one manifold admits a metric of non-negative Ricci curvature and admits a metric of positive Ricci curvature if and only if its fundamental group is finite. So cohomogeneity one manifolds provide a good setting for examples of manifolds with certain curvature restrictions. Cohomogeneity one actions are also of independent interest in the field of group actions as they are the simplest examples of inhomogeneous actions. They also arise in physics as new examples of Einstein and Einstein-Sasaki manifolds (see \cite{Conti} or \cite{GHY}) and as examples of manifolds with $G_2$ and $\Spin(7)$-holonomy (see \cite{CS} and \cite{CGLP}). It is then an interesting question how big the class of cohomogeneity one manifolds is. Such manifolds where classified in dimensions 4 and lower in \cite{Parker} and \cite{Ne}. Physicists are interested in those of dimension 5, 7 and 8 and many of the most interesting examples appearing in \cite{GZ1} were 7-dimensional. The main result in this paper is a classification of compact simply connected cohomogeneity one manifolds in dimensions 5, 6 and 7.

Before we state the theorem we will review some basic facts about cohomogeneity one manifolds. Recall that a compact cohomogeneity one manifold with finite fundamental group has a description in terms of its \emph{group diagram}
  $$G\,\,\supset \,\, K^-\!,K^+ \,\, \supset \,\, H$$
where $G$ is the group that acts, which is assumed to be compact, $H$ is a principal isotropy subgroup, and $K^\pm$ are certain nonprincipal isotropy subgroups which both contain $H$ (see Section \ref{basic_structure} for details). We will henceforth describe actions in terms of their corresponding group diagrams.

If the group $G$ is disconnected then the identity component still acts by cohomogeneity one. Further, since the isometry group of a compact Riemannian manifold is a compact Lie group, it is natural to restrict our attention to actions by compact groups. So we will always assume that $G$ is compact and connected.

There is a simple way to build cohomogeneity one actions from lower dimensional actions by taking products. Say $G_1$ acts by cohomogeneity one on $M_1$ and $G_2$ acts transitively on $M_2$. Then it is clear that $G_1\times G_2$ acts by cohomogeneity one on  $M_1\times M_2$, as a product. Such actions are referred to as \emph{product actions}.

We call a cohomogeneity one action of $G$ on $M$ \emph{reducible} if there is a proper normal subgroup of $G$ that still acts by cohomogeneity one with the same orbits. This gives a way of reducing these actions to simpler actions. Conversely there is a natural way of extending an arbitrary cohomogeneity one action to an action by a possibly larger group. Such extensions, called \emph{normal extensions}, are described in more detail in Section \ref{EaR}. It turns out that every reducible action is a normal extension of its restricted action. Therefore it is natural to restrict ourselves to nonreducible actions in the classification.

Recall that a cohomogeneity one action is \emph{nonprimitive} if all the isotropy subgroups, $K^-$, $K^+$ and $H$ for some group diagram representation, are all contained in some proper subgroup $L$ in $G$. Such a nonprimitive action is well known to be equivalent to the usual $G$ action on $G\times_L M_L$, where $M_L$ is the cohomogeneity one manifold given by the group diagram $L\supset K^-, K^+\supset H$. With these definitions in place, we are ready to state the main result.

\begin{main}\label{simpleclass}
Every nonreducible cohomogeneity one action on a compact simply connected manifold of dimension 5, 6 or 7 by a compact connected group is equivalent to one of the following:
\begin{enumerate}
\item[i.] an isometric action on a symmetric space,
\item[ii.] a product action,
\item[iii.] the $\SO(2)\SO(n)$ action on the Brieskorn variety, $B^{2n-1}_d$,
\item[iv.] one of the primitive actions listed in Table \ref{class:prim} or a nonprimitive action from Table \ref{class:nonprim}.
\end{enumerate}
Hence every cohomogeneity one action on such a manifold by a compact connected group is a normal extension of one of these actions.
\end{main}

\begin{rem*}
When reading the tables below, we observe the following conventions and notations. In all cases we denote $H_\pm:=H\cap K^\pm_0$. $H_+$ is either $H_0$ in the case that $\dim K^+/H >1$, or $H_0\cdot \Z_n$, for some $n\in \Z$, in the case $\dim K^+/H=1$, and similarly for $H_-$. Here and throughout, $L_0$ denotes the identity component of a given Lie group $L$. Next, we always know that $H\subset K^\pm$, and this puts some unstated restrictions on the groups in the tables. Furthermore, $H_0$ should be understood to be trivial unless otherwise stated. In the tables, we also assume, when we have a group of the form $\set{(e^{ip\theta},e^{iq\theta})}$, that $\gcd(p,q)=1$, and similarly for other such groups. Furthermore $\theta, \phi\in \R$ and $z,w\in S^1\subset\C\subset\H$ are taken to vary arbitrarily, while integers $a,b,c,m,n,p,q,r,s,\lambda,\mu$ are understood to be fixed within a given group diagram. Finally, $i,j,k\in S^3\subset \H$ are the usual unit quaternions.

Notice that many of the diagrams are not effective, since $G$ and $H$ share a finite normal subgroup. We have allowed this possibility so that our descriptions are simpler. The effective version of each action can always be determined by quotienting each group in the diagram by $Z(G)\cap H$, where $Z(G)$ is the center of $G$ (see Section \ref{basic_structure} for details).
\end{rem*}

In Section \ref{others}, we collect some basic facts about each family of actions in Tables \ref{class:prim} and \ref{class:nonprim}. This section would be of interest to the reader who wants to quickly know what can be easily said about these actions. For example some of these actions are of types (i.), (ii.) and (iii.) of Theorem \ref{simpleclass} for special choices of parameters. We also describe the isometric actions on symmetric spaces from Theorem \ref{simpleclass} in Section \ref{symspactions} and the Brieskorn actions in Section \ref{brieskorn}.

In the process of proving this theorem we get a complete list of the possible nonreducible actions in these dimensions. Many of these actions are either sum actions or product actions. Such actions are easily understood and easily identified from their group diagrams, as described in Section \ref{special_types} and summarized in Table \ref{SaP}. The remainder of the actions are listed in Tables \ref{5dimlist} through \ref{7dimlist2} of the appendix, for the convenience of the reader.

{\setlength{\tabcolsep}{0.2cm}
\renewcommand{\arraystretch}{1.6}
\stepcounter{equation}
\begin{table}[!h]
\begin{center}
\begin{tabular}{|c|c|}
\hline
$P^5$ & $S^3\times S^1 \,\, \supset \,\, \set{(e^{i\theta},1)}\cdot H, \,\, \set{(e^{jp\theta}, e^{i\theta})}  \,\, \supset \,\,  \langle(j,i)\rangle$\\
 &  where $p\equiv 1 \mod 4$\\
\hline
$P^7_A$ & $S^3\times S^3 \,\, \supset \,\, \set{(e^{ip_-\theta},e^{iq_-\theta})}, \,\, \set{(e^{jp_+\theta},e^{jq_+\theta})}\cdot H \,\, \supset \,\,  \langle (i,i) \rangle$\\
 & where $p_-,q_-\equiv 1 \mod 4$\\
\hline
$P^7_B$ & $S^3\times S^3 \,\, \supset \,\, \set{(e^{ip_-\theta},e^{iq_-\theta})} \cdot H, \,\, \set{(e^{jp_+\theta},e^{jq_+\theta})}\cdot H \,\, \supset \,\,  \langle (i,i), (1,-1) \rangle$\\
 & where $p_-,q_-\equiv 1 \mod 4$, $p_+$ even\\
\hline
$P^7_C$ & $S^3\times S^3 \,\, \supset \,\, \set{(e^{ip_-\theta},e^{iq_-\theta})} \cdot H, \,\, \set{(e^{jp_+\theta},e^{jq_+\theta})}\cdot H \,\, \supset \,\,  \Delta Q$\\
 & where $p_\pm,q_\pm\equiv 1 \mod 4$\\
\hline
$P^7_D$ & $S^3\times S^3 \,\, \supset \,\, \set{(e^{ip\theta},e^{iq\theta})}, \,\, \Delta S^3 \cdot \Z_n \,\, \supset \,\,  \Z_n$\\
 & where $n=2$ and $p$ or $q$ even; or $n=1$ and $p$ and $q$ arbitrary\\
\hline
\end{tabular}
\end{center}
\vspace{0.1cm}
\caption{Primitive cohomogeneity one manifolds of Theorem \ref{simpleclass}}\label{class:prim}
\end{table}}
{\setlength{\tabcolsep}{0.2cm}
\renewcommand{\arraystretch}{1.302}
\stepcounter{equation}
\begin{table}[!h]
\begin{center}
\begin{tabular}{|c|c|}
\hline
$N^5$ & $S^3\times S^1 \,\, \supset \,\, \set{(e^{ip_-\theta},e^{iq_-\theta})}\cdot H, \,\, \set{(e^{ip_+\theta},e^{iq_+\theta})}\cdot H  \,\, \supset \,\,  H_-\cdot H_+$\\
 &  $K^-\ne K^+$, $(q_-,q_+)\ne \mathbf{0}$,  $\gcd(q_-,q_+,d)=1$\\
 &  where $d= \#(K^-_0\cap K^+_0)/\#(H\cap K^-_0\cap K^+_0)$\\
\hline
$N^6_A$ & $S^3\times T^2 \,\, \supset \,\, \set{(e^{i a_-\theta},e^{ib_-\theta},e^{ic_-\theta})}\cdot H, \,\, \set{(e^{i a_+\theta},e^{ib_+\theta},e^{ic_+\theta})}\cdot H \,\, \supset \,\, H$\\
 & where $K^-\ne K^+$, $H=H_-\cdot H_+$, $\gcd(b_\pm,c_\pm)=1$,\\
 & $a_\pm=rb_\pm+sc_\pm$, and $K^-_0\cap K^+_0\subset H$\\
\hline
$N^6_B$ & $S^3\times S^3 \,\, \supset \,\, \set{(e^{i\theta},e^{i\phi})}, \,\, \set{(e^{i\theta},e^{i\phi})} \,\, \supset \,\, \set{(e^{ip\theta},e^{iq\theta})} \cdot \Z_n$\\
\hline
$N^6_C$ & $S^3\times S^3 \,\, \supset \,\, T^2, \,\, S^3\times \Z_n \,\, \supset \,\, S^1 \times \Z_n$\\
\hline
$N^6_D$ & $S^3\times S^3 \,\, \supset \,\, T^2, \,\, S^3\times S^1 \,\, \supset \,\, \set{(e^{ip\theta},e^{i\theta})}$\\
\hline
$N^6_E$ & $S^3\times S^3 \,\, \supset \,\, S^3\times S^1, \,\, S^3\times S^1 \,\, \supset \,\, \set{(e^{ip\theta},e^{i\theta})}$\\
\hline
$N^6_F$ & $\SU(3) \,\, \supset \,\, \S(\U(2)\U(1)), \,\, \S(\U(2)\U(1)) \,\, \supset \,\, \SU(2)\SU(1)\cdot \Z_n$\\
\hline
$N^7_A$ & $S^3\times S^3 \,\, \supset \,\, \set{(e^{ip_-\theta},e^{iq_-\theta})}\cdot H_+, \,\, \set{(e^{ip_+\theta},e^{iq_+\theta})} \cdot H_- \,\, \supset \,\,  H_-\cdot H_+$\\
\hline
$N^7_B$ & $S^3\times S^3 \,\, \supset \,\, \set{(e^{ip\theta},e^{iq\theta})}\cdot H_+, \,\, \set{(e^{j\theta},1)}\cdot H_- \,\, \supset \,\,  H_-\cdot H_+$\\
 & where $H_\pm = \Z_{n_\pm} \subset K^\pm_0$, $n_+\le 2$, $4|n_-$ and $p_- \equiv \pm \frac{n_-}{4} \mod n_-$\\
\hline
$N^7_C$ & $S^3\times S^3 \,\, \supset \,\, \set{(e^{ip\theta},e^{iq\theta})}, \,\, S^3\times \Z_n \,\, \supset \,\,  \Z_n$\\
 & where $(q,n)=1$\\
\hline
$N^7_D$ & $S^3\!\!\times\! S^3\!\!\times\! S^1 \supset  \set{(z^pw^{\lambda m}\!,z^qw^{\mu m}\!,w)},   \set{(z^pw^{\lambda m}\!,z^qw^{\mu m}\!,w)}   \supset  H_0\cdot\Z_n $\\
 & where $H_0=\set{(z^p,z^q,1)}$, $p\mu-q\lambda=1$ and $\Z_n\subset \set{(w^{\lambda m}\!,w^{\mu m}\!,w)}$\\
\hline
$N^7_E$ & $S^3\!\!\times\! S^3\!\!\times\! S^1\! \supset\!  \set{\!(z^pw^{\lambda m_-}\!,z^qw^{\mu m_-}\!,w^{n_-})\!}\! H,   \set{\!(z^pw^{\lambda m_+}\!,z^qw^{\mu m_+}\!,w^{n_+})\!}\! H \! \supset \! H$\\
 & where $H=H_-\cdot H_+$, $H_0=\set{(z^p,z^q,1)}$, $K^-\ne K^+$, $p\mu-q\lambda=1$,\\
 & $\gcd(n_-,n_+,d)=1$ where $d$ is the index of $H\cap K^-_0 \cap K^+_0$ in $K^-_0 \cap K^+_0$\\
\hline
$N^7_F$ & $S^3 \! \times \! S^3 \! \times \! S^1 \,\, \supset \,\, \set{(e^{ip\phi}e^{ia\theta},e^{i\phi},e^{i\theta})}, \,\, S^3\! \times \! S^1\! \times\! \Z_n  \,\, \supset \,\, \set{(e^{ip\phi},e^{i\phi},1)}\cdot \Z_n$\\
 & $\Z_n\subset \set{(e^{ia\theta},1,e^{i\theta})}$\\
\hline
$N^7_G$ & $\SU(3) \,\, \supset \,\, \S(\U(1)\U(2)), \,\, \S(\U(1)\U(2)) \,\, \supset \,\, T^2$\\
\hline
$N^7_H$ & $\SU(3)\times S^1 \,\, \supset \,\, \set{(\beta(m_-\theta),e^{in_-\theta})}\cdot H, \,\, \set{(\beta(m_+\theta),e^{in_+\theta})}\cdot H  \,\, \supset \,\, H$\\
 & $H_0 = \SU(1)\SU(2)\times 1$, $H=H_-\cdot H_+$, $K^-\ne K^+$,\\
 & $\beta(\theta)= \diag(e^{-i\theta},e^{i\theta},1)$, $\gcd(n_-,n_+,d)=1$\\
 & where $d$ is the index of $H\cap K^-_0 \cap K^+_0$ in $K^-_0 \cap K^+_0$\\
\hline
$N^7_I$ & $\Sp(2)\,\, \supset \,\, \Sp(1)\Sp(1), \,\, \Sp(1)\Sp(1) \,\, \supset \,\, \Sp(1)\SO(2)$\\
\hline
\end{tabular}
\end{center}
\vspace{0.1cm}
\caption{Nonprimitive cohomogeneity one manifolds of Theorem \ref{simpleclass}.}\label{class:nonprim}
\end{table}}

The next theorem addresses the issue of nonnegative sectional curvature. In \cite{Verdiani} and \cite{GWZ}, simply connected cohomogeneity one manifolds admitting invariant metrics of positive sectional curvature were classified.
As it is very difficult to distinguish between manifolds admitting positive curvature and those that merely admit nonnegative curvature, it is interesting to see which cohomogeneity one manifolds can admit invariant metrics of nonnegative curvature.

One particularly interesting example in this context is the Brieskorn variety $B^{2n-1}_d$ with the cohomogeneity one action by $\SO(n)\SO(2)$ (see Section \ref{brieskorn}). In \cite{GVWZ}, it was shown that $B^{2n-1}_d$ admits an invariant metric of nonnegative sectional curvature if and only if $n\le 3$ or $d\le 2$. So most of these actions do not admit invariant nonnegatively curved metrics.

On the other hand, a construction for building metrics of nonnegative sectional curvature on a large class of cohomogeneity one manifolds was described in \cite{GZ1}. They showed that every cohomogeneity one manifold with two nonprincipal orbits of codimension 2 admits an invariant metric of nonnegative sectional curvature. The following theorem relies heavily on that result.

\begin{main}\label{curvature}
Every nonreducible cohomogeneity one action of a compact connected group on a compact simply connected manifold of dimension 7 or less admits an invariant metric of nonnegative sectional curvature, except the Brieskorn variety, $B^{7}_d$ for $d\ge 3$, and possibly some of the members of the following family of actions
 $$S^3\times S^3 \,\, \supset \,\, \set{(e^{ip\theta},e^{iq\theta})}, \,\, \Delta S^3 \cdot \Z_n \,\, \supset \,\,  \Z_n$$
where $(p,q)=1$ and either $n=1$; or $p$ or $q$ even and $n=2$.
\end{main}

\begin{rem*}
In the case $n=2$ and $q=p+1$, these actions are isometric actions on certain positively curved Eschenburg spaces (\cite{Ziller_notes} or \cite{GWZ}). So in fact, many of the members of this exceptional family are already known to admit invariant metrics of \emph{positive} sectional curvature. It is then reasonable to expect many more of them to admit nonnegative curvature as well.
\end{rem*}

Determining the full topology of all the manifolds appearing in the classification would be quite difficult. However, in dimension 5 we can give a complete answer.

\begin{main}\label{5class}
Every compact simply connected cohomogeneity one manifold of dimension 5 must be diffeomorphic to one of the following: $S^5$, $\SU(3)/\SO(3)$ or one of the two $S^3$ bundles over $S^2$.

In particular, the actions of type $P^5$ are all actions on $\SU(3)/\SO(3)$, and actions of type $N^5$ are either on $S^3\times S^2$ or the nontrivial $S^3$ bundle over $S^2$, depending on the parameters.
\end{main}

The paper is organized as follows. In Section \ref{prelims} we will discuss cohomogeneity one manifolds in general and develop some basic facts that will be useful throughout the paper. The classification will take place in Sections \ref{5dim} to \ref{7dim}. Next, in Sections \ref{identify} through \ref{5dimtop} we look at some of the actions in more detail and prove the main theorems. Finally, we give several convenient tables in the appendix.


This work was completed as a Ph.D. thesis while the author was a graduate student at the University of Pennsylvania. The author would like to extend his deepest thanks to W.~Ziller for his tremendous help throughout the entire process.

%
%
%
%
%
%

\renewcommand{\thetable}{\theequation}

\section{Cohomogeneity one manifolds}\label{prelims}

In this section we will discuss the cohomogeneity one action of a Lie group $G$ on a manifold $M$ in general. We first review the basic structure of such actions, and see that they are completely determined by certain isotropy subgroups. We will then discuss how we can determine the fundamental group of the manifold from these isotropy groups. We will also give some helpful restrictions on the possible groups that can occur in our situation, and we will review some basic facts about Lie groups that will be important for us. We end this section with a review of the classification of cohomogeneity one manifolds in dimensions 4 and lower in the compact simply connected case.

Throughout this section, \emph{$G$ will denote a compact connected Lie group and $M$ will be a closed and connected manifold}, unless explicitly stated otherwise.

\subsection{Basic Structure}\label{basic_structure}

The action of a compact Lie group on a manifold is said to be cohomogeneity one if the orbit space is one dimensional, or equivalently if there are orbits of codimension 1. Let $M$ be a compact, connected cohomogeneity one manifold for the compact, connected group $G$. It follows from the general theory of group actions that there is an open dense connected subset $M_0\subset M$ consisting of codimension 1 orbits and that the map $M_0 \to M_0/G$ is a locally trivial fiber bundle. Since $M_0/G$ is then one dimensional it must either be an open interval or a circle. In the case that $M_0/G \approx S^1$ it follows that $M = M_0$ and so $M$, being fibred over a circle, must have infinite fundamental group. Since we are only interested in simply connected manifolds we will henceforth restrict our attention to those $M$ with $M_0 / G \approx (-1,1)$. Since $M$ is compact it follows that $M/G \approx [-1,1]$. We will refer to such manifolds as {\em interval cohomogeneity one manifolds}.

To understand the well known structure of $M$ further, choose an arbitrary but fixed $G$-invariant Riemannian metric on $M$, and let $\pi: M \to M/G \approx [-1,1]$ denote the projection. Let $c:[-1,1]\to M$ be a minimal geodesic between the two non-principal orbits $\pi^{-1}(-1)$ and $\pi^{-1}(1)$ and reparameterize the quotient interval $M/G\approx [-1,1]$ such that $\pi\circ c = \id_{[-1,1]}$. Denote the isotropy groups by $H=G_{c(0)}$ and $K^\pm=G_{c(\pm1)}$ and let $D_{\pm}$ denote the disk of radius $1$ normal to the non-principal orbits $\pi^{-1}(\pm1)=G\cdot c(\pm1)$ at $c(\pm1)$. One can see then that $K^\pm$ acts on $D_\pm$ and does so transitively on $\partial D_\pm$ with isotropy $H$ at $c(0)$. Therefore $S^{l_\pm}:=\del D_\pm = K^\pm\cdot c(0) \approx K^\pm/H$. The slice theorem \cite{Br} tells us that the tubular neighborhoods of the non-principal orbits have the form $\pi^{-1}[-1,0] \approx G\times_{K^-}D_-$ and $\pi^{-1}[0,1] \approx G\times_{K^+}D_+$. Therefore we have decomposed our manifold into two disk bundles $G\times_{K^\pm}D_\pm$ glued along their common boundary $\pi^{-1}(0)=G\cdot c(0) \approx G/H$. That is
\begin{equation}\label{decomp}
  M\approx G\times_{K^-}D_- \,\, \cup_{G/H} \,\, G\times_{K^+}D_+ \text{ \hspace{1em} where \hspace{1em} } S^{l_\pm}=\del D_\pm \approx K^\pm/H.
\end{equation}
This gives a description of $M$ entirely in terms of $G$ and the isotropy groups $K^\pm \supset H$. The collection of $G$ with its isotropy groups $G\supset K^+,K^-\supset H$ is called the \emph{group diagram} of the cohomogeneity one manifold. Note: in the group diagram we understand that $G$ contains \emph{both} subgroups $K^-$ and $K^+$ and that \emph{both} $K^-$ and $K^+$ contain $H$ as a subgroup.

Conversely, let $G\supset K^+,K^- \supset H$ be compact groups with $K^{\pm}/H \approx S^{l_\pm}$. We know from the classification of transitive actions on spheres (see \cite{Besse-geos} page 195) that the $K^\pm$ action on $S^{l_\pm}$ must be linear and hence it extends to an action on the disk $D_\pm$ bounded by $S^{l_\pm}$, for each $\pm$. Therefore one can construct a cohomogeneity one manifold $M$ using \eqref{decomp}. \emph{So a cohomogeneity one manifold $M$ with $M/G \approx [-1,1]$ determines a group diagram $G\supset K^+,K^- \supset H$ with $K^{\pm}/H \approx S^{l_\pm}$ and conversely, such a group diagram determines a cohomogeneity one action}. This reduces the question of classifying such cohomogeneity one manifolds to a question of finding subgroups of compact groups with certain properties.

We sometimes record group diagrams as
\begin{equation}
 \label{diagram}
\begin{split}
\xymatrix{
& G \ar@{-}[dr] \ar@{-}[dl] & \\
K^{-}\ar@{-}[dr] & & K^+ \ar@{-}[dl] \\
{}\save[]+<-45pt,15pt>* {S^{l_{-}} = K^- /H }\restore & H &
{}\save[]+<45pt,15pt>* {S^{l_{+}} = K^+ /H}
\restore}
\end{split}
\end{equation}

The question of whether or not two different group diagrams determine the same action will be important to understand. We say the action of $G_1$ on $M_1$ is \emph{equivalent} to the action of $G_2$ on $M_2$ if there is a diffeomorphism $f:M_1\to M_2$ and an isomorphism $\phi:G_1\to G_2$ such that $f(g\cdot x)=\phi(g)\cdot f(x)$ for all $x\in M_1$ and $g\in G_1$. We will classify cohomogeneity one action up to this type of equivalence. However when $G_1=G_2$ there is a stronger type of equivalence that is sometimes preferred. We also say a map $f:M_1\to M_2$ between $G$-manifolds is $G$-\emph{equivariant} if $f(g\cdot x)=g\cdot f(x)$, for all $x\in M_1$ and $g\in G$. The next proposition, taken from \cite{GWZ}, applies to $G$-equivariant diffeomorphism.

\begin{prop}\label{normalizer}
Let a cohomogeneity one action of $G$ on $M$ be given by the group diagram $G\supset K^-,K^+\supset H$. Then any of the following operations on the group diagram will result in a $G$-equivariantly diffeomorphic manifold.
\begin{itemize}
\item[i.]Switching $K^-$ and $K^+$,
\item[ii.]Conjugating each group in the diagram by the same element of $G$,
\item[iii.]Replacing $K^-$ with $aK^-a^{-1}$ for $a\in N(H)_0$.
\end{itemize}

Conversely, the group diagrams for two $G$-equivariantly diffeomorphic manifolds must be taken to each other by some combination of these three operations.
\end{prop}

The next corollary is particularly helpful when finding the group diagram for a given cohomogeneity one action. Recall that the groups $K^-, K^+$ and $H$ are the isotropy groups along a minimal geodesic between nonprincipal orbits, with respect to some $G$ invariant metric on $M$. In most cases it is not convenient to explicitly find such a metric and geodesic. The following corollary solves this problem.

\begin{cor}\label{can_choose_c}
Let $M$ be an interval cohomogeneity one manifold for the group $G$ and let $\gamma:[-1,1]\to M$ be any continuous curve between nonprincipal orbits which meets each orbit precisely once and which is differentiable at the nonprincipal orbits with derivative transverse to these orbits. If $\gamma$ satisfies the further property $G_{\gamma(t)}=G_{\gamma(0)}$, for all $t\in (-1,1)$, then $G\supset G_{\gamma(-1)}, G_{\gamma(1)}\supset G_{\gamma(0)}$ is a valid group diagram for the action of $G$ on $M$.
\end{cor}

\begin{proof}
Fix a $G$ invariant metric on $M$ and let $c:[-1,1]\to M$ be a minimal geodesic between nonprincipal orbits. Then the group diagram along $c$ is given by $K^\pm=G_{c(\pm1)}$ and $H=G_{c(0)}$. We can assume that $\gamma(t)$ and $c(t)$ are in the same $G$-orbit for each $t\in[-1,1]$, after reparameterizing $\gamma$. This reparametrization will not effect the property of the derivative of $\gamma$ at the nonprincipal orbits. After applying some element of $G$ to $c$ we can also assume that $c(0)=\gamma(0)$.

Since the principal part of $M$ is $G$ equivariantly diffeomorphic to $G/H\times (-1,1)$, it follows that we can write $\gamma(t)=g(t)c(t)$ for some continuous curve $g:(-1,1)\to G$, with $g(0)=e$. The fact the $\gamma'(\pm1)$ exists and is transverse to the nonprincipal orbits means that we can extend $g(t)$ to a continuous function on $[-1,1]$ with $\gamma(t)=g(t)c(t)$.

Notice that we have $G_{\gamma(t)}=G_{g(t)c(t)}=g(t)G_{c(t)}g(t)^{-1}$ for all $t$. Since $\gamma(0)=c(0)$ and by our hypothesis we also know $G_{\gamma(t)}=G_{\gamma(0)}=H$ for all $t\in(-1,1)$. Therefore $H=G_{\gamma(t)} = g(t)G_{c(t)}g(t)^{-1} = g(t)Hg(t)^{-1}$ for all $t\in(-1,1)$. Since $g(0)=e$ it follows that $g(t)\in N(H)_0$ for all $t\in(-1,1)$. By continuity it follows that $n_\pm:=g(\pm1)\in N(H)_0$, as well. Putting this together we find that the diagram $G\supset G_{\gamma(-1)}, G_{\gamma(1)}\supset G_{\gamma(0)}$ is the diagram $G\supset n_-K^-n_-^{-1}, n_+K^+n_+^{-1}\supset H$, which represents our original action by Proposition \ref{normalizer}.
\end{proof}

Recall an action of $G$ on $M$ is \emph{effective} if no element $g\in G$ fixes $M$ pointwise, except $g=1$. We claim that a cohomogeneity one action, as above, is effective if and only if $G$ and $H$ do not share any nontrivial normal subgroups. It is clear that if $N$ is the ineffective kernel of the $G$ action, i.e. $N=\ker(G\to \Diff(M))$, then $N$ will be a normal subgroup of both $G$ and $H$. Conversely, let $N$ be the largest normal subgroup shared by $G$ and $H$. Then, as above, $N$ fixes the entire geodesic $c$ pointwise. Therefore, since $N$ is normal, it fixes all of $M$ pointwise. So it is easy to determine the effective version of any cohomogeneity one action from its group diagram alone. Because of this, we will generally allow our actions to be ineffective, however we will be most interested in the \emph{almost effective} actions, that is, actions with at most a finite ineffective kernel. In this case, $N$ is a discrete normal subgroup and hence $N\subset Z(G)$, where $Z(G)$ is the center of $G$. Then in fact $N=H\cap Z(G)$ in this case, by what we said above.

We now introduce the important idea of primitivity and nonprimitivity.

\begin{defn}\label{def:prim}
We say the cohomogeneity one manifold $M_G$ is \emph{nonprimitive} if it has some group diagram representation $G\supset K^-,K^+\supset H$ for which there is a proper connected closed subgroup $L\subset G$ with $L\supset K^-,K^+$. It then follows that $L\supset K^-,K^+\supset H$ is a group diagram which determines some cohomogeneity one manifold $M_L$.
\end{defn}

\begin{example*}
As an example, consider the group diagram $S^3\supset \set{e^{i\theta}}, \set{e^{j\theta}} \supset \set{\pm 1}$. There is no proper subgroup $L$ which contains both $K^-=\set{e^{i\theta}}$ and $K^+= \set{e^{j\theta}}$, however, by \ref{normalizer} this action is equivalent to the action with group diagram $S^3\supset \set{e^{i\theta}}, \set{e^{i\theta}} \supset \set{\pm 1}$. So in fact, this action is primitive.
\end{example*}

This next proposition shows the importance of the idea of primitivity.

\begin{prop}\label{prop:prim}
Take a nonprimitive cohomogeneity one manifold $M_G$ with $L$ and $M_L$ as in \ref{def:prim}. Then $M_G$ is $G$-equivariantly diffeomorphic to $G\times_L M_L = (G\times M_L)/L$ where $L$ acts on $G\times M_L$ by $\ell \star (g,x) = (g\ell^{-1}, \ell x)$. Hence there is a fiber bundle
  $$M_L\to M_G\to G/L.$$
\end{prop}
\begin{proof}
Let $c$ be a minimal geodesic in $M_L$ between nonprincipal orbits. Then it is clear that the curve $\tilde c(t)= (1,c(t))\in G\times M_L$ is a geodesic where we equip $G\times M_L$ with the product metric, for the biinvariant metric on $G$. It is also clear that $\tilde c$ is perpendicular to the $L$ orbits in $G\times M_L$. Therefore $c$ descends to a geodesic $\hat c$ in $G\times_L M_L$ which is perpendicular to the $G$ orbits. The isotropy groups of the $G$ actions on $G\times_L M_L$ are clearly given by $G_{\hat c(t)}= L_{c(t)}$ and hence this $G$ action on $G\times_L M_L$ is cohomogeneity one with group diagram $G\supset K^-,K^+\supset H$. This proves the proposition.
\end{proof}

\subsection{The Fundamental Group}\label{fundamental_group}

We will generally be looking at cohomogeneity one actions in terms of their group diagrams. Since this paper is concerned with simply connected cohomogeneity one manifolds, it will be important to be able to determine the fundamental group of the manifold using only the group diagram. In this section we will show how to do this and give strong but simple conditions on which group diagrams can give simply connected manifolds. Recall we are assuming that \emph{$G$ is compact and connected} throughout this section.

This first proposition is take from \cite{GWZ} (Lemma 1.6).

\begin{prop}\label{no_exceptionals}
Let $M$ be a compact simply connected cohomogeneity one manifold for the group $G$ as above. Then $M$ has no exceptional orbits, and hence, in the notation above, $l_\pm\ge 1$, or equivalently $\dim K^\pm > \dim H$.
\end{prop}

This next proposition can be considered as the van Kampen theorem for cohomogeneity one manifolds. It tells us precisely how to compute the fundamental groups from the group diagrams alone.

\begin{prop}[van Kampen]\label{fundgrp}
Let $M$ be the cohomogeneity one manifold given by the group diagram $G\supset K^+,K^- \supset H$ with $K^{\pm}/H = S^{l_\pm}$ and assume $l_\pm\ge 1$. Then $\pi_1(M)\approx \pi_1(G/H)/N_-N_+$ where
  $$N_\pm=\ker \{\pi_1(G/H)\to\pi_1(G/K^\pm)\}=\Im\{ \pi_1(K^\pm/H)\to \pi_1(G/H) \}.$$
In particular $M$ is simply connected if and only if the images of $K^{\pm}/H=S^{l_\pm}$ generate $\pi_1(G/H)$ under the natural inclusions.
\end{prop}
\begin{proof}
We will compute the fundamental group of $M$ using van Kampen's theorem. With the notation of Section \ref{basic_structure}, we know that $M$ can be decomposed as $\pi^{-1}([-1,0])\cup\pi^{-1}([0,1])$ where $\pi^{-1}([-1,0])\cap\pi^{-1}([0,1]) = G\cdot x_0$. Here, we know that, with a slight abuse of notation, $\pi^{-1}([0,\pm1])$ deformation retracts to $\pi^{-1}(\pm1) = G\cdot x_\pm \approx G/K^\pm$. So in fact we have the homotopy equivalence $\pi^{-1}([0,\pm1]) \to G/K^\pm : g\cdot c(t) \mapsto gK^\pm$. Therefore we have the commutative diagram of pairs:
\begin{equation}\label{orbitcomdiagram}
\xymatrix{
{(G\cdot x_0, x_0)}\ar[r]\ar[d]_\wr&{(\pi^{-1}([0,\pm1]), x_0)}\ar[d]_\wr\\
{(G/H,H)}\ar[r]&{(G/K^\pm, K^\pm)}
}
\end{equation}
where the vertical maps are both homotopy equivalences, the top map is the inclusion, and the bottom map is the natural quotient. This gives the corresponding diagram of fundamental groups:
$$\xymatrix{
{\pi_1(G\cdot x_0, x_0)}\ar[r]\ar[d]_\wr&{\pi_1(\pi^{-1}([0,\pm1]), x_0)}\ar[d]_\wr\\
{\pi_1(G/H,H)}\ar[r]&{\pi_1(G/K^\pm, K^\pm)}
}
$$
Therefore we may freely use $\pi_1(G/H)\to\pi_1(G/K^\pm)$ in place of $\pi_1(G\cdot x_0, x_0) \to \pi_1(\pi^{-1}([0,\pm1]), x_0)$ for van Kampen's theorem.

Now look at the fiber bundle
\begin{equation}\label{KH_fib_bund}
K^\pm/H\to G/H \to G/K^\pm \text{ \hspace{1em} where \hspace{1em} } K^\pm/H\approx S^{l_\pm}.
\end{equation}
This gives the long exact sequence of homotopy groups:
\begin{equation}\label{les_KH}
\begin{CD}
\cdots\to \pi_i(S^{l_\pm})@>i_*^\pm>> \pi_i(G/H) @>\rho_*^\pm>>\pi_i(G/K^\pm)@>\partial_*>> \pi_{i-1}(S^{l_\pm})\to\cdots\\
\cdots\to \pi_1(S^{l_\pm})@>i_*^\pm>> \pi_1(G/H) @>\rho_*^\pm>>\pi_1(G/K^\pm)@>\partial_*>> \pi_{0}(S^{l_\pm})
\end{CD}
\end{equation}
Notice that this implies $\rho_*^\pm:\pi_1(G/H) \to \pi_1(G/K^\pm)$ is onto, since $l_\pm>0$. In fact it follows $G/H\to G/K^\pm$ is $l_\pm$-connected, but we will not need this.

Then, by van Kampen's theorem, $\pi_1(M)\approx \pi_1(G/H)/N_-N_+$ where $N_\pm = \ker(\rho_*^\pm)$. Finally, by \ref{les_KH}, we see $N_\pm=\ker(\rho_*^\pm)=\Im(i_*^\pm)$, and this concludes the proof.
\end{proof}

With this language and notation we now give a reformulation of Lemma 1.6 from \cite{GWZ}. This corollary will be very convenient for dealing with the case that $l_-$ or $l_+$ is greater than 1.

\begin{cor}\label{fundgrp2}
Let $M$ be the cohomogeneity one manifold given by the group diagram $G\supset K^+,K^- \supset H$ with $K^{\pm}/H = S^{l_\pm}$.
\begin{enumerate}
\item[i.]If $l_+>1$ and $l_-\ge1$ then $\pi_1(M)\approx \pi_1(G/K^-)$ and $H\cap K^+_0=H_0$. In particular, if $M$ is simply connected, $K^-$ is connected.
\item[ii.]If $l_-,l_+>1$ then $\pi_1(M)\approx \pi_1(G/H)\approx \pi_1(G/K^\pm)$. In particular, if $M$ is simply connected,  all of $H$, $K^-$ and $K^+$ are connected.
\end{enumerate}
\end{cor}
\begin{proof}
First say $l_+>1$ and $l_-\ge 1$. Then $N_+$ is trivial and $\pi_1(M)\approx \pi_1(G/H)/N_-$ where $N_-=\ker \{\rho_*^-:\pi_1(G/H)\to\pi_1(G/K^-)\}$. Then, since $\rho_*^-$ is onto in the case $l_-\ge 1$, $\pi_1(M)\approx \pi_1(G/H)/\ker\rho_*^-\approx \Im\rho_*^-=\pi_1(G/K^-)$. If, in addition, $l_-> 1$ then $N_-$ is also trivial. Then $\pi_1(M)\approx \pi_1(G/H)$ and by the same argument $\pi_1(M)\approx \pi_1(G/K^+)$.

The connectedness of the groups involved follows from the following fact: If $L\subset J$ are compact lie groups, with $J$ connected, and $J/L$ simply connected then $L$ must be connected. Otherwise, $J/L_0\to J/L$ would be a non-trivial cover. This also implies that $K^+_0\cap H=H_0$ since $K^+_0/(H\cap K^+_0)$ is a simply connected sphere since $l_+>1$.
\end{proof}

This corollary tells us how to deal with the case that $l_-$ or $l_+$ is greater than one. For the case that both $l_\pm=1$ the following lemma will be very helpful.

\begin{lem}\label{Hgens}
Let $M$ be the cohomogeneity one manifold given by the group diagram $G\supset K^+,K^- \supset H$. Denote $H_\pm=H\cap K^\pm_0$ and let $\alpha_\pm:[0,1]\to K^\pm_0$ be curves which generate $\pi_1(K^\pm/H)$, with $\alpha_\pm(0)=1\in G$. $M$ is simply connected if and only if
\begin{enumerate}
\item[i.]$H$ is generated as a subgroup by $H_-$ and $H_+$, \emph{and}
\item[ii.]$\alpha_-$ and $\alpha_+$ generate $\pi_1(G/H_0)$.
\end{enumerate}
\end{lem}
\begin{rem*}
The curves $\alpha_\pm$ are, in general, not loops in $G/H_0$. However we can compose them in $G/H_0$ either via pointwise multiplication in $G$ or via lifting their compositions in $G/H$, where they are loops, to $G/H_0$. When we say $\alpha_\pm$ generate $\pi_1(G/H_0)$ we mean the combinations of these curves which form loops in $G/H_0$ generate $\pi_1(G/H_0)$. Also notice that if $\dim(K^-/H) > 1$ we can simply take $\alpha_-(t)= 1$, $\forall t$, and similarly for $\alpha_+$.
\end{rem*}

The fact that $(\pi_1(M)=0) \Rightarrow \text{(i.)}$ is equivalent to Lemma 1.7 of \cite{GWZ}. However we give an independent proof here which leads to the full version of this lemma.

\begin{proof}By \ref{fundgrp}, $M$ is simply connected if and only if $\pi_1(G/H)=\langle \alpha_-\rangle\cdot\langle \alpha_+\rangle$, where $\alpha_\pm$ are considered as loops in $G/H$. Furthermore, since $\alpha_\pm$ are loops in $K^\pm/H$ it follows that $a_\pm := \alpha_\pm(1) \in K^\pm_0\cap H=H_\pm$. It is clear that the group generated by $H_-$ and $H_+$ is the same as the group generated by $a_-$, $a_+$ and $H_0$. Therefore condition (i.) from the lemma is equivalent to the statement that $H$ is generated by $a_-$, $a_+$ and $H_0$.

First assume that $M$ is simply connected so that $\pi_1(G/H)=\langle \alpha_-\rangle\cdot\langle \alpha_+\rangle$. Since the map $G/H_0\to G/H$ is a cover, it is clear that $\alpha_-$ and $\alpha_+$ generate $\pi_1(G/H_0)$. We must only show that $H$ is generated by $a_-$, $a_+$ and $H_0$.

Choose an arbitrary component, $hH_0$, of $H$. We claim that some product of $a_-$ and $a_+$ will lie in $hH_0$. For this, let $\gamma:[0,1]\to G$ be an arbitrary path with $\gamma(0)=1$ and $\gamma(1)\in hH_0$. Then $\gamma$ represents a loop in $G/H$ and since $\pi_1(G/H)=\langle \alpha_-\rangle\cdot\langle \alpha_+\rangle$ we must have that $[\gamma]=[\alpha_-]^{n}\cdot[\alpha_+]^{m}$ for some $m,n\in\Z$, where $[\cdot]$ denotes the corresponding class in $\pi_1(G/H)$.

We now make use of the following observation. In general, for compact Lie groups $J\subset L$, take paths $\beta_\pm:[0,1]\to L$, with $\beta_\pm(0)=1$ and $\beta_\pm(1)\in J$. Then we see that $\left(\beta_-\cdot\beta_+(1)\right)\circ \beta_+$ is fixed endpoint homotopic to $\beta_-\cdot\beta_+$ in $L$, where $\beta_-\cdot\beta_+(1)$ is the path $t\mapsto \beta_-(t)\cdot\beta_+(1)$; $\circ$ denotes path composition; and $\beta_-\cdot\beta_+$ is the path $t\mapsto \beta_-(t)\cdot\beta_+(t)$. Therefore $[\beta_-]\cdot[\beta_+] = [\beta_-\cdot\beta_+]$ as classes in $\pi_1(J/L)$.

In our case, this implies
  $$[\gamma]=[\alpha_-]^{n}\cdot[\alpha_+]^{m} = [\alpha_-^{n}\cdot\alpha_+^{m}]$$
in $\pi_1(G/H)$. Now look at the cover $G/H_0\to G/H$. Since the paths $\gamma$ and $\alpha_-^{n}\cdot\alpha_+^{m}$ both start at $1\in G$, it follows that $\gamma$ and $\alpha_-^{n}\cdot\alpha_+^{m}$ both end in the same component of $H$. Hence
  $$(\alpha_-^{n}\cdot\alpha_+^{m})(1)  = \alpha_-(1)^{n}\cdot\alpha_+(1)^{m} = a_-^{n}\cdot a_+^{m} \in hH_0.$$
Therefore, $a_-$, $a_+$ and $H_0$ generate $H$, proving (i.).

Next suppose (i.) and (ii.) hold. Again, since the map $G/H_0\to G/H$ is a cover, $\pi_1(G/H)$ is generated by $\pi_1(G/H_0)$ and a collection of curves in $G/H_0$ which go from $H_0$ to each component of $H$. $\alpha_-$ and $\alpha_+$ already generate $\pi_1(G/H_0)$ by assumption. Saying $a_-$, $a_+$ and $H_0$ generate $H$ is equivalent to saying that combinations of $\alpha_-$ and $\alpha_+$ can reach any component of $H$, when considered as paths in $G/H_0$. Hence $\alpha_\pm$ generate $\pi_1(G/H)$ and $M$ is simply connected by \ref{fundgrp}.
\end{proof}

\subsection{Extensions and Reductions}\label{EaR}

In this section we will describe a natural way of reducing certain cohomogeneity one actions to actions by smaller groups with the same orbits. We will also describe a way of extending actions to larger groups and we will see that these two processes are inverses of each other.

\begin{prop}\label{reductions}
Let $M$ be the cohomogeneity one manifold given by the group diagram $G\supset K^+,K^- \supset H$ and suppose $G=G_1\times G_2$ with $\proj_2(H)=G_2$. Then the subaction of $G_1\times 1$ on $M$ is also by cohomogeneity one, with the same orbits, and with isotropy groups $K_1^\pm=K^\pm\cap (G_1\times 1)$ and $H_1=H\cap (G_1\times 1)$.
\end{prop}
\begin{proof}
Recall that the action of $G$ on each orbit $G\cdot x$ is equivalent to the $G$ action on $G/G_x$. So it is enough to test the claim on each type of orbit: $G/K^+$, $G/K^-$ and $G/H$. Let $G/G_x$ be one such orbit and notice that $H\subset G_x$. Then for each element $(g_1,g_2)G_x\in G/G_x$ there is some element of $H$ of the form $(h_1,g_2)$ since $\proj_2(H)=G_2$. Then $(g_1,g_2)G_x = (g_1h_1^{-1}, 1)\cdot (h_1,g_2) G_x = (g_1h_1^{-1}, 1)G_x$ and hence an arbitrary point $(g_1,g_2)G_x$ is in the $G_1\times 1$ orbit of $(1,1)G_x$. This proves $G_1\times 1$ acts on $M$ with the same orbits as $G$ and hence still acts by cohomogeneity one. The fact that the isotropy groups of the $G_1\times 1$ action are $K_1^\pm=K^\pm\cap (G_1\times 1)$ and $H_1=H\cap (G_1\times 1)$ is then clear.
\end{proof}

We will now describe a way of extending a given cohomogeneity one action to an action by a possibly larger group. Let $M$ be a cohomogeneity one manifold with group diagram $G_1\supset K^-_1, K^+_1\supset H_1$ and let $L$ be a compact connected subgroup of $N(H_1)\cap N(K^-_1)\cap N(K^+_1)$. Then notice $L\cap H_1$ is normal in $L$ and define $G_2:= L/(L\cap H_1)$. Then define an action by $G_1\times G_2$ on $M$ orbitwise by $(\hat g_1,[l])\star g_1 (G_1)_x= \hat g_1 g_1 l^{-1}(G_1)_x$ on each orbit $G_1/(G_1)_x$, for $(G_1)_x=H_1$ or $K^\pm_1$.

\begin{defn}\label{nomral_extensions}
Such an extension is called a \emph{normal extension}.
\end{defn}

\begin{prop}\label{p:normal_extension}
This extension describes a smooth action of $G:=G_1\times G_2$ on $M$ with the same orbits as $G_1$ and with group diagram
\begin{equation}\label{d:normal}
  G_1\times G_2 \,\, \supset \,\, (K^-_1\times 1)\cdot \Delta L, \,\, (K^+_1\times 1)\cdot \Delta L \,\, \supset \,\, (H_1\times 1)\cdot \Delta L
\end{equation}
where $\Delta L = \set{(l,[l])|\,l\in L}$.
\end{prop}
\begin{proof}
It is clear that this action is well defined and has the same orbits as the original $G_1$ action. Now let $c:[-1,1]\to M$ be a minimal geodesic between nonprincipal orbits in $M$ such that $(G_1)_{c(t)}=H_1$ for $t\in (-1,1)$ and $(G_1)_{c(\pm 1)}=K^\pm_1$. Then it is clear that the isotropy subgroups of $G=G_1\times G_2$ are
  $$H:=G_{c(t)}= H_1\cdot \Delta L \, \text{ for $t\in(-1,1)$ and } \, K^\pm:=G_{c(\pm1)}= K^\pm_1\cdot \Delta L$$
where we are identifying $G_1$ with $G_1\times 1$. So if we can show that the action is smooth and that there is a $G$-invariant metric on $M$ such that $c$ is a minimal geodesic then we will be done.

To do this let $\overline M$ be the manifold given by the group diagram $G \supset K^-, K^+\supset H$ with the corresponding geodesic $\overline c$, as above. Notice that $\proj_2(H)=\proj_2(H_1\cdot \Delta L)=G_2$. Hence by \ref{reductions}, $G_1$ still acts isometrically on $\overline M$ by cohomogeneity one with isotropy groups $(H_1\cdot \Delta L)\cap (G_1\times 1)=H_1\times 1$ and $(K^\pm_1\cdot \Delta L)\cap (G_1\times 1) = K^\pm_1\times 1$. Therefore $\overline M$ and $M$ are $G_1$-equivariantly diffeomorphic, via the map $\phi: g_1\cdot \wbar c(t) \mapsto g_1\cdot c(t)$.

We now claim that the map $\phi$ is also $G$-equivariant. To see this define the set theoretic map $\psi: \wbar M\to M: g\cdot \wbar c(t)\mapsto g\cdot c(t)$. This is well defined since $G$ has the same isotropy group at $\wbar c(t)$ as at $c(t)$. It is also clear that this set map is $G$-equivariant by definition. By restricting to elements of the form $g_1\cdot \wbar c(t)$ for $g_1\in G_1$ we see that $\psi(g_1\cdot \wbar c(t))= \phi(g_1\cdot \wbar c(t))$. Since the $G_1$ orbits are equal to the $G$ orbits in $\wbar M$ by \ref{reductions}, it follows that $\psi=\phi$ as maps. In particular $\psi$ is a diffeomorphism, since $\phi$ is. Therefore $\wbar M$ is $G$ equivariantly diffeomorphic to $M$. This completes the proof.
\end{proof}

\begin{prop}\label{extensions}
For $M$ as in Proposition \ref{reductions}, the action by $G=G_1\times G_2$ occurs as the normal extension of the reduced action of $G_1\times 1$ on $M$.
\end{prop}
\begin{proof}
We first claim that we can assume  $H\cap (1\times G_2)=1$, which will be useful for technical reasons. To see this, suppose $H_2:=H\cap (1\times G_2)$ is nontrivial. $H_2$ is obviously normal in $H$ and it is also normal in $G$ since $\proj_2(H)=G_2$. Then there is a more effective version of the same action by $(G_1\times G_2)/H_2\approx G_1 \times (G_2/H_2)=: G_1\times \tilde G_2$. We still have $\proj_2(\tilde H)=\tilde G_2$ for this action, where $\tilde H$ is the new principal isotropy group, and this time $\tilde H \cap (1\times \tilde G_2)=1$. So assume $H\cap (1\times G_2)=1$.

Consider the reduced action with diagram $G_1\times 1\supset K^-_1\times 1,K^+_1\times 1\supset H_1\times 1$ from \ref{reductions}. Let $L=\proj_1(H_0)\subset G_1$. We claim that the original $G_1\times G_2$ action is equivalent to the normal extension of the $G_1$ action via $L$. First notice that since $H_1$ is normal in $H$, $H_1$ is also normal in $L=\proj_1(H_0)$. Similarly $L$ is in the normalizer of $K^\pm_1$. So in fact
  $$L\subset N(H_1)\cap N(K^-_1)\cap N(K^+_1).$$

Now, notice the map $\proj_1:H_0\to L= \proj_1(H_0)$ is onto with trivial kernel, since we assumed $H\cap (1\times G_2)=1$. Therefore $\proj_1$ is a Lie group isomorphism and hence has an inverse $\psi: L\to H_0$ which must have the form $\psi(l)=(l,\phi(l))$ for some map $\phi:L\to G_2$. Notice that $\phi$ maps $L$ onto $G_2$ with kernel $H_1\cap L$. Therefore $G_2\approx L/(H_1\cap L)$, via $\phi$.

Notice that $H_0=\psi(L)=\set{(l,\phi(l))}$. It is also clear that $H=H_1\cdot H_0$ and similarly $K^\pm=K^\pm_1\cdot H_0$. Therefore we can write the group diagram for our original $G_1\times G_2$ action as
  $$G_1\times G_2\,\, \supset\,\,  K^-_1\cdot H_0, \,\, K^+_1\cdot H_0\,\,  \supset \,\, H_1\cdot H_0.$$
Then, after the isomorphism $G_1\times G_2 \to G_1\times \big(L/(H_1\cap L)\big): (g_1,\phi(l))\mapsto (g_1,[l])$, $H_0$ becomes $\Delta L:=\set{(l,[l])}$ and this diagram becomes exactly the diagram in \ref{d:normal}. Therefore the original action by $G_1\times G_2$ is equivalent to the normal extension of the $G_1$ action along $L$.
\end{proof}

\begin{defn}\label{irred}
We say the cohomogeneity one group diagram $G\supset K^-,K^+\supset H$ is \emph{nonreducible} if $H$ does not project onto any factor of $G$.
\end{defn}

If our group $G$ has the form $G_1\times \cdots \times G_l\times T^n$, where $G_i$ are all simple Lie groups, then we claim the condition that $H$ projects onto some factor of $G$ is equivalent to the condition that some proper normal subgroup of $G$ acts by cohomogeneity one with the same orbits. \ref{reductions} proves this claim in one direction. Conversely, suppose that some normal subgroup $N$ of $G$ acts by cohomogeneity one with the same orbits. Since the orbits of $G$ are connected we can assume that $N$ is connected. Therefore $N=\prod_{i\in I} G_i\times T^p$, for some subset $I\subset \set{1,\dots,l}$ and some $T^p\subset T^n$. Then let $L=\prod_{i\notin I} G_i\times T^q$, where $T^p\times T^q=T^n$, so that $G=N\times L$. The assumption that $N$ acts on $M$ with the same orbits means that $N$ acts transitively on $G/H=(N\times L)/H$. This means we can write any element $(n,\ell)H\in G/H$ as $(\tilde n, 1)H$ and hence $H$ must project onto $L$.

Since every compact Lie group has a cover of the form $G_1\times \cdots \times G_l\times T^n$, every cohomogeneity one action can be given almost effectively by such a group. So it makes sense to call an arbitrary cohomogeneity one action reducible if there is some proper normal subgroup that still acts by cohomogeneity one.

Most importantly, this section shows that the classification of cohomogeneity one manifolds is quickly reduced to the classification of the nonreducible ones. Therefore we will assume in our classification that all our actions are nonreducible and we will loose little generality, since every other cohomogeneity one action will be a normal extension of a nonreducible action.

\subsection{More limits on the groups}\label{more_limits}
In this section we give a few more restrictions on the groups that can act by cohomogeneity one on simply connected manifolds. The first addresses the case that the group has an abelian factor.

\begin{prop}\label{abfactor}
Let $M$ be the cohomogeneity one manifold given by the group diagram $G\supset K^+,K^- \supset H$ where $G=G_1\times T^m$ acts almost effectively and nonreducibly and $G_1$ is semisimple. Then we know $H_0=H_1\times 1\subset G_1\times 1$. Further, if $M$ is simply connected then $m\le 2$ and
\begin{enumerate}
  \item[i.]if $m=1$ then at least one of $\proj_2(K^\pm_0)=S^1$, say $\proj_2(K^-_0)=S^1$. Then $K^-/H \approx S^1$ and $K^-=S^1_-\cdot H$ for a circle group $S^1_-$, with $\proj_2(S^1_-)=S^1$. Furthermore, if $\rank(H)=\rank(G_1)$ or if $H_1$ is maximal-connected in $G_1$, then $H$, $K^-$ and $K^+$ are all connected; $K^-=H_1\times S^1$; and $K^+$ is either $H_1\times S^1$ or has the form $K_1\times 1$, for $K_1/H_1\approx S^{l_+}$.
  \item[ii.]if $m=2$ then both $K^\pm/H$ are circles and $K^\pm=S^1_\pm\cdot H$ for circle groups $S^1_\pm$, with $\proj_2(S^1_-)\cdot\proj_2(S^1_+)= T^2$. Furthermore, if $\rank(H)=\rank(G_1)$ then the $G$ action is equivalent to the product action of $G_1\times T^2$ on $(G_1/H_1)\times S^3$, where $T^2$ acts on $S^3\subset \C^2$ by component-wise multiplication.
\end{enumerate}
\end{prop}
\begin{proof}
Notice that in all cases $\proj_2(K^\pm_0)$ is a compact connected subgroup of $T^m$. Now say $\proj_2(K^-_0)$ is nontrivial. It must then be a torus, $T^n\subset T^m$. Then we have $\proj_2:K^-_0\twoheadrightarrow T^n$ with kernel $K^-_0\cap (G_1\times 1)$. Therefore we have the fiber bundle
  $$\big(K^-_0\cap (G_1\times 1)\big)/(H_1\times 1) \to K^-_0/(H_1\times 1) \to K^-_0/\big(K^-_0\cap (G_1\times 1)\big) \approx T^n$$
which gives the piece of the long exact sequence
  $$\pi_1(K^-_0/(H_1\times 1)) \to \pi_1(T^n) \to \pi_0(\big(K^-_0\cap (G_1\times 1)\big)/(H_1\times 1)).$$
The last group in this sequence is finite and the middle group is infinite. This means that $K^-_0/H_0$ has infinite fundamental group. Given that this space is a sphere, it follows that $K^-/H\approx S^1$. Therefore $K^-_0=H_0\cdot S^1_-$, for some circle group $S^1_-$ with $\proj_2(S^1_-)=S^1\subset T^m$. Similarly, if $\proj_2(K^+_0)$ is nontrivial then $K^+_0=H_0\cdot S^1_+$, for $S^1_+$ with $\proj_2(S^1_+)=S^1\subset T^m$.


We know that $\proj_2(K^-_0)$ and $\proj_2(K^+_0)$ generate some torus $T^n$ in $T^m$, with $n\le 2$. It is clear that if $m>n$ then $K^-/H$ and $K^+/H$ will not generate $\pi_1(G/H)$ and hence $M$ will not be simply connected, by \ref{fundgrp}. Therefore, $m\le2$, and if $m=1$ then one of $K^\pm$ must be a circle as above, and if $m=2$ then both $K^\pm$ must be circles as above. This proves the first part of the proposition.

For the second part, if $\rank H_1=\rank G_1$ or if $H_1$ is maximal-connected in $G_1$, we first claim that $\proj_1(K^-_0) = H_1$, if $K^-/H\approx S^1$. In the case that $H_1$ is maximal in $G_1$ this is clear since if $\proj_1(K^-_0)$ is larger than $H_1$ it would be all of $G_1$. Yet there is no compact semisimple group $G_1$ with subgroup $H_1$ where $G_1/H_1\approx S^1$. For the case that $\rank(H_1)=\rank(G_1)$, recall that for a general compact Lie group, the rank and the dimension have the same parity modulo 2. Since $K^-=S^1_-\cdot H$, $\proj_1(K^-_0)$ is at most one dimension larger than $H_1$. But if  $\proj_1(K^-_0)$ is of one higher dimension than $H_1$ it would follow that $\rank(\proj_1(K^-_0)) = \rank(H_1)+1 = \rank(G_1)+1$, a contradiction since $\proj_1(K^-_0)\subset G_1$. Therefore $\proj_1(K^-_0) = H_1$ in either case. Then since $K^- = S^1_-\cdot H$ it follows that $K^-_0 = H_1\times S^1_- \subset G_1\times T^m$. Similarly if $K^+/H\approx S^1$ then $K^+_0=H_1\times S^1_+$.

To see that all the groups are connected in this case, we notice that if $K^-_0\cap H$ is not $H_0$ then $H\cap 1\times S^1$ is nontrivial and there is a more effective action for the same groups with $H\cap 1\times S^1 = 1$. So we can assume that $K^-_0\cap H = H_0$. If, in addition, $K^+/H\approx S^1$ then by the same argument $K^+_0\cap H = H_0$ as well. If $\dim(K^+/H)>1$ then $K^+_0\cap H = H_0$ already, since $K^+_0/(H\cap K^+_0)$ would be a simply connected sphere. In any case we know that $K^\pm_0 \cap H = H_0$. Then, by \ref{Hgens}, $H$ must be connected, forcing $K^-$ and $K^+$ to be connected as well.

Now, it only remains to prove the last statement of (\emph{ii.}). In this case we already know $K^-=H_1\times S^1_-$ and $K^+=H_1\times S^1_+$. It is then clear that $K^-/H$ and $K^+/H$ generate $\pi_1(G/H)\approx \pi_1((G_1/H_1)\times T^2)$ if and only if $S^1_-$ and $S^1_+$ generate $\pi_1(T^2)$. This happens precisely when there is an automorphism of $T^2$ taking $S^1_-$ to $S^1\times 1$ and $S^1_+$ to $1\times S^1$. From \ref{fundgrp} we can assume this automorphism exists. After this automorphism the group diagram has the form
  $$G_1\times S^1\times S^1 \,\, \supset \,\, H_1 \times S^1\times 1, \,\, H_1\times 1\times S^1 \,\,\supset\,\, H_1\times 1\times 1.$$
It is easy to check that this action is the action described in the proposition (see Section \ref{special_types} for more details).
\end{proof}

The next two propositions give the possible dimensions that the group $G$ can have, if it acts by cohomogeneity one.

\begin{prop}\label{dimG}
If a Lie group $G$ acts almost effectively and by cohomogeneity one on the manifold $M^n$ then
  $$n-1 \le \dim(G) \le n(n-1)/2.$$
\end{prop}
\begin{proof}
Recall that for a principal orbit $G\cdot x\approx G/H$, $\dim G/H=n-1$, so the first inequality is trivial. Now we claim that $G$ also acts almost effectively on a principal orbit $G\cdot x\approx G/H$. To see this, suppose an element $g\in G$ fixes $G\cdot x$ pointwise. Then in particular $g\in H$. We saw above that $H$ fixes the geodesic $c$ pointwise and hence $g$ fixes all of $M$ pointwise. So, in fact, $G$ acts almost effectively on $G/H$. Now equip $G/H$ with a $G$ invariant metric. It then follows that $G$ maps into $\Isom G/H$ with finite kernel. Since $\dim G/H=n-1$, we know $\dim (\Isom G/H) \le n(n-1)/2$ and this proves the second inequality.
\end{proof}

The case where $G$ has the largest possible dimension is special. For this we have the following proposition.

\begin{prop}\label{topdim}
Let $G$ be a compact Lie group that acts almost effectively and by cohomogeneity one on the manifold $M^n$, $n>2$, with group diagram $G\supset K^-,K^+\supset H$ and no exceptional orbits. If $\dim(G) = n(n-1)/2$ and $G$ is simply connected then $G$ is isomorphic to $\Spin(n)$ and the action is equivalent to the $\Spin(n)$ action on $S^n\subset \R^n\times \R$ where $\Spin(n)$ acts on $\R^n$ via $\SO(n)$, leaving $\R$ pointwise fixed.
\end{prop}
\begin{proof}
Notice first that since $G$ acts on $M$ almost effectively we know that $G$ also acts almost effectively on the principal orbits which are equivariantly diffeomorphic to $G/H$. Now endow $G/H$ with the metric induced from a biinvariant metric on $G$, so that $G$ acts by isometry. Therefore we have a Lie group homomorphism $G \to \Isom G/H$ with finite kernel. Since $\dim G = n(n-1)/2$ and $\dim G/H = n-1$ it follows that $G/H$ must be a space form (see \cite{Petersen}). Further, since $G$ is simply connected it follows that $G/H_0$ is a compact simply connected space form. Hence $G/H_0$ is isometric to $S^{n-1}$ and $G$ still acts almost effectively and by isometry on $S^{n-1}$. So in fact $G \to \Isom S^{n-1} = \SO(n)$ as a Lie group homomorphism with finite kernel. Since $\dim G = \dim \SO(n)$ it follows that $G$ is isomorphic to $\Spin(n)$. We also know that the only way $\Spin(n)$ can act transitively on an $(n-1)$-sphere is with $\Spin(n-1)$ isotropy (see \cite{Besse-geos} page 195). Therefore there is an isomorphism $G\to\Spin(n)$ taking $H_0$ to $\Spin(n-1)$.

We also see that $\Spin(n-1)$ is maximal among connected subgroups of $\Spin(n)$. Hence $K^\pm$ must both be $\Spin(n)$ and hence $H$ is connected since $n>2$. Therefore the group diagram for this action is $\Spin(n)\supset \Spin(n),\Spin(n)\supset \Spin(n-1)$. It is easy to check that the $\Spin(n)$ action on $S^n$ described in the proposition also gives this diagram. Hence the two actions are equivalent.
\end{proof}

\subsection{Special types of actions}\label{special_types}

There are several types of actions that are easily described and easily recognized from their group diagrams. We will discuss these here so that we can exclude them in our classification. We summarize the results of this section in Table \ref{SaP} of the appendix.

\subsubsection{Product actions}
Say $G$ acts on $M$ by cohomogeneity one with group diagram $G\supset K^-,K^+\supset H$, and $L$ acts transitively on the homogeneous space $L/J$. Then it is clear that the action of $G\times L$ on $M\times (L/J)$ as a product, i.e.\ $(g,l)\star (p, \ell J)=(g p, l\ell J)$, is cohomogeneity one. Suppose $c$ is a minimal geodesic in $M$ between nonprincipal orbits, which gives the group diagram above. If we fix an $L$-invariant metric on $L/J$ then in the product metric on $M\times (L/J)$ the curve $\tilde c=(c,1)$ is a minimal geodesic between nonprincipal orbits. It is easy to see that the resulting group diagram is
\begin{equation}\label{product_diagram}
   G\times L \,\, \supset \,\, K^-\times J, \,\, K^+\times J \,\, \supset \,\, H\times J.
\end{equation}
Conversely, any diagram of this form will give a product action as described above. These diagrams are easy to recognize from the $J$ factor that appears in each of the isotropy groups.

\subsubsection{Sum actions}
Suppose $G_i$ acts transitively, linearly and isometrically on the sphere $S^{m_i}\subset \R^{m_i+1}$ with isotropy subgroup $H_i$, for $i=1,2$. Then we have an action of $G:=G_1\times G_2$ on $S^{m_1+m_2+1}\subset \R^{m_1+1}\times \R^{m_2+1}$ by taking the product action: $(g_1,g_2)\star (x,y)=(g_1\cdot x,g_2\cdot y)$. Such actions are called sum actions. Now, fix two unit vectors $e_i\in S^{m_i}$ with $(G_i)_{e_i}=H_i$, for $i=1,2$, and define $c(\theta)=(\cos(\theta)e_1,\sin(\theta)e_2)$. Upon computing the isotropy groups we find that the orbits through $c(\theta)$ for $\theta\in (0,\pi/2)$ are codimension one and hence this action is cohomogeneity one. We easily find the group diagram to be
\begin{equation}\label{sum_diagram}
   G_1\times G_2 \,\, \supset \,\, G_1\times H_2, \,\, H_1\times G_2 \,\, \supset \,\, H_1\times H_2.
\end{equation}
Conversely, take a group diagram of this form. Then $G_i/H_i$ are spheres and hence by the classification of transitive actions on spheres, $G_i$ actually acts linearly and isometrically on $S^{m_i}\subset \R^{m_i+1}$. Hence this action is a sum action as described above. Diagrams of the form \ref{sum_diagram} are easy to recognize from the $H_1$ and $H_2$ factors in the ``middle'' and the $G_1$ and $G_2$ factors on the ``outside'' of the pair $K^-,K^+$. In particular these actions are always isometric actions on symmetric spheres.

\subsubsection{Fixed point actions}
Here we will completely characterize the cohomogeneity one actions that have a fixed point. In fact we will not put any dimension restrictions on the actions in this subsection. Say $G$ acts effectively and by cohomogeneity one on the simply connected manifold $M$ and assume there is a fixed point $p_-\in M$, i.e.\ $G_{p_-}=G$. It is clear that the point $p_-$ cannot be in a principal orbit, so we can assume that $K^-=G$. Therefore the group diagram for this action will have the form
\begin{equation}\label{fixed_point_diagram}
   G \,\, \supset \,\, G, \,\, K^+ \,\, \supset \,\, H.
\end{equation}
Conversely, such a diagram clearly gives an action with a fixed point. So to classify fixed point cohomogeneity one actions we must only classify diagrams of type \ref{fixed_point_diagram}.

Because we assumed the action is effective, it follows that the $G$ action on $G/H\approx S^{l_-}$ is an effective transitive action on a sphere. Such actions were classified by Montgomery, Samelson and Borel (see \cite{Besse-geos} page 195). Up to equivalence, this gives us the possibilities for $G$ and $H$. In particular $H$ and hence $K^+$ must be connected. In Section 2 of \cite{GZ2}, the authors list all possible closed connected subgroups $K^+$ between $H$ and $G$, for each pair $G,H$.

In the case where $K^+=G$, we have
\begin{equation}\label{2fixed_points-diagram}
   G \,\, \supset \,\, G, \,\, G \,\, \supset \,\, H.
\end{equation}
To see what this action is, identify $G/H$ with the unit sphere $S^l\subset \R^{l+1}$. We know from the classification of transitive actions on spheres mentioned above, that $G$ acts linearly and isometrically on $\R^{l+1}$. It is easy to check that $M=S^{l+1}\subset \R^{l+1}\times \R$ with the action given by $g\star (x,t)= (gx,t)$. We will call such actions \emph{two-fixed-point} actions. In particular this is an isometric action on the sphere $S^{l+1}$. Notice that if $H_0$ is maximal among connected subgroups of $G$ then \ref{2fixed_points-diagram} is the only possible diagram for this $G$ and $H_0$, assuming there are no exceptional orbits. This gives the following convenient proposition.

\begin{prop}\label{Hmax}
Let $M$ be a simply connected cohomogeneity one manifold for the group $G$, with principal isotropy group $H$, as above. If $H_0$ is maximal among connected subgroups of $G$ then the action is equivalent to an isometric two-fixed-point action on a sphere.
\end{prop}

Therefore we must only consider the case in which $K^+$ is a subgroup strictly between $H$ and $G$. Following the tables given in \cite{GZ2}, we address these cases one by one. We first list the diagram, then the corresponding action. In each case it is easy to check that the action listed gives the corresponding diagram.

\begin{enumerate}
\item[$\cdot$] $\SU(n) \,\, \supset \,\, \SU(n), \,\, \S(\U(n-1)\U(1)) \,\, \supset \,\, \SU(n-1)$:

$\SU(n)$ on $\CP^n$ given by $A\star[z_0,z_1,\dots,z_n] = [z_0,A(z_1,\dots,z_n)]$.

\item[$\cdot$] $\U(n) \,\, \supset \,\, \U(n), \,\, \U(n-1)\U(1) \,\, \supset \,\, \U(n-1)$:

$\U(n)$ on $\CP^n$ given by $A\star[z_0,z_1,\dots,z_n] = [z_0,A(z_1,\dots,z_n)]$.

\item[$\cdot$] $\Sp(n) \,\, \supset \,\, \Sp(n), \,\, \Sp(n-1)\Sp(1) \,\, \supset \,\, \Sp(n-1)$:

$\Sp(n)$ on $\HP^n$ given by $A\star[x_0,x_1,\dots,x_n] = [x_0,A(x_1,\dots,x_n)]$.

\item[$\cdot$] $\Sp(n) \,\, \supset \,\, \Sp(n), \,\, \Sp(n-1)\U(1) \,\, \supset \,\, \Sp(n-1)$:

$\Sp(n)$ on $\CP^{2n+1}=S^{4n+3}/S^1$ for $S^{4n+3}\subset \H^{n+1}$ given by $A\star[x_0,x_1,\dots,x_n] = [x_0,A(x_1,\dots,x_n)]$

\item[$\cdot$] $\Sp(n)\times \Sp(1) \,\, \supset \,\, \Sp(n)\times \Sp(1), \,\, \Sp(n-1)\Sp(1)\times \Sp(1) \,\, \supset \,\, \Sp(n-1)\Delta \Sp(1)$:

$\Sp(n)\times \Sp(1)$ on $\HP^n$ given by $(A,p)\star[x_0,x_1,\dots,x_n] = [px_0,A(x_1,\dots,x_n)]$.

\item[$\cdot$] $\Sp(n)\times \U(1) \,\, \supset \,\, \Sp(n)\times \U(1), \,\, \Sp(n-1)\Sp(1)\times \U(1) \,\, \supset \,\, \Sp(n-1)\Delta \U(1)$:

$\Sp(n)\times \U(1)$ on $\HP^n$ given by $(A,z)\star[x_0,x_1,\dots,x_n] = [zx_0,A(x_1,\dots,x_n)]$.

\item[$\cdot$] $\Sp(n)\times \U(1) \,\, \supset \,\, \Sp(n)\times \U(1), \,\, \Sp(n-1)\U(1)\times \U(1) \,\, \supset \,\, \Sp(n-1)\Delta \U(1)$:

$\Sp(n)\times \U(1)$ on $\CP^{2n+1}=S^{4n+3}/S^1$ for $S^{4n+3}\subset \H^{n+1}$ given by \\
$(A,z)\star[x_0,x_1,\dots,x_n] = [zx_0,A(x_1,\dots,x_n)]$.                           

\item[$\cdot$] $\Spin(9) \,\, \supset \,\, \Spin(9), \,\, \Spin(8) \,\, \supset \,\, \Spin(7)$:

$\Spin(9)$ on $\CaP=\F_4/\Spin(9)$ (see \cite{Iwata}).
\end{enumerate}

In conclusion, we have shown the following.

\begin{prop}\label{fixedpt}
Every cohomogeneity one action on a compact simply connected manifold with a fixed point is equivalent to one of the isometric action on a compact rank one symmetric space described above.
\end{prop}

\subsection{Important Lie Groups}\label{lie_groups}

It is well known that every compact connected Lie group has a finite cover of the form $G_{ss}\times T^k$, where $G_{ss}$ is semisimple and simply connected and $T^k$ is a torus. The classification of simply connected semisimple Lie groups is also well know and all the possibilities are listed in Table \ref{posgrps} for dimension 21 and less. If an arbitrary compact group $G$ acts on a manifold $M$, then every cover $\tilde G$ of $G$ still acts on $M$, although less effectively. So allowing for a finite ineffective kernel, and because $G$ will always have dimension 21 or less by \ref{dimG}, we can assume that $G$ is the product of groups from Table \ref{posgrps}.

{\setlength{\tabcolsep}{0.40cm}
\renewcommand{\arraystretch}{1.6}
\stepcounter{equation}
\begin{table}[!h]
      \begin{center}
          \begin{tabular}{|c c c|}
\hline
Group    & Dimension    & Rank   \\
\hline \hline
$S^1\approx \U(1)\approx \SO(2)$   &   1  & 1 \\
\hline
$S^3\approx \SU(2)\approx \Sp(1)\approx \Spin(3)$   &   3  & 1    \\
\hline
$\SU(3)$   &   8  & 2   \\
\hline
$\Sp(2)\approx \Spin(5)$   &   10  & 2   \\
\hline
$\G_2$   &   14  & 2    \\
\hline
$\SU(4)\approx \Spin(6)$   &   15  & 3    \\
\hline
$\Sp(3)$   &   21  & 3   \\
\hline
$\Spin(7)$   &   21  & 3   \\
\hline
          \end{tabular}
      \end{center}
      \vspace{0.1cm}
      \caption{Classical compact groups in dimensions 21 and less, up to cover.}\label{posgrps}
\end{table}}

For the classifications of cohomogeneity one diagrams we will also need to know the subgroups of the groups listed in Table \ref{posgrps}, for certain dimensions. These subgroups are well known (see for example \cite{Dynkin}). Table \ref{grplist} lists closed connected subgroups of the indicated groups in the dimensions that will be relevant for our study.

{\setlength{\tabcolsep}{0.40cm}
\renewcommand{\arraystretch}{1.6}
\stepcounter{equation}
\begin{table}[!h]
      \begin{center}
          \begin{tabular}{|c||c c|}
\hline
Group    & Dimensions   & Subgroups   \\
\hline \hline
$T^2$   &   $\dim\ge 1$:  & $\set{(e^{ip\theta},e^{iq\theta})}$ \\
\hline
$S^3$   &  $\dim\ge 1$:  & $\set{e^{x\theta}=\cos\theta+x\sin\theta}$ for $x\in\Im(S^3)$   \\
\hline
$S^3\times S^3$   &  $\dim\ge 1$:  & $S^1 \subset T^2$; $T^2$; \\
                  &    &$S^3\times 1$; $1\times S^3$; $\Delta S^3=\set{(g,g)}$;\\
                  &    & $S^3\times S^1$; $S^1\times S^3$\\
\hline
$\SU(3)$   &   $\dim\ge 1$:  & $S^1\subset T^2$; $T^2$; $\SO(3)$; \\
           &        & $\SU(2)$; $\U(2)=\S(\U(2)\U(1))$  \\
\hline
$\Sp(2)$   &   $\dim\ge 4$: & $\U(2)$; $\Sp(1)\SO(2)$; $\Sp(1)\Sp(1)$\\
\hline
$\G_2$   &   $\dim\ge 8$: & $\SU(3)$   \\
\hline
$\SU(4)$   &  $\dim\ge 9$: & $\U(3)$; $\Sp(2)$   \\
\hline
          \end{tabular}
      \end{center}
      \vspace{0.1cm}
      \caption{Compact connected proper subgroups in the specified dimensions, up to conjugation.}\label{grplist}
\end{table}}

We can use this information about the subgroups of the classical Lie groups to make the following claim.

\begin{prop}\label{prods}
Let $M$ be the cohomogeneity one manifold given by the group diagram $G\supset K^+,K^- \supset H$, where $G$ acts nonreducibly on $M$. Suppose $G$ is the product of groups
  $$G = \prod_{t=1}^i(\SU(4)) \times \prod_{t=1}^j(\G_2) \times \prod_{t=1}^k(\Sp(2)) \times \prod_{t=1}^l(\SU(3)) \times \prod_{t=1}^m(S^3) \times (S^1)^n$$
where $i,j,k,l,m$ and $n$ are allowed to be zero and where we imagine most of them are zero. Then
  $$\dim(H)\le 10i + 8j + 6k + 4l + m.$$
\end{prop}
Of course the most important applications of this proposition will be in the case that $i,j,k,l,m$ and $n$ are all small and mostly zero. Although this might not seem helpful, we will see that it is very helpful in ruling out many product groups.
\begin{proof} Since the action is nonreducible, we know that $H$ does not project onto any of the factors in this product. That means each $\proj_\nu(H)$ is a proper subgroup. Therefore we have
$$
  \nonumber H_0  \subset  \prod_{t=1}^i(I_t) \times \prod_{t=1}^j(J_t) \times \prod_{t=1}^k(K_t) \times \prod_{t=1}^l(L_t) \times \prod_{t=1}^m(S^1_t) \times 1.
$$
Table \ref{grplist} gives us the largest possible dimension of each of these subgroups. In particular $\dim(I_t)\le 10$, $\dim(J_t)\le 8$, $\dim(K_t)\le 6$ and $\dim(L_t)\le 4$.
\end{proof}


\subsection{Low dimensional classification}\label{2-4dim}

Recall that cohomogeneity one manifolds of dimension 4 and lower were classified in \cite{Ne} and \cite{Parker}. For the convenience of the reader, and in order to correct an omission in \cite{Parker}, we will now reproduce this classification in the case of compact simply connected cohomogeneity one manifolds.

Suppose the compact connected group $G$ acts almost effectively and nonreducibly on the compact simply connected manifold $M$ by cohomogeneity one with group diagram $G\supset K^-, K^+\supset H$.

\subsubsection{2 Dimensional manifolds}

Suppose $\dim(M)=2$. We know from \ref{dimG} that $\dim(G)=1$, and hence $G=S^1$. Then, for the action to be effective, $H$ must be trivial. By \ref{no_exceptionals}, $K^\pm$ must both be $S^1$ and hence the group diagram is $S^1\supset S^1, S^1\supset 1$. It is easy to see that this is the action of $S^1$ on $S^2$ via rotation about some fixed axis.

\subsubsection{3 Dimensional manifolds}

Now suppose $\dim M =3$ so that $\dim G$ is either 2 or 3. If $\dim G=2$ then $G=T^2$ and $H$ is discrete. For the action to be effective, $H$ must be trivial. Then it is clear that both $K^\pm$ are circle groups in $T^2$. From \ref{fundgrp}, $M$ will be simply connected if and only if $K^-$ and $K^+$ generate $\pi_1(T^2)$. This happens precisely when there is an automorphism of $T^2$ taking $K^-$ to $S^1\times 1$ and $K^+$ to $1\times S^1$. So our group diagram in this case is $T^2\supset S^1\times 1, 1\times S^1\supset 1$, up to automorphism. It is easy to check that this is the action of $T^2$ on $S^3\subset \C^2$ via $(z,w)\star(x,y)=(zx,wy)$.

Next, if $\dim G=3$ then, by \ref{abfactor}, $G=S^3$ up to cover. Then \ref{topdim} says this action is a two-fixed-point action on a sphere.

\subsubsection{4 Dimensional manifolds}\label{dim4class}

Finally say $\dim M=4$, so that $3\le \dim G\le 6$. Up to cover, \ref{abfactor} says that $G$ must be one of the following four groups: $S^3$, $S^3\times S^1$, $S^3\times T^2$ or $S^3\times S^3$. If $G=S^3\times S^3$ then \ref{topdim} again implies that the action is a two-fixed-point action on $S^4$. If $G=S^3\times T^2$ then $H_0$ would have to be a two dimensional subgroup of $S^3\times 1$ for the action to be nonreducible, which is impossible. Next, if $G=S^3\times S^1$, then $H_0$ would be a one dimensional subgroup of $S^3\times 1$, say $H_0=S^1\times 1$. Then by \ref{abfactor}, $K^-$, $K^+$ and $H$ are all connected and we can assume $K^-=S^1\times S^1$. This proposition also says $K^+$ is either $S^3\times 1$ or $S^1\times S^1$. So there are two possibilities for diagrams in this case:
\begin{eqnarray}\label{M4_4}
   \nonumber S^3\times S^1 \,\, \supset \,\, S^1\times S^1, \,\, S^1\times S^1  \,\, \supset \,\,  S^1\times 1\\
   \nonumber S^3\times S^1 \,\, \supset \,\, S^1\times S^1, \,\, S^3\times 1  \,\, \supset \,\,  S^1\times 1.
\end{eqnarray}
The first is a product action on $S^2\times S^2$ and the second is a sum action on $S^4$.

The only remaining case is $G=S^3$, where $H$ is discrete. From \ref{fixedpt}, we can assume that $K^\pm$ are both proper subgroups of $G$, and hence they are both circle groups. After conjugation we can assume $K^-_0=\set{e^{i\theta}}$, and say $K^+_0=\set{e^{x\theta}}$ for some $x\in \Im \H \cap S^3$. If $x=\pm i$ then it is clear from \ref{Hgens} that $H$ must be a cyclic subgroup of $K^-=K^+$. In this case we have the group diagram:
\begin{equation}\label{M4_3a}
  S^3 \,\, \supset \,\, S^1, \,\, S^1  \,\, \supset \,\,  \Z_n
\end{equation}
One can show that this is either an action on $S^2\times S^2$ or on $\CP^2\#-\CP^2$, depending on whether $n$ is even or odd, respectively (see \cite{Parker}).

Next suppose $x\ne \pm i$. We know that for an arbitrary closed subgroup, $L\subset N(L_0)$. In particular $K^-\subset (\set{e^{i\theta}}\cup\set{je^{i\theta}})$ and $K^+\subset (\set{e^{x\theta}}\cup\set{ye^{x\theta}})$ for some $y\in x^\bot\cap \Im \H \cap S^3$. Hence $H$ must be a subgroup of the intersection of these two sets. If $x\notin i^\bot$ then these sets intersect in $\set{\pm1}\cup (\set{je^{i\theta}}\cap\set{ye^{x\theta}})$. However, by \ref{Hgens}, $H$ must be generated by its intersections with $K^-_0$ and $K^+_0$. Therefore $H\subset \set{\pm1}$ in this case, and hence $N(H)=G$. Then by \ref{normalizer}, we could conjugate $K^+$ by some element of $G$ to make $K^+_0=K^-_0$, which is the case we already considered.

So we can assume that $i\perp x$. Recall that conjugation of $S^3$ by the element $e^{i\theta_0}$ rotates the $jk$-plane by the angle $2\theta_0$ and fixes the $1i$-plane. After such a conjugation of $G$ we can assume that $K^+_0=\set{e^{j\theta}}$ and $K^-_0=\set{e^{i\theta}}$. This time
  $$H\subset N(K^-_0)\cap N(K^+_0)= \set{\pm1, \pm i, \pm j, \pm k}=:Q.$$
And we can also assume $H\nsubseteq \set{\pm1}$ by the argument above. Then $H$ contains some element of $Q \setminus \set{\pm1}$, and since $H$ is generated by its intersection with $K^-_0$ and $K^+_0$, we can assume $i\in H$.

If $H=\langle i\rangle$ then we have the following diagram:
\begin{equation}\label{M4_3b}
  S^3 \,\, \supset \,\, \set{e^{i\theta}}, \,\, \set{e^{j\theta}}\cup\set{ie^{j\theta}}  \,\, \supset \,\,  \langle i\rangle.
\end{equation}
One can easily check that this is the $\SO(3)$ action on $\CP^2$, via $\SO(3)\subset \SU(3)$. This action was missing from the classification in \cite{Parker}, along with other reducible actions.
If $H$ contains any other element of $Q$, in addition to $\langle i\rangle$, then $H=Q$, and we have the following diagram:
\begin{equation}\label{M4_3c}
  S^3 \,\, \supset \,\, \set{e^{i\theta}}\cup\set{je^{i\theta}}, \,\, \set{e^{j\theta}}\cup\set{ie^{j\theta}}  \,\, \supset \,\,  \set{\pm1, \pm i, \pm j, \pm k}.
\end{equation}
This is the action of $\SO(3)$ on $S^4\subset V:= \set{A\in \R^{3\times 3} | A=A^t, \tr(A)=0}$ by conjugation, where $V$ carries the inner product $\langle A,B \rangle= \tr(AB)$, as described in \cite{GZ1}.

The full classification of 4-dimensional compact cohomogeneity one manifolds given in \cite{Parker}, finds over 60 families of actions. The above work shows how important the assumption $\pi_1(M)=0$ is in simplifying the classification.

%
%
%
%
%
%

\section{Classification in Dimension Five}\label{5dim}

In this section we will go through the five dimensional classification. Throughout this section, $M$ will denote a 5-dimensional compact simply connected cohomogeneity one manifold for the compact connected group $G$ which acts almost effectively and nonreducibly, with group diagram $G \supset K^-, K^+ \supset H$ where $K^\pm/H\approx S^{l_\pm}$.

We will complete the classification by finding all such group diagrams which give simply connected manifolds. The first step is to find the possibilities for $G$. Since we are allowing the action to have finite ineffective kernel, after lifting the action to a covering group of $G$, we can assume that $G$ is a product of groups from \ref{posgrps}. In fact we have the following proposition.

\begin{prop}\label{5dimgrpsprop}
$G$ and $H_0$ must be one of the pairs of groups listed in Table \ref{5dimgrps}, up to equivalence.
\end{prop}
\begin{proof}
We will first show that all the possibilities for $G$ are listed in the table. We know from \ref{dimG} that $4\le \dim G \le 10$ and $\dim H = \dim G - 4$ since the principal orbits $G/H$ are codimension one in $M$. Further, since $G$ is a product of groups from \ref{posgrps}, $G$ must have the form $(S^3)^m\times T^n$, $\SU(3)\times T^n$ or $\Spin(5)$. From \ref{abfactor} we can assume $n\le2$ in all cases. First suppose that $G=(S^3)^m\times T^n$. Then by \ref{prods} we have $3m+n-4= \dim H \le m$ which means $0\le 4 - 2m -n$. Hence $m\le 2$ and if $m=2$ then $n=0$. So all the possibilities for groups of the form $(S^3)^m\times T^n$ are in fact listed in the table. Next suppose $G=\SU(3)\times T^n$. Then by \ref{prods} again we know that $\dim H \le 4$ which means $\dim G = \dim H + 4 \le 8$. Hence $\SU(3)$ is the only possibility of this form. Therefore all of the possible groups $G$ are listed in Table \ref{5dimgrps}.

Now we will show that for each possible $G$ described above we have listed all the possible subgroups $H_0$ of the right dimension. It is clear that if $G=S^3\times S^1$ then $H$ is discrete. Next, if $G=S^3\times T^2$, then for the action to be nonreducible, $\proj_2(H)\subset T^2$ must be trivial. Hence, $H_0$ is a closed connected one dimensional subgroup of $S^3$, as stated. If $G=S^3\times S^3$ then it is clear that $T^2$ must be a maximal torus in $G$. If $G=\SU(3)$, we see from \ref{grplist} that $H_0$ must be $\U(2)$ up to conjugation. Finally, \ref{topdim} deals with the last case where $\dim G= 10$.
\end{proof}

{\setlength{\tabcolsep}{0.40cm}
\renewcommand{\arraystretch}{1.6}
\stepcounter{equation}
\begin{table}[!h]
      \begin{center}
          \begin{tabular}{|c||c c|}
\hline
No.        & $G$   & $H_0$  \\
\hline \hline
1     & $S^3\times S^1$ & $\set{1}$   \\
\hline
2     & $S^3\times T^2$    &  $S^1\times 1$  \\
\hline
3     & $S^3\times S^3$    &   $T^2$  \\
\hline
4     & $\SU(3)$    &   $\U(2)$  \\
\hline
5     & $\Spin(5)$    &   $\Spin(4)$  \\
\hline
          \end{tabular}
      \end{center}
      \vspace{0.1cm}
      \caption{Possibilities for $G$ and $H_0$, in the 5-dimensional case.}\label{5dimgrps}
\end{table}}

In the rest of the section we proceed case by case to find all possible diagrams for the pairs of groups listed in \ref{5dimgrps}. We will do this by finding the possibilities for $K^\pm$, with $K^\pm/H$ a sphere. Recall, from \ref{no_exceptionals} and \ref{fixedpt}, that we can assume
  $$\dim G > \dim K^\pm > \dim H.$$

\subsection{Case 1: $G = S^3\times S^1$}\label{dim5.4}

In this section we consider the case that $G= S^3\times S^1$. In this case $H$ must be discrete. It then follows that for $K/H$ to be a sphere $K_0$ itself must be a cover of a sphere. Then from \ref{grplist}, the only compact connected subgroups of $S^3\times S^1$ that cover spheres are $S^3\times 1$ or circle groups of the form $\set{(e^{xp\theta},e^{iq\theta})}$ where $x\in \Im(\H)$. Further, from \ref{abfactor}, we know that at least one of $K^\pm_0$ is a circle. This leads us into the following cases: both $K^\pm_0$ are circles or $K^-_0$ is a circle and $K^+_0=S^3\times 1$.

\begin{case}[A] $K^-_0$ is a circle and $K^+_0=S^3\times 1$.

First, from \ref{fundgrp2}, $K^-$ must be connected with $H\subset K^-$ and $H\cap K^+_0=1$. After conjugation of $G$, we may assume $K^-=\set{(e^{ip\theta},e^{iq\theta})}$ and $p,q\ge0$. We also know from \ref{fundgrp2} that for $M$ to be simply connected $G/K^-=S^3\times S^1 / \set{(e^{ip\theta},e^{iq\theta})}$ must also be simply connected. It is not hard to see that this happens precisely when $q=1$. Finally, if $H=\Z_n\subset K^-$ the condition that $H\cap (1\times S^1)=e$ means $(p,n)=1$. Then $K^+=K^+_0\cdot H=S^3\times\Z_n$. In conclusion, such an action must have the following group diagram:
\begin{equation}\label{M5_4A}
   S^3\times S^1 \,\, \supset \,\, \set{(e^{ip\theta},e^{i\theta})}, \,\, S^3\times\Z_n  \,\, \supset \,\,  \Z_n
\end{equation}
Conversely, given such groups they clearly determine a simply connected cohomogeneity one manifold, by \ref{fundgrp}. This is action $Q^5_C$ of the appendix.
\end{case}

\begin{case}[B] Both $K^\pm_0$ are circle groups.

After conjugation we can take
\begin{equation}\label{dim5.4B_K}K^-_0 = \set{(e^{ip_-\theta},e^{iq_-\theta})} \text{ and } K^+_0=\set{(e^{xp_+\theta},e^{iq_+\theta})}
\end{equation}
for some $x\in\Im(\H)\cap S^3$ and $(p_\pm,q_\pm)=1$. From \ref{Hgens}, we know that $H$ must be generated by $H_-=H\cap K^-_0$ and $H_+=H\cap K^+_0$, which are cyclic subgroups of the circles $K^-_0$ and $K^+_0$, respectively.

We will now have to break this into two more cases, depending on whether $K^-_0$ and $K^+_0$ are both contained in a torus $T^2$ in $G$.
\begin{case}[B1] $K^-_0$ and $K^+_0$ are not both contained in any torus $T^2\subset G$.

Here we can assume that $x\ne\pm i$. Further, from \ref{dim5.4B_K}, we see $p_\pm\ne 0$ in this case, since otherwise $K^\pm_0$ would be contained in the same torus. Further, from \ref{abfactor}, we know that at least one of $q_\pm$ must be nonzero, say $q_+\ne0$. A computation shows that $N(K^+_0) = \set{(e^{x\theta},e^{i\phi})}$, since $p_+q_+\ne0$. Therefore $K^+\subset \set{(e^{x\theta},e^{i\phi})}$, since every compact subgroup of a Lie group is contained in the normalizer of its identity component. Similarly,
  $$K^-\subset N(K^-_0) = \left\{ \begin{array}{cc}
     \set{(e^{i\theta},e^{i\phi})} & \text{ if } q_-\ne0\\
     \set{(e^{i\theta},e^{i\phi})}\cup \set{(je^{i\theta},e^{i\phi})} & \text{ if } q_-=0
  \end{array} \right.$$

Therefore $H\subset K^-\cap K^+\subset N(K^-_0)\cap N(K^+_0)$. If $q_-\ne 0$ then this means $H \subset \set{(e^{i\theta},e^{i\phi})} \cap \set{(e^{x\theta},e^{i\phi})} = \set{(\pm1, e^{i\phi})}$. Then $H$ lies in the center of $G$ and so by \ref{normalizer}, we can conjugate $K^+$ by any element of $G$, and still have the same action. In particular we can conjugate $K^+$ to lie in the same torus as $K^-$, hence reducing such actions to Case B2. So we can assume that $q_-=0$ and hence $K^-_0 = \set{(e^{i\theta},1)}$.

Then, we have $H \subset N(K^-_0)\cap N(K^+_0)=\big( \set{(e^{i\theta},e^{i\phi})}\cup \set{(je^{i\theta},e^{i\phi})} \big) \cap \set{(e^{x\theta},e^{i\phi})}$. This intersection will again be $\set{(\pm1, e^{i\phi})}$ unless $x\perp i$. As above, we can again assume $x\perp i$. Further, after conjugation of $G$ by $(e^{i\theta_0}, 1)$, for a certain value of $\theta_0$, $K^-_0$ will remain fixed and $K^+_0$ will be taken to $\set{(e^{jp_+\theta}, e^{iq_+\theta})}$, with $p_+,q_+>0$. So we can assume
  $$K^-_0 = \set{(e^{i\theta},1)} \text{ and } K^+_0 = \set{(e^{jp_+\theta}, e^{iq_+\theta})}.$$
And therefore, $H \subset N(K^-_0)\cap N(K^+_0)=\set{\pm1, \pm j}\times S^1 \subset S^3\times S^1$. We saw above that we can assume $H$ is not contained in $\set{(\pm1, e^{i\phi})}$, and hence $H$ must contain an element of the form $(j,z_0)$, which we can assume also lies in $K^+_0$, by \ref{Hgens}. We can also assume that $H\cap (1\times S^1)=1$, so that $\#(z_0)|\#(j)=4$ and hence $z_0\in \set{\pm 1, \pm i}$, where $\#(g)$ denotes the order of the element $g$. So $H\cap K^+_0$ is generated by $(j,z_0)$. Similarly, $H\cap K^-_0$ must also be a subset of $\set{\pm1, \pm j}\times 1$ and is therefore either trivial or $\set{(\pm1,1)}$. For convenience we break this up into three more cases, depending on the order of $z_0$.

\begin{case}[B1a]The order of $z_0$ is one, i.e. $z_0=1$.

In this case $H=\langle(j,1)\rangle$, $K^-_0\cap H = \set{(\pm1, 1)}$ and $K^+_0\cap H=H$. Hence $K^+$ is connected. The condition that $H\subset K^+$ means $4|q_+$ and $p_+$ is odd. We can represent $\pi_1(K^+/H)$ with the curve $\alpha_+:[0,1] \to K^+ : t \mapsto (e^{2\pi jp_+t/4},e^{2\pi iq_+t/4})$ and we can represent $\pi_1(K^-/H)$ with the curve $\alpha_-:[0,1] \to K^- : t \mapsto (e^{2\pi it/2},1)$. From \ref{Hgens}, $M$ will be simply connected if and only if $\alpha_\pm$ generate $\pi_1(G)$. We see that the possible loops in $G$ that $\alpha_\pm$ can form are combinations of $\alpha_-^2$, $\alpha_+^4$ and $\alpha_-\circ\alpha_+^2$. Yet each of these loops can only give an even multiple of the loop $1\times S^1 \subset S^3\times S^1 = G$, which generates $\pi_1(G)$. Hence $M$ will never be simply connected in this case.
\end{case}

\begin{case}[B1b]The order of $z_0$ is two, i.e. $z_0=-1$.

In this case $H=\langle(j,-1)\rangle$, and again $K^-_0\cap H = \set{(\pm1, 1)}$ and $K^+_0\cap H=H$, so that $K^+$ is connected. This time, the condition that $H\subset K^+$ means that $p_+$ is odd and $q_+\equiv 2 \mod 4$. Then, in this case we can represent $\pi_1(K^+/H)$ with the curve $\alpha_+:[0,1] \to K^+ : t \mapsto (e^{2\pi jp_+t/4},e^{2\pi iq_+t/4})$ and $\pi_1(K^-/H)$ with $\alpha_-:[0,1] \to K^- : t \mapsto (e^{2\pi it/2},1)$, and again $M$ will be simply connected if and only if $\alpha_\pm$ generate $\pi_1(G)$, by \ref{Hgens}. The loops that $\alpha_\pm$ can generate are again combinations of $\alpha_-^2$, $\alpha_+^4$ and $\alpha_-\circ \alpha_+^2$ but in this case $\alpha_-^2$ corresponds to zero times around the loop $1\times S^1$; $\alpha_+^4$ corresponds to $q_+$ times around $1\times S^1$; and $\alpha_-\circ \alpha_+^2$ corresponds to $q_+/2$ times around $1\times S^1$. Together with the constraints $q_+\equiv 2 \mod 4$ and $q_+>0$, we see that $M$ will be simply connected if and only if $q_+=2$. Therefore this case gives the family of actions:
\begin{eqnarray}\label{M5_4B1b}
   S^3\times S^1 \,\, \supset \,\, \set{(e^{i\theta},1)}\cdot H, \,\, \set{(e^{jp_+\theta}, e^{2i\theta})}  \,\, \supset \,\,  \langle(j,-1)\rangle\\
   \nonumber\text{  where  } p_+>0 \text{ is odd.}
\end{eqnarray}
These are the actions of family $Q^5_B$ in the appendix.
\end{case}

\begin{case}[B1c]The order of $z_0$ is four, i.e. $z_0=\pm i$.

After a conjugation of $G$ which will not effect the form of $K^\pm$ we can assume $z_0=i$, that is $(j,i)\in H\cap K^+_0$, although we can no longer assume $p_+>0$. As explained above, $H\cap K^-_0\subset \set{(\pm1,1)}$. Yet if $(-1,1)\in H$ then $(-1,-1)\cdot(-1,1)=(1,-1)\in H$, violating our assumption that $H\cap 1\times S^1 = 1$. Therefore $H=\langle(j,i)\rangle\subset K^+$, $K^+$ is connected and $K^-_0\cap H = 1$. This also implies that $p_+$ and $q_+$ are odd and $p_+\equiv q_+ \mod 4$. In this case $\pi_1(K^+/H)$ can be represented by the curve $\alpha_+:[0,1] \to K^+ : t \mapsto (e^{2\pi jp_+t/4},e^{2\pi iq_+t/4})$ and $\pi_1(K^-/H)$ can be represented by $\alpha_-:[0,1] \to K^- : t \mapsto (e^{2\pi it},1)$. As above, by \ref{Hgens}, $M$ will be simply connected if and only if $\alpha_\pm$ generate $\pi_1(G)$. We see that the only loops in $G$ that $\alpha_\pm$ can generate are $\alpha_-$ and $\alpha_+^4$ where $\alpha_-$ in trivial in $\pi_1(G)$ and $\alpha_+^4$ represents $q_+$ times around the loop $1\times S^1$, which generates $\pi_1(G)$. Together with our assumption that $q_+>0$, we get that $M$ is simply connected if and only if $q_+=1$. This case gives precisely the following family of actions:
\begin{eqnarray}\label{M5_4B1b}
   S^3\times S^1 \,\, \supset \,\, \set{(e^{i\theta},1)}\cdot H, \,\, \set{(e^{jp_+\theta}, e^{i\theta})}  \,\, \supset \,\,  \langle(j,i)\rangle\\
   \nonumber\text{  where  } p_+\equiv 1 \mod 4
\end{eqnarray}
These actions make up the family $P^5$, of the appendix.
\end{case}

\end{case}
\begin{case}[B2] $K^-_0$ and $K^+_0$ are both contained in a torus $T^2\subset G$.

After conjugation of $G$ we may assume that $K^\pm_0\subset \set{(e^{i\theta},e^{i\phi})}$. It then follows from \ref{Hgens} that $K^\pm, H \subset \set{(e^{i\theta},e^{i\phi})}$. Here again there will be two cases depending on whether or not $K^-_0$ and $K^+_0$ are distinct circles.
\begin{case}[B2a] $K^-_0=K^+_0$.

Then we can take $K^-_0=K^+_0=\set{(e^{ip\theta},e^{iq\theta})}$ with $q>0$. From \ref{Hgens} it follows that $H$ is a cyclic subgroup of $K^\pm_0$ and $K^\pm$ is connected. It is then clear from \ref{Hgens} that $q=1$ and hence we have the following family of actions:
\begin{equation}\label{M^5_4B2a}
  S^3\times S^1\,\, \supset \,\, \set{(e^{ip\theta},e^{i\theta})}, \,\, \set{(e^{ip\theta},e^{i\theta})} \,\, \supset \,\, \Z_n.
\end{equation}
Conversely, the resulting manifolds will all be simply connected by \ref{Hgens}. This family is labeled $Q^5_A$ in the appendix.
\end{case}
\begin{case}[B2b] $K^-_0\ne K^+_0$.

Here, say $K^\pm_0= \set{(e^{ip_\pm\theta},e^{iq_\pm\theta})}$ and then $H=H_-\cdot H_+$ for cyclic subgroups $H_\pm \subset K^\pm_0$, by \ref{Hgens}. Now let $\alpha_\pm:[0,1]\to K^\pm_0$ be curves with $\alpha_\pm(0)=1$ which represent $\pi_1(K^\pm/H)$. Then by \ref{Hgens}, $M$ will be simply connected if and only if the combinations of $\alpha_\pm$ which form loops in $G$, generate $\pi_1(G)$. Let $\delta_\pm:[0,1] \to K^\pm_0:t\mapsto (e^{2\pi ip_\pm t},e^{2\pi iq_\pm t})$ be curves that pass once around the circles $K^\pm_0$. Then $\delta_\pm = \alpha_\pm^{m_\pm}$ are two such loops.

To find all such loops, consider the covering map $\wp: \R^2\to T^2\subset S^3\times S^1: (x,y)\mapsto (e^{2\pi ix}, e^{2\pi iy})$ and let $\widetilde K^\pm$ be lines through the origin and through the points $(p_\pm,q_\pm)$ in $\R^2$. Then it is clear that $\wp^{-1} (K^\pm_0) = \widetilde K^\pm + \Z^2$.

Now let $\widetilde \gamma$ be a path in $\R^2$ which starts at $(0,0)$ follows $\widetilde K^-$ until the first point of the intersection $\widetilde K^- \cap\big( \widetilde K^+ + \Z^2 \big)$ then follows $\widetilde K^+ + \Z^2$ to the first integer lattice point $(\lambda, \mu)$. Then $\gamma:=\wp(\widetilde \gamma)$ gives a loop in $G$. Notice that $K^-_0\cap K^+_0$ is a cyclic subgroup of both $K^\pm_0$ and any curve in $K^-_0\cup K^+_0$ is homotopic within $T^2$ to a curve in $K^-_0$ followed by a curve in $K^+_0$. It then follows that $\delta_-$, $\delta_+$ and $\gamma$ generate all possible loops in $K^-_0 \cap K^+_0$.
Similarly, if $d$ is the index of $H\cap K^-_0\cap K^+_0$ in $K^-_0\cap K^+_0$, then $\gamma^d$ can be imagined as a curve that starts at 1, travels along $K^-_0$ to the first element of $H$ in $H\cap K^-_0\cap K^+_0$, and then follows $K^+_0$ back to the identity. Then $\delta_-$, $\delta_+$ and $\gamma^d$ generate the same homotopy classes of loops as $\alpha_-$ and $\alpha_+$. Therefore $M$ will be simply connected if and only if $\delta_-$, $\delta_+$ and $\gamma^d$ generate $\pi_1(G)$.

Let $c:[0,1]\to G=S^3\times S^1:t\mapsto (1,e^{2\pi it})$ represent the generator of $\pi_1(G)$. Then it is clear that $\delta_\pm$ is homotopic to $c^{q_\pm}$ in $G$ and that $\gamma$ is homotopic to $c^\mu$ in $G$. Therefore $M$ is simply connected if and only if $\langle c^{q_-}, c^{q_+}, c^{d\mu} \rangle = \langle c \rangle$. Notice further that by the construction of $\widetilde \gamma$, $(\lambda, \mu)$ is an integer lattice point which is closest to the line $\widetilde K^+_0$. Therefore $(p_+,q_+)$ and $(\lambda,\mu)$ generate all of $\Z^2$ and in particular $q_+$ and $\mu$ are relatively prime. Hence $\langle c^{q_-}, c^{q_+}, c^{d\mu} \rangle = \langle c^{q_-}, c^{q_+}, c^d \rangle$ and therefore $M$ is simply connected if and only if $\gcd(q_-,q_+,d)=1$. Therefore, we get precisely the following family of simply connected diagrams
\begin{eqnarray}\label{M5_4B1b}
   S^3\times S^1 \,\, \supset \,\, \set{(e^{ip_-\theta},e^{iq_-\theta})}\cdot H, \,\, \set{(e^{ip_+\theta},e^{iq_+\theta})}\cdot H  \,\, \supset \,\,  H_-\cdot H_+\\
   \nonumber K^-\ne K^+,  \gcd(q_-,q_+,d)=1 \text{ where } d= \#(K^-_0\cap K^+_0)/\#(H\cap K^-_0\cap K^+_0).
\end{eqnarray}
This is the family labeled $N^5$ in the appendix.
\end{case}

\end{case}
\end{case}

\subsection{The remaining cases}

In this section we will examine the remaining possibilities for $G$: Cases 2-5 from Table \ref{5dimgrps}.

\subsubsection{Cases 2, 4 and 5}

In Case 2, $G=S^3\times T^2$ and $H_0=S^1\times 1$ where $\rank S^3=\rank S^1$. Proposition \ref{abfactor}, then says that the resulting action must be a product action. In Case 4, we know from \ref{grplist} that $H_0=\U(2)$ is maximal among connected subgroups of $G=\SU(3)$. Hence any action with these groups would be an isometric two-fixed-point action on a sphere, by \ref{Hmax}. Finally, Case 5 is fully described by \ref{topdim}.

\subsubsection{Case 3: $G=S^3\times S^3$}\label{dim5.6}

Now $G=S^3\times S^3$ and $H_0=S^1\times S^1$. Then from \ref{grplist}, any proper connected subgroup $K$ of $G$, containing $H_0$ and of higher dimension, must be $S^3\times S^1$ or $S^1\times S^3$. Then, since our only possibilities for $K$ have $\dim(K/H)=2$, \ref{fundgrp2} implies that all of $H$, $K^-$ and $K^+$ are connected. Therefore, up to equivalence, we only have the following possible diagrams:
\begin{equation}\label{sofar3.3_1}
   S^3\times S^3 \,\, \supset \,\, S^3\times S^1, \,\, S^3\times S^1 \,\, \supset \,\, S^1\times S^1
\end{equation}
\begin{equation}\label{sofar3.3_2}
   S^3\times S^3 \,\, \supset \,\, S^3\times S^1, \,\, S^1\times S^3 \,\, \supset \,\, S^1\times S^1
\end{equation}
Conversely, it is clear that these both give simply connected manifolds. The first of these actions is a product action and the second is a sum action.

%
%
%
%
%
%

\section{Classification in Dimension Six}\label{6dim}

In this section we will carry out the 6-dimensional classification. Throughout this section we will keep the notations and conventions established at the beginning of Section \ref{5dim}, this time for a 6-dimensional manifold $M$. As in the previous case we have the following result to describe the possible groups.

%

\begin{prop}\label{6dimgrpsprop}
$G$ and $H_0$ must be one of the pairs of groups listed in Table \ref{6dimgrps}, up to equivalence.
\end{prop}
\begin{proof}
We first show that all the possibilities for $G$ are listed in the table. We know from \ref{dimG} that $5 \le \dim G \le 15$ in this case and $\dim G = \dim H + 5$ since $\dim(G/H)=5$. From \ref{posgrps}, $G$ must have the form $(S^3)^m\times T^n$, $\SU(3)\times (S^3)^m\times T^n$, $\Sp(2)\times (S^3)^m\times T^n$, $\G_2\times T^n$, or $\Spin(6)$. Further, by \ref{abfactor}, we can assume $n\le2$ in all cases.

First suppose $G = (S^3)^m\times T^n$. Then by \ref{prods}, $3m + n - 5 = \dim H \le m$ and hence $0\le 5 - 2m - n$. Therefore $m\le 2$ and if $m=2$ then $n\le1$. We see that all of these possibilities are recorded in the table. Next assume $G= \SU(3) \times (S^3)^m \times T^n$. Then by \ref{prods} again we know $8 + 3m +n - 5 = \dim H \le 4 + m$ or $0\le 1 - 2m - n$. Hence $m=0$ and $n\le 1$. Note again that these two possibilities for $G$ are listed in the table. Next if $G= \Sp(2)\times (S^3)^m\times T^n$, \ref{prods} gives $0\le 1 - 2m - n$ again. So again $m=0$ and $n\le 1$. However, if $G= \Sp(2)$ then $\dim H = 5$ and $\rank H \le \rank G=2$. Yet there are no 5-dimensional compact groups of rank two or less, by \ref{posgrps}. So in fact, $\Sp(2)$ is not a possibility for $G$. Finally suppose that $G = \G_2\times T^n$. Then by \ref{prods}, $\dim H\le 8$ and yet $H$ would have to be $9+n$ dimensional in this case. Hence this is not a possibility either.

Now we will show that for each $G$ in the table, all the possibilities for $H_0$ are listed. First, if $G=G_1\times T^m$, then for the action to be nonreducible, we can assume $\proj_2(H_0)$ is trivial in these cases. Then, we use \ref{grplist} to find the possibilities for $H_0$ in each case, up to conjugation. For the last case, \ref{topdim} tells us the full story.
\end{proof}

{\setlength{\tabcolsep}{0.40cm}
\renewcommand{\arraystretch}{1.6}
\stepcounter{equation}
\begin{table}[!h]
      \begin{center}
          \begin{tabular}{|c||c c|}
\hline
No.   & $G$   & $H_0$  \\
\hline \hline
1     & $S^3\times T^2$ & $\set{1}$   \\
\hline
2     & $S^3\times S^3$    &  $\set{(e^{ip\theta},e^{iq\theta})}$  \\
\hline
3     & $S^3\times S^3\times S^1$    &   $T^2\times 1$  \\
\hline
4     & $\SU(3)$    &   $\SU(2)$, $\SO(3)$  \\
\hline
5     & $\SU(3)\times S^1$    &   $\U(2)\times 1$  \\
\hline
6     & $\Sp(2)\times S^1$    &   $\Sp(1)\Sp(1)\times 1$  \\
\hline
7     & $\Spin(6)$    &   $\Spin(5)$  \\
\hline
          \end{tabular}
      \end{center}
      \vspace{0.1cm}
      \caption{Possibilities for $G$ and $H_0$, in the 6-dimensional case.}\label{6dimgrps}
\end{table}}

We will now continue with the classification on a case by case basis. As in dimension 5, the case that $H$ is discrete is the most difficult.

\subsection{Case 1: $G=S^3\times T^2$}

Here $G=S^3\times T^2$ and $H$ is discrete. By \ref{abfactor}, we see that $K^\pm_0$ must both be circle groups in $G$, say $K^\pm_0 = \set{(e^{x_\pm a_\pm\theta},e^{ib_\pm\theta},e^{ic_\pm\theta})}$ for $x_\pm \in \Im S^3$, where $(b_-,c_-)$ and $(b_+,c_+)$ are linearly independent. After conjugation we can assume that $x_-=i$ and we claim we can also assume that $x_+=i$. If one of $a_\pm$ is zero then this is clear. Otherwise we have $N(K^-_0)=\set{(e^{i\theta},e^{i\phi},e^{i\psi})}$ and $N(K^+_0)=\set{(e^{x_+\theta},e^{i\phi},e^{i\psi})}$ and $H \subset N(K^-_0)\cap N(K^+_0)= \set{(\pm1,e^{i\phi},e^{i\psi})}$ if $x_+\ne\pm i$. But then $H$ would be normal in $G$ and by \ref{normalizer}, we would be able to conjugate $K^+$ to make $x_+=i$ without affecting the resulting manifold. So we can assume
  $$K^\pm_0 = \set{(e^{i a_\pm\theta},e^{ib_\pm\theta},e^{ic_\pm\theta})}.$$

Let $\alpha_\pm:[0,1]\to K^\pm_0$ be curves with $\alpha_\pm(0)=1$ which represent $\pi_1(K^\pm/H)$. \ref{Hgens} says that $M$ will be simply connected if and only if $H$ is generated by $\alpha_\pm(1)$ as a group and $\alpha_\pm$ generate $\pi_1(G)$. Assume that $H$ is generated by $\alpha_\pm(1)$ and we will find the conditions under which $\alpha_\pm$ generate $\pi_1(G)$.

Notice that \ref{Hgens} implies that $K^\pm$ and $H$ must all be contained in $T^3=\set{(e^{i\theta},e^{i\phi},e^{i\psi})}$, in order for $M$ to be simply connected. Now consider the cover $\wp: \R^3\to T^3: (x,y,z)\mapsto (e^{2\pi i x},e^{2\pi iy},e^{2\pi iz})$. In $\R^3$, let $\widetilde K^\pm$ but the lines through the origin and the points $(a_\pm,b_\pm,c_\pm)$. Then it is clear that $\wp^{-1}(K^\pm_0)=\widetilde K^\pm + \Z^3$. Next, denote the plane spanned by $\widetilde K^\pm$ by $Q$ and the lattice $Q\cap \Z^3$ by $L$.

We then see that any loop generated by $\alpha_\pm$ will lift to a path in $Q$ from the origin to a point in $L=Q\cap \Z^3$. Finally define the map $\fp: \R^3 \to \R^2: (x,y,z)\mapsto (y,z)$. Then we have an isomorphism of $\pi_1(G)\to \Z^2$ given as follows: for $[c]\in \pi_1(G)$ lift $c$ to a curve $\tilde c$ in $\R^3$ starting from the origin via $\wp$, then $[c]\mapsto \fp(\tilde c(1))$. It is clear that for the combinations of $\alpha_\pm$ which form loops in $G$ to generate $\pi_1(G)$, we must at least have $\fp(L)=\Z^2$. This means that $L$ must have the form
  $$L=\set{(f(i,j),i,j)|i,j\in \Z}$$
for some function $f$ of the form $f(i,j)=ri+sj$ with fixed $r,s\in \Z$. In particular, it follows that $a_\pm=f(b_\pm,c_\pm)=rb_\pm+sc_\pm$ since $(a_\pm,b_\pm,c_\pm)\in L$. Hence $\gcd(b_\pm,c_\pm)=1$ since we assumed that $\gcd(a_\pm,b_\pm,c_\pm)=1$.

Now define the curve $\widetilde \gamma:[0,1]\to \R^3$ as follows: $\widetilde \gamma$ starts at the origin, follows $\widetilde K^-$ to the first point of intersection in $\big(\widetilde K^+ + \Z^3\big)\cap \widetilde K^-$, then follows $\widetilde K^+ + \Z^3$ to the first integer lattice point $(f(\lambda, \mu),\lambda, \mu)$ in $\Z^3$.
We claim that $(b_+,c_+)$ and $(\lambda,\mu)$ generate $\Z^2$. To see this note that, by the construction of $\widetilde \gamma$, the point $(f(\lambda, \mu),\lambda, \mu)$ is a point in $L$ which is closest to the line $\widetilde K^+$. Hence $(a_+,b_+,c_+)$ and $(f(\lambda, \mu),\lambda, \mu)$ generate $L$ and so $(b_+,c_+)$ and $(\lambda, \mu)$ generate $\Z^2$.

Define $\widetilde \delta_\pm:[0,1]\to \R^3:t\mapsto t(a_\pm,b_\pm,c_\pm)$ and let $\gamma=\wp(\widetilde \gamma)$ and $\delta_\pm=\wp(\widetilde \delta_\pm)$. If $d$ denotes the index of $H\cap K^-_0\cap K^+_0$ in $K^-_0\cap K^+_0$ then we claim that $\delta_-$, $\delta_+$ and $\gamma^d$ generate the same subgroup of $\pi_1(G)$ as $\alpha_-$ and $\alpha_+$. To see this, notice that $\alpha_\pm$ can be taken to be paths in $K^\pm_0$ from the identity to the first element of $H\cap K^\pm_0$ and that any combination of $\alpha_\pm$ which forms a loop in $G$ can be expressed as a curve in $K^-_0$, from the identity to an element of $H\cap K^-_0\cap K^+_0$, followed by a curve in $K^+_0$ back to the identity. We see from the construction of $\widetilde \gamma$ that $\gamma^d$ is a loop from the identity, along $K^-_0$ to the first element of $H\cap K^-_0\cap K^+_0$, then around $K^+_0$ some number of times before returning to the identity. Since $H\cap K^-_0\cap K^+_0$ is a cyclic subgroup $K^-_0$, we see that any loop generated by $\alpha_-$ and $\alpha_+$ can be expressed as a power of $\gamma^d$ followed by a power of $\delta_+$. So  $\delta_-$, $\delta_+$ and $\gamma^d$ generate the same subgroup of $\pi_1(G)$ as $\alpha_-$ and $\alpha_+$.

Then by \ref{Hgens}, $M$ will be simply connected if and only if $\delta_-$, $\delta_+$ and $\gamma^d$ generate $\pi_1(G)$. Via the isomorphism $\pi_1(G)\to \Z^2$ described above, $\delta_-$, $\delta_+$ and $\gamma^d$ correspond to $(b_-,c_-)$, $(b_+,c_+)$ and $d(\lambda,\mu)$, respectively. So $M$ is simply connected if and only if $d(\lambda,\mu)$ generates $\Gamma = \Z^2/\langle(b_-,c_-),(b_+,c_+)\rangle$. We saw above that $(b_+,c_+)$ and $(\lambda,\mu)$ themselves generate $\Z^2$. Therefore $\Gamma$ is a cyclic group generated by $(\lambda,\mu)$. We see the order of $\Gamma = \Z^2/\langle(b_-,c_-),(b_+,c_+)\rangle$ is $n=\pm(b_-c_+-c_-b_+)=\#(K^-_0\cap K^+_0)$. Then from the definition of $d$ we get that $d|n$. Hence $d(\lambda,\mu)$ generates $\Gamma$ if and only if $d=1$. Therefore $M$ is simply connected if and only if $K^-_0\cap K^+_0\subset H$.

Hence we have precisely the following family of diagrams:
\begin{eqnarray}\label{M6_5}
   S^3\times T^2 \,\, \supset \,\, \set{(e^{i a_-\theta},e^{ib_-\theta},e^{ic_-\theta})}\cdot H, \,\, \set{(e^{i a_+\theta},e^{ib_+\theta},e^{ic_+\theta})}\cdot H \,\, \supset \,\, H\\
   \nonumber\text{where $K^-\ne K^+$, $H=H_-\cdot H_+$, $\gcd(b_\pm,c_\pm)=1$,}\\
   \nonumber\text{$a_\pm=rb_\pm+sc_\pm$, and $K^-_0\cap K^+_0\subset H$.}
\end{eqnarray}

This is the family $N^6_A$ of the appendix. Notice that we can eliminate several parameters from the expression of the group diagram above. After an automorphism of $G$ we can assume that $K^-_0\subset S^3\times S^1\times 1$ and hence $(a_-,b_-,c_-)=(r,1,0)$. However the symmetric presentation in \ref{M6_5} will be preferred for our purposes.

\subsection{The remaining cases}

We will now address Cases 2-7 from Table \ref{6dimgrps}.

\subsubsection{Case 2: $G=S^3\times S^3$}

Here $G = S^3\times S^3$ and $H_0 = \set{(e^{ip\theta},e^{iq\theta})}$, for $(p,q)=1$ and $p,q\ge0$, after conjugation of $G$. Then, from \ref{grplist}, the possible compact proper subgroups $K$ containing $H$ with $K/H\approx S^l$ are: any torus $T^2\supset H$; $S^3\times 1$ if $q=0$; $1\times S^3$ if $p=0$; $\Delta S^3$ if $p=q=1$; $S^3\times S^1$ where $S^3\times S^1/H\approx S^3$ if and only if $q=1$; or $S^1\times S^3$ where $S^1\times S^3/H\approx S^3$ if and only if $p=1$.

We will now break this into cases by pairing together all of the possibilities for $K^\pm$, remembering that we can switch the places of $K^-$ and $K^+$ without effecting the resulting action.

\begin{case}[A]$K^-_0$ and $K^+_0$ are both tori.

Here we need to break this up further into two more cases depending on whether or not $K^-_0$ and $K^+_0$ are the same torus.

\begin{case}[A1]$K^-_0$ and $K^+_0$ the same torus.

Here $K^-_0 = K^+_0 = T^2$ and hence, by \ref{Hgens}, $H \subset K^\pm_0$ and $K^\pm$ are both connected. We also see from \ref{Hgens} that any $H \subset K^\pm$ with $H_0=S^1$ will give a simply connected manifold. In general such groups $H$ will have the form $\set{(e^{ip\theta},e^{iq\theta})} \cdot \Z_n$ after a conjugation of $G$. Therefore we get the following family of actions in this case:
\begin{equation}\label{M5_6A1}
   S^3\times S^3 \,\, \supset \,\, \set{(e^{i\theta},e^{i\phi})}, \,\, \set{(e^{i\theta},e^{i\phi})} \,\, \supset \,\, \set{(e^{ip\theta},e^{iq\theta})} \cdot \Z_n
\end{equation}
This is the family $N^6_B$ of the appendix.
\end{case}
\begin{case}[A2]$K^-_0$ and $K^+_0$ are different tori.

For $K^-_0$ and $K^+_0$ to be different tori, both containing the circle $H_0 = \set{(e^{ip\theta},e^{iq\theta})}$, it follows that either $p$ or $q$ must be zero. Suppose, without loss of generality that $q = 0$, so that $H_0 = \set{(e^{i\theta},1)}$. It then follows that $K^\pm_0$ must have the form $K^\pm_0 = \set{(e^{i\theta},e^{x_\pm\phi})}$ for some $x_\pm\in \Im(S^3)$. Notice that for $M$ to be simply connected, by \ref{Hgens}, $H \subset \set{(e^{i\theta},g)} = S^1\times S^3$. Therefore $H$ and $K^\pm$ all have the form $H=S^1\times \widehat H$ and $K^\pm = S^1 \times \widehat K^\pm$, where $\widehat K^\pm/\widehat H \approx S^1$. That means $S^3 \supset \widehat K^+, \widehat K^- \supset \widehat H$ gives a four dimensional cohomogeneity one manifold. Further, from \ref{Hgens}, it follows that this 4-manifold will be simply connected if and only if $M$ is. Hence our action is a product action with some simply connected 4-manifold.
\end{case}
\end{case}

\begin{case}[B]$K^-_0 = T^2$ and $K^+_0 = S^3\times 1$. Then $q=0$ and $H_0 = \set{(e^{i\theta},1)}$.

From \ref{Hgens}, $H$ must be of the form $S^1\times \Z_n \subset T^2$. This gives the following family of group diagrams:
\begin{equation}\label{M5_6B}
   S^3\times S^3 \,\, \supset \,\, T^2, \,\, S^3\times \Z_n \,\, \supset \,\, S^1 \times \Z_n
\end{equation}
Conversely, these diagrams obviously determines simply connected manifolds, by \ref{Hgens}. This is the family $N^6_C$ of the appendix.
\end{case}

\begin{case}[C]$K^-_0 = T^2$ and $K^+_0 = \Delta S^3$. Then $p=q=1$ and $H_0 = \Delta S^1$.

Again, by \ref{Hgens}, $H$ will have the form $\Delta S^1 \cdot \Z_n$. Yet, every compact Lie group is contained in the normalizer of its identity component. In particular $\Delta S^3 \subset N(\Delta S^3) = \pm\Delta S^3=\set{(g, \pm g)}$. This means that $n$ is at most two. Therefore, we have the following two possibilities for group diagrams:
\begin{eqnarray}\label{M5_6C}
   S^3\times S^3 \,\, \supset \,\, T^2, \,\, \Delta S^3\cdot \Z_n \,\, \supset \,\, \Delta S^1 \cdot \Z_n\\
   \nonumber\text{  where  } n = 1\text{ or }2
\end{eqnarray}
From \ref{Hgens}, we see that these are both in fact simply connected. These are the actions in family $Q^6_A$.
\end{case}

\begin{case}[D]$K^-_0 = T^2$ and $K^+_0 = S^3\times S^1$. Then $q=1$ and $H_0 = \set{(e^{ip\theta},e^{i\theta})}$.

It is clear in this case that $K^-_0 \subset K^+_0$. Further, for $K^+/H$ to be a 3-sphere, $H\cap K^+_0 = H_0$. Therefore $H$ and $K^\pm$ are all connected. We then have the following family of diagrams which all give simply connected manifolds by \ref{Hgens}:
\begin{equation}\label{M5_6D}
   S^3\times S^3 \,\, \supset \,\, T^2, \,\, S^3\times S^1 \,\, \supset \,\, \set{(e^{ip\theta},e^{i\theta})}.
\end{equation}
We see that these are the actions of type $N^6_D$.
\end{case}

\begin{case}[E]$K^-_0 = S^3\times 1$ and $K^+_0 = S^3\times 1$. Then $q=0$ and $H_0 = \set{(e^{i\theta},1)}$.

From \ref{fundgrp2}, we know that $H$ and $K^\pm$ must all be connected in this case. We then have the following group diagram which gives a simply connected manifold by \ref{Hgens}:
\begin{equation}\label{M5_6E}
   S^3\times S^3 \,\, \supset \,\, S^3 \times 1, \,\, S^3\times 1 \,\, \supset \,\, S^1\times 1.
\end{equation}
We quickly notice that this is a product action.
\end{case}

\begin{case}[F]$K^-_0 = S^3\times 1$ and $K^+_0 = S^1\times S^3$. Then $q=0$, $p=1$ and $H_0 = \set{(e^{i\theta},1)}$.

As in the previous case, we get the following simply connected group diagram:
\begin{equation}\label{M5_6F}
   S^3\times S^3 \,\, \supset \,\, S^3 \times 1, \,\, S^1\times S^3 \,\, \supset \,\, S^1\times 1.
\end{equation}
We see that this is a sum action.
\end{case}

\begin{case}[G]$K^-_0 = \Delta S^3$ and $K^+_0 = \Delta S^3$. Then $p=q=1$ and $H_0 = \Delta S^1$.

As above we have the following group diagram:
\begin{equation}\label{M5_6G}
   S^3\times S^3 \,\, \supset \,\, \Delta S^3, \,\, \Delta S^3 \,\, \supset \,\, \Delta S^1.
\end{equation}
This is action $Q^6_B$ of the appendix.
\end{case}

\begin{case}[H]$K^-_0 = \Delta S^3$ and $K^+_0 = S^3\times S^1$. Then $p=q=1$ and $H_0 = \Delta S^1$.

Again, we have:
\begin{equation}\label{M5_6H}
   S^3\times S^3 \,\, \supset \,\, \Delta S^3, \,\, S^3\times S^1 \,\, \supset \,\, \Delta S^1
\end{equation}
which is action $Q^6_C$ of the appendix.
\end{case}

\begin{case}[I]$K^-_0 = S^3\times S^1$ and $K^+_0 = S^3\times S^1$. Then $q=1$ and $H_0 = \set{(e^{ip\theta},e^{i\theta})}$.

Here, as above, we have:
\begin{equation}\label{M5_6I}
   S^3\times S^3 \,\, \supset \,\, S^3\times S^1, \,\, S^3\times S^1 \,\, \supset \,\, \set{(e^{ip\theta},e^{i\theta})}
\end{equation}
which is the family $N^6_E$ of the appendix.
\end{case}

\begin{case}[J]$K^-_0 = S^3\times S^1$ and $K^+_0 = S^1\times S^3$. Then $p=q=1$ and $H_0 = \Delta S^1$.

Our last possibility in this case is the following diagram:
\begin{equation}\label{M5_6J}
   S^3\times S^3 \,\, \supset \,\, S^3\times S^1, \,\, S^1\times S^3 \,\, \supset \,\, \Delta S^1,
\end{equation}
labeled $Q^6_D$ in the appendix.
\end{case}

\subsubsection{Case 3: $G=S^3\times S^3\times S^1$}

Here $G = S^3\times S^3\times S^1$ and $H_0 = T^2\times 1\subset S^3\times S^3\times S^1$. By \ref{abfactor}, one of $K^\pm/H$ must be a circle, say $K^-/H\approx S^1$. Furthermore, since $\rank(H)=\rank(S^3\times S^3)$, \ref{abfactor} says $K^-=T^2\times S^1$, and all of $H$, $K^-$ and $K^+$ are connected. We will now find the possibilities for $K^+$. Notice that if $\proj_3(K^+)$ is nontrivial, then by \ref{abfactor}, $K^+=T^2\times S^1$, giving one possibility. Otherwise $K^+\subset S^3\times S^3\times 1$. In this case, $K^+$, which must contain $H$, must be one of $S^3\times S^1 \times 1$, $S^1\times S^3 \times 1$ or $S^3\times S^3 \times 1$. But $S^3\times S^3\times 1 / S^1\times S^1\times 1 \approx S^2\times S^2$ which is not a sphere. Putting this together, we see our only possible group diagrams, up to automorphism, are:
\begin{equation}\label{M^6_7a}
   S^3\times S^3 \times S^1 \,\, \supset \,\, S^1\times S^1 \times S^1, \,\, S^1\times S^1 \times S^1 \,\, \supset \,\, S^1\times S^1 \times 1,
\end{equation}
\begin{equation}\label{M^6_7b}
   S^3\times S^3 \times S^1 \,\, \supset \,\, S^1\times S^1 \times S^1, \,\, S^3\times S^1 \times 1 \,\, \supset \,\, S^1\times S^1 \times 1.
\end{equation}
We see, however, that both of these actions are product actions.

\subsubsection{Case 4: $G=\SU(3)$}

In this case $G=\SU(3)$ and $H_0$ must be $\SO(3)$ or $\SU(2)$. Since $\SO(3)$ is maximal-connected in $\SU(3)$, we may disregard this case by \ref{Hmax}. So assume $H_0=\SU(2)=\set{\diag(A,1)}$.

Then a proper closed subgroup $K$ with $K/H\approx S^l$ must be a conjugate of $\U(2)$, by \ref{grplist}. Now notice that the only conjugate of $\U(2)$ which contains $\SU(2) = \set{\diag(A,1)}$ is $\U(2) = \set{\diag(A,\det(\bar A))}$. So we can assume $K^\pm_0=\U(2) = \set{\diag(A,\det(\bar A))}$.

Recall that $H$ must be generated by a subgroup of $K^\pm_0 = \set{\diag(A,\det(\bar A))}$. Therefore $H = H_0\cdot \Z_n\subset K^\pm_0$ and $K^\pm$ are connected. We then get the following possible diagrams:
\begin{equation}\label{M6_8A}
   \SU(3) \,\, \supset \,\, \set{\diag(A,\det(\bar A))}, \,\, \set{\diag(A,\det(\bar A))} \,\, \supset \,\, \set{\diag(A,1)}\cdot \Z_n
\end{equation}
Conversely these all give simply connected manifolds by \ref{Hgens}. These actions make up the family $N^6_F$ of the appendix.

\subsubsection{Case 5: $G=\SU(3)\times S^1$}

Now $G=\SU(3)\times S^1$ and $H_0 = \U(2) \times 1 = \set{\diag(A,\det(\bar A))}\times 1$. To find the possible connected subgroups $K$ with $K/H$ a sphere, notice that such a $K\subset \SU(3)\times S^1$ would have $\proj_1(K)$ equal either $\U(2)$ or $\SU(3)$, by \ref{grplist}. Since $\rank(\U(2))=\rank(\SU(3))$, it follows that such a proper subgroup $K\supset H$ of higher dimension must be one of $\U(2)\times S^1$ or $\SU(3)\times 1$. The only such group with $K/H_0$ a sphere is $\U(2)\times S^1$. Hence we must have $K^\pm_0=\U(2)\times S^1$. Notice now that if $H$ were to have another component in $\U(2)\times S^1$ then $H$ would intersect $1\times S^1$ nontrivially, giving us a more effective action with the same groups. So we can assume this does not happen. Therefore we have the following possibility:
\begin{equation}\label{M6_8B}
   \SU(3)\times S^1 \,\, \supset \,\, \U(2)\times S^1, \,\, \U(2)\times S^1 \,\, \supset \,\, \U(2)\times 1
\end{equation}
which we see is simply connected by \ref{Hgens}. However, we see this is a product action.

\subsubsection{Case 6: $G=\Sp(2) \times S^1$}

Here $G = \Sp(2) \times S^1$ and $H_0 = \Sp(1)\Sp(1)\times 1$. To find the possibilities for connected groups $K$ with $K/H\approx S^l$ note that if $\proj_2(K)\subset S^1$ is nontrivial then $K = \Sp(1)\Sp(1)\times S^1$, by \ref{abfactor}. Otherwise $K\subset\Sp(2) \times 1$ and hence by \ref{grplist}, $K = \Sp(2) \times 1$. In either case $K/H\approx S^l$, in fact. Further by \ref{abfactor}, we can assume $K^- = \Sp(1)\Sp(1)\times S^1$ and all of $K^\pm$ and $H$ are connected. Therefore we have the following two possibilities:
\begin{equation}\label{M6_10A}
   \Sp(2) \times S^1 \,\, \supset \,\, \Sp(1)\Sp(1)\times S^1, \,\, \Sp(1)\Sp(1)\times S^1 \,\, \supset \,\, \Sp(1)\Sp(1)\times 1
\end{equation}
\begin{equation}\label{M6_10B}
   \Sp(2) \times S^1 \,\, \supset \,\, \Sp(1)\Sp(1)\times S^1, \,\, \Sp(2) \times 1 \,\, \supset \,\, \Sp(1)\Sp(1)\times 1
\end{equation}
both of which are simply connected by \ref{Hgens}. We easily see that the first action is a product action and the second action is a sum action.

\subsubsection{Case 7: $G = \Spin(6)$}

In this case we know from \ref{topdim}, that this gives a two-fixed-point action on a sphere.

%
%
%
%
%
%

\section{Classification in Dimension Seven}\label{7dim}

In this section we complete the classification in dimension 7. As in the previous section we keep the notation and conventions established in Section \ref{5dim}, this time for a 7 dimensional manifold manifold $M$. In this case, the next proposition gives us the possibilities for $G$ and $H_0$.


\begin{prop}
Table \ref{7dimgrps} list all the possibilities for $G$ and $H_0$, up to equivalence.
\end{prop}
\begin{proof}
We will first show that all the possibilities for $G$ are listed in \ref{7dimgrps}. Recall that $6\le\dim(G)\le 21$ by \ref{dimG} and that $\dim H = \dim G - 6$, since $\dim G/H= \dim M-1 =6$ in this case. A priori, by \ref{posgrps} we need to check all of the possibilities for $G$ of the form $(S^3)^m\times T^n$, $(\SU(3))^l \times (S^3)^m \times T^n$, $(\Sp(2))^k\times (S^3)^m \times T^n$, $\G_2 \times (S^3)^m \times T^n$, $\SU(4) \times (S^3)^m \times T^n$,  $\Sp(2) \times \SU(3) \times (S^3)^m \times T^n$, $\Sp(3)$ and $\Spin(7)$. Note that by \ref{abfactor} we can assume that $n\le2$ in all cases.

First suppose $G = (S^3)^m \times T^n$. By \ref{prods}, $3m + n - 6 = \dim(H)\le m$ which means $0\le 6 - 2m - n$ and so $m\le 3$ and if $m=3$ then $n=0$. Notice that all of these possibilities are listed in the table. Next if $G= (\SU(3))^l \times (S^3)^m \times T^n$ for $l>0$, then as before $8l+3m+n-6 = \dim(H) \le 4l+m$ or $0\le 6- 4l-2m-n$. Hence $l=1$, and: $m=1$ and $n=0$ or $m=0$ and $n\le 2$. All of these possibilities are listed in the table. Next suppose $G= (\Sp(2))^k \times (S^3)^m \times T^n$. Then we get $10k+3m+n-6= \dim(H)\le 6k + m$ or $0\le 6- 4k - 2m -n$. As before $k=1$, and: $m=1$ and $n=0$ or $m=0$ and $n\le 2$. However, if $G = \Sp(2)\times S^1$ then $\dim H= 5$ and by \ref{prods}, $H_0\subset \Sp(2)\times 1$. Then $\rank H \le \rank \Sp(2)=2$ and yet there are no compact 5-dimensional groups of rank 2 or less. So $\Sp(2)\times S^1$ is not a possibility for $G$. Next, if $G = \Sp(2) \times \SU(3) \times (S^3)^m \times T^n$ then $0\le -2 - 2m -n$, which is impossible. Now say $G = \G_2 \times (S^3)^m \times T^n$. We get $14 + 3m + n -6 \le 8 + m$ or $0\le -2m -n$ and hence $m=n=0$. Lastly, if $G = \SU(4) \times (S^3)^m \times T^n$ then $15 + 3m + n - 6 \le 10 + m$ or $0\le 1 - 2m - n$. Therefore $m=0$ and $n\le1$. Finally if $\dim(G)=21$ we know from \ref{topdim} that $G$ must be isomorphic to $\Spin(7)$ and in this case $H$ will be $\Spin(6)$.

Next we check that in the rest of the cases, we have listed all the possibilities for $H_0$. Again, we can assume that $H_0\subset G_1\times 1$ in the cases that $G=G_1\times T^m$. Then we use \ref{grplist} to find the possibilities for $H_0$. The only exceptional cases are 9 and 11, where $G=G_1\times S^3$. By \ref{prods}, $H_0\subset L\times S^1$ where $L$ is of dimension 4 or less in Case 9 and dimension 6 or less in Case 11. However, since $\dim H= \dim G-6$, we see that $H_0= L\times S^1$ where $L$ is of maximal dimension in each case. From \ref{grplist}, we see that $H_0$ must be one of the groups listed below.
\end{proof}

{\setlength{\tabcolsep}{0.40cm}
\renewcommand{\arraystretch}{1.6}
\stepcounter{equation}
\begin{table}[!h]
      \begin{center}
          \begin{tabular}{|c||c c|}
\hline
No.    &  $G$   & $H_0$  \\
\hline \hline
1   &    $S^3\times S^3$ & $\set{1}$   \\
\hline
2   &    $S^3\times S^3\times S^1$    &  $\set{(e^{ip\theta},e^{iq\theta})}\times 1$  \\
\hline
3   &    $S^3\times S^3\times T^2$    &   $T^2\times 1$  \\
\hline
4   &    $\SU(3)$    &   $T^2$  \\
\hline
5   &    $S^3\times S^3\times S^3$    &   $T^3$  \\
\hline
6   &    $\SU(3)\times S^1$    &   $\SU(2)\times 1$, $\SO(3)\times 1$  \\
\hline
7   &    $\SU(3)\times T^2$    &   $\U(2)\times 1$   \\
\hline
8   &    $\Sp(2)$    &   $\U(2)_{max}$, $\Sp(1)\SO(2)$   \\
\hline
9   &    $\SU(3)\times S^3$    &   $\U(2)\times S^1$  \\
\hline
10   &   $\Sp(2)\times T^2$    &   $\Sp(1)\Sp(1)\times 1$  \\
\hline
11   &   $\Sp(2)\times S^3$    &   $\Sp(1)\Sp(1)\times S^1$  \\
\hline
12   &   $\G_2$    &   \SU(3)  \\
\hline
13   &   $\SU(4)$    &   $\U(3)$  \\
\hline
14   &   $\SU(4)\times S^1$    &   $\Sp(2)\times 1$   \\
\hline
15   &   $\Spin(7)$    &   $\Spin(6)$  \\
\hline
          \end{tabular}
      \end{center}
      \vspace{0.1cm}
      \caption{Possibilities for $G$ and $H_0$, in the 7-dimensional case.}\label{7dimgrps}
\end{table}}

As in the previous sections we proceed to find all possible diagrams, by taking each case, one at a time.

\subsection{Cases 2 and 6}

Here we present Cases 2 and 6 which are similar since they both involve the same difficulty that arises in the case of $G=S^3\times S^1$ in dimension 5. In each case we will use the following lemma to deal with this difficulty.

\begin{lem}\label{7dim_exception}

Let $M$ be a simply connected cohomogeneity one manifold given by the group diagram $G\supset K^-,K^+\supset H$, with $G=G_1\times S^1$, $G_1$ simply connected, and $H_0=H_1\times 1$. Suppose further that there is a compact subgroup $L\subset G_1$ of the form $L=H_1\cdot\set{\beta(\theta)}$ where $\set{\beta(\theta)}$ is a circle group of $G$ parameterized once around by $\beta:[0,1]\to G_1$ and $\set{\beta(\theta)}\cap H_1=1$. Define $\delta:[0,1]\to G:t\mapsto (1,e^{2\pi it})$ be a loop once around $1\times S^1$. If $K^\pm_0\subset L\times S^1$ then the group diagram for $M$ has one of the following forms, all of which give simply connected manifolds:

\begin{eqnarray}\label{d:exception_ne}
   G_1\times S^1 \,\, \supset \,\, H_+\cdot \set{(\beta(m_-\theta),\delta(n_-\theta))}, \,\, H_-\cdot \set{(\beta(m_+\theta),\delta(n_+\theta))}  \,\, \supset \,\, H\\
   \nonumber \text{where $H=H_-\cdot H_+$, $K^-\ne K^+$, $\gcd(n_-,n_+,d)=1$}\\
   \nonumber \text{and $d$ is the index of $H\cap K^-_0 \cap K^+_0$ in $K^-_0 \cap K^+_0$},
\end{eqnarray}
\begin{eqnarray}\label{d:exception_e}
   G_1\times S^1 \,\, \supset \,\, \set{(\beta(m\theta),\delta(\theta))}\cdot H_0, \,\, \set{(\beta(m\theta),\delta(\theta))} \cdot H_0 \,\, \supset \,\, H_0\cdot \Z_n\\
   \nonumber \text{where $\Z_n\subset \set{(\beta(m\theta),\delta(\theta))}$}.
\end{eqnarray}

\end{lem}
\begin{proof}
It is clear, as in \ref{abfactor}, that $K^\pm/H$ must be circles and hence $K^\pm_0= H_0\cdot\set{(\beta(m_\pm\theta),\delta(n_\pm\theta))}$. From \ref{Hgens}, $H$ must have the form $H=H_-\cdot H_+$ for $H_\pm = K^\pm_0\cap H= H_0\cdot \Z_{k_\pm}$ for $\Z_{k_\pm}\subset \set{(\beta(m_\pm\theta),\delta(n_\pm\theta))}$. Then with the notation of \ref{Hgens}, we see that $\alpha_\pm$ can be taken as $\alpha_\pm(t) = (\beta(m_\pm t/k_\pm),\delta(n_\pm t/k_\pm))$. Then \ref{Hgens} says that $M$ is simply connected if and only if $\alpha_\pm$ generate $\pi_1(G/H_0)$. Since $\set{\beta(\theta)}\cap H_1=1$ we see that $\set{\beta(\theta)}$ injects onto a circle in $G/H_0$ which is contractible since $G_1$ is simply connected. We also see that $\delta$ generates $\pi_1(G/H_0)$ since $H_0\subset G_1\times 1$.

This brings us precisely to the situation we encountered in Case (B2) of Section \ref{dim5.4}. The argument given there shows that if $K^-_0=K^+_0$ we get the second diagram from the lemma, and in the case $K^-_0\ne K^+_0$ then $M$ is simply connected if and only if $\gcd(n_-,n_+,d)=1$ where $d$ is the index of $H/H_0\cap K^-_0/H_0 \cap K^+_0/H_0$ in $K^-_0/H_0 \cap K^+_0/H_0$. We can also write $d$ as the index of $H\cap K^-_0 \cap K^+_0$ in $K^-_0 \cap K^+_0$.
\end{proof}

We will now address Cases 2 and 6 individually, making use of the lemma when needed.

\subsubsection{Case 2: $G=S^3\times S^3\times S^1$}

Here $G=S^3\times S^3\times S^1$ and $H_0=\set{(e^{ip\theta},e^{iq\theta},1)}$. After an automorphism of $G$ we can assume that $p\ge q\ge 0$ and in particular $p\ne 0$.  We know from \ref{abfactor} that $K^-_0$, say, is a two torus. After conjugation we can assume that $K^-_0 = \set{(e^{ia_-\theta},e^{ib_-\theta},e^{ic_-\theta})}\cdot \set{(e^{ip\theta},e^{iq\theta},1)}$ even if $q=0$. From \ref{abfactor}, if $\proj_3(K^+_0)$ is nontrivial then $K^+_0$ is also a torus. Otherwise $K^+_0\subset S^3\times S^3\times 1$. Therefore from \ref{grplist}, we see that $K^+_0$ must be one of the following groups: $T^2$, $S^3\times 1\times 1$ if $q=0$, $\Delta S^3\times 1$ if $p=q=1$, or $S^3\times S^1\times 1$ if $q=1$ and allowing arbitrary $p$. We will now break this up into cases depending on what $K^+_0$ is.

\begin{case}[A] $\dim K^+ > 2$.

Then by \ref{fundgrp2}, $K^-$ is connected and $M$ is simply connected if and only if $G/K^-$ is. So we can assume $K^- = \set{(e^{ia\theta},e^{ib\theta},e^{i\theta})}\cdot \set{(e^{ip\theta},e^{iq\theta},1)}$, that is $c=1$. We also know that $H$ is a subgroup of $K^-$ of the form $\set{(e^{ip\theta},e^{iq\theta},1)}\cdot \Z_n$ for $\Z_n\subset \set{(e^{ia\theta},e^{ib\theta},e^{i\theta})}$, such that $\Z_n\cap K^+_0 = 1$, which is automatic, and $\Z_n\subset N(K^+_0)$.

\begin{case}[A1] $K^+_0=S^3\times 1\times 1$ and $q=0$.

Then $H_0=S^1\times 1\times 1$ and $K^-=\set{(1,e^{ib\theta},e^{i\theta})}\cdot \set{(e^{i\theta},1,1)}$. Then we see that the $\Z_n$ in $H$ can be arbitrary and we get the following family:
\begin{eqnarray}\label{M7_7A1}
   S^3\times S^3\times S^1 \,\, \supset \,\, \set{(e^{i\phi},e^{ib\theta},e^{i\theta})}, \,\, S^3\times 1\times 1\cdot \Z_n  \,\, \supset \,\, S^1\times 1\times 1 \cdot \Z_n\\
   \nonumber \Z_n\subset \set{(1,e^{ib\theta},e^{i\theta})}.
\end{eqnarray}
This is family $Q^7_C$ of the appendix.
\end{case}

\begin{case}[A2] $K^+_0=\Delta S^3\times 1$ and  $p=q=1$.

Here $H_0=\set{(e^{i\theta},e^{i\theta}, 1)}$ and we can take $K^-= \set{(1,e^{ib\theta},e^{i\theta})}\cdot \set{(e^{i\theta},e^{i\theta},1)}$ for a new $b$. Then for $\Z_n \subset \set{(1,e^{ib\theta},e^{i\theta})}$ to satisfy $\Z_n\subset N(K^+_0)$ simply means that $n|2b$. Then the further condition that $H\cap 1\times 1\times S^1=1$, for the action to be effective, means that $n$ is 1 or 2. Therefore we have the following diagrams in this case:
\begin{eqnarray}\label{M7_7A2}
   S^3\times S^3\times S^1 \,\, \supset \,\, \set{(e^{i\phi},e^{i\phi}e^{ib\theta},e^{i\theta})}, \,\, \Delta S^3\times 1\cdot \Z_n  \,\, \supset \,\, \Delta S^1\times 1\cdot \Z_n\\
   \nonumber \Z_n\subset \set{(1,e^{ib\theta},e^{i\theta})} \text{ where $n$ is 1 or 2}.
\end{eqnarray}
This family of actions gives the actions of type $Q^7_D$ from the appendix.
\end{case}

\begin{case}[A3] $K^+_0=S^3\times S^1\times 1$, $q=1$ and $p$ arbitrary.

Here $H_0=\set{(e^{ip\theta},e^{i\theta},1)}$ and we can take $K^-=\set{(e^{ia\theta},1,e^{i\theta})}\cdot \set{(e^{ip\theta},e^{i\theta},1)}$ for a new $a$. Then the $\Z_n\subset\set{(e^{ia\theta},1,e^{i\theta})}$ in $H$ automatically satisfies the condition $\Z_n\subset N(K^+_0)$. Hence we have the following diagrams:
\begin{eqnarray}\label{M7_7A3}
   S^3 \! \times \! S^3 \! \times \! S^1 \,\, \supset \,\, \set{(e^{ip\phi}e^{ia\theta},e^{i\phi},e^{i\theta})}, \,\, S^3\! \times \! S^1\! \times\! \Z_n  \,\, \supset \,\, \set{(e^{ip\phi},e^{i\phi},1)}\cdot \Z_n\\
   \nonumber \Z_n\subset \set{(e^{ia\theta},1,e^{i\theta})}.
\end{eqnarray}
This is the family $N^7_F$.
\end{case}
\end{case}

\begin{case}[B] $\dim K^+ =2$ so $K^+_0\approx T^2$.

Here $H_0=\set{(e^{ip\theta},e^{iq\theta},1)}$ again where we assume $p\ge q\ge0$ and $K^-_0 = \set{(e^{ia_-\theta},e^{ib_-\theta},e^{ic_-\theta})}\cdot \set{(e^{ip\theta},e^{iq\theta},1)}$. We now break this into two cases depending on whether or not $q$ is zero.

\begin{case}[B1] $q=0$

Here $H_0=S^1\times 1\times 1$ and so we know that $K^\pm_0= S^1\times \bar K^\pm_0$ for some groups $\bar K^\pm_0\subset S^3\times S^1$. Then from \ref{Hgens}, $H$ must have the form $S^1\times \bar H$ for a subgroup $\bar H$ generated by $\bar H\cap \bar K^-_0$ and $\bar H\cap \bar K^+_0$. Similarly, by \ref{Hgens}, the manifold $M$ will be simply connected if and only if the 5-manifold $\bar M$ given by the group diagram $S^3\times S^1\supset\bar K^-, \bar K^+ \supset \bar H$ is simply connected. So these actions are product actions with some simply connected 5 dimensional cohomogeneity one manifold.
\end{case}

\begin{case}[B2] $p,q\ne0$

Here we can take $K^\pm_0 = \set{(e^{ia_\pm\theta},e^{ib_\pm\theta},e^{ic_\pm\theta})}\cdot \set{(e^{ip\theta},e^{iq\theta},1)}$ although there is a more convenient way to write these groups in our case. Notice that for $p\mu-q\lambda=1$, we can write any element of the torus $T^2$ uniquely as $(e^{ip\theta},e^{iq\theta})(e^{i\lambda\phi},e^{i\mu\phi}) = (z^p,z^q)(w^\lambda,w^\mu)$. Then we can write
  $$K^\pm_0=\set{(z^p,z^q,1)(w^{m_\pm\lambda},w^{m_\pm\mu},w^{n_\pm})}$$
for some $m_\pm,n_\pm\in \Z$ with $\gcd(m_\pm,n_\pm)=1$. Then letting $\beta(t)=(e^{2\pi i\lambda t},e^{2\pi i\mu t})$ we see this satisfies the conditions of \ref{7dim_exception}. By that lemma, we have precisely the following two families of diagrams:
\begin{eqnarray}\label{M7_7B2a}         
   S^3\!\!\times\! S^3\!\!\times\! S^1\! \supset\!  \set{\!(z^pw^{\lambda m_-}\!,z^qw^{\mu m_-}\!,w^{n_-})\!}\! H,   \set{\!(z^pw^{\lambda m_+}\!,z^qw^{\mu m_+}\!,w^{n_+})\!}\! H \! \supset \! H\\
   \nonumber \text{where $H=H_-\cdot H_+$, $H_0=\set{(z^p,z^q,1)}$, $K^-\ne K^+$, $p\mu-q\lambda=1$,}\\
   \nonumber \text{$\gcd(n_-,n_+,d)=1$ where $d$ is the index of $H\cap K^-_0 \cap K^+_0$ in $K^-_0 \cap K^+_0$},
\end{eqnarray}  
\begin{eqnarray}\label{M7_7B2b}         
   S^3\!\!\times\! S^3\!\!\times\! S^1 \supset  \set{(z^pw^{\lambda m}\!,z^qw^{\mu m}\!,w)},   \set{(z^pw^{\lambda m}\!,z^qw^{\mu m}\!,w)}   \supset  H_0\cdot\Z_n \\
   \nonumber \text{where $H_0=\set{(z^p,z^q,1)}$, $p\mu-q\lambda=1$ and $\Z_n\subset \set{(w^{\lambda m}\!,w^{\mu m}\!,w)}$}.
\end{eqnarray}  
These two families are $N^7_E$ and $N^7_D$, respectively.
\end{case}
\end{case}

\subsubsection{Case 6: $G=\SU(3)\times S^1$}

Here $G=\SU(3)\times S^1$ and $H_0$ is either $\SU(2)\times 1$ or $\SO(3)\times 1$. First, if $H_0= \SO(3)\times 1$ then $H_1=\SO(3)$ is maximal in $\SU(3)$ and so by \ref{abfactor}, $H$, $K^-$ and $K^+$ are all connected, $K^-$, say, is $\SO(3)\times S^1$ and $K^+$ is either $\SO(3)\times S^1$ or $\SU(3)\times 1$. Since $\SU(3)/\SO(3)$ is not a sphere we see we have only one possible diagram:
\begin{equation}\label{M7_9a}
   \SU(3)\times S^1 \,\, \supset \,\, \SO(3)\times S^1, \,\, \SO(3)\times S^1 \,\, \supset \,\, \SO(3)\times 1,
\end{equation}
which comes from a product action.

For the other case assume $H_0 = \SU(2)\times 1$ where $\SU(2)=\SU(1)\SU(2)$ is the lower right block. Notice from \ref{grplist}, that $\proj_1(K^\pm_0)$ is either $\SU(2)$, $\U(2)$ or $\SU(3)$. We know, as in \ref{abfactor}, that if $\proj_2(K^\pm_0)$ is nontrivial then $K^\pm_0= H_0\cdot S^1$ and hence has the from $\set{(\beta(m_\pm\theta),e^{in_\pm\theta})}\cdot H_0$ where $\beta(\theta)= \diag(e^{-i\theta},e^{i\theta},1)\in \SU(3)$. In fact $K^-_0$ must have this form, so assume $K^-_0=\set{(\beta(m_-\theta),e^{in_-\theta})}\cdot H_0$. The other possibility for $K^+_0$ is $\SU(3)\times 1$ which does give $K^+_0/H_0\approx S^5$.

First suppose $K^+_0 =\SU(3)\times 1$. Then from \ref{fundgrp2}, $K^-$ is connected and $\pi_1(M)\approx \pi_1(G/K^-)$. It then follows that $n_-=1$ so $K^-=\set{(\beta(m\theta),e^{i\theta})}\cdot H_0$ in this case. From \ref{Hgens}, $H=H_0\cdot \Z_n$ for $\Z_n\subset \set{(\beta(m\theta),e^{i\theta})}$. The condition that $H\cap \SU(3)\times 1=1$ means that $\gcd(m,n)=1$. Therefore we get the following family of diagrams in this case:
\begin{eqnarray}\label{M7_9b}
   \SU(3)\times S^1 \,\, \supset \,\, \set{(\beta(m\theta),e^{i\theta})}\cdot H_0, \,\, \SU(3)\times \Z_n  \,\, \supset \,\, H_0\cdot \Z_n\\
   \nonumber H_0 = \SU(1)\SU(2)\times 1,\,\, \Z_n\subset \set{(\beta(m\theta),e^{i\theta})},\\
   \nonumber \beta(\theta)= \diag(e^{-i\theta},e^{i\theta},1),\,\, \gcd(m,n)=1.
\end{eqnarray}
This is family $Q^7_G$ of the appendix.

Next assume $K^\pm_0=\set{(\beta(m_\pm\theta),e^{in_\pm\theta})}\cdot H_0$. Notice that $\set{\beta(\theta)}\cap H_0=1$ and hence this situation satisfies the hypotheses of \ref{7dim_exception}, for $L=\U(2)$. Then, by that lemma, we have precisely the following two families of diagrams:
\begin{eqnarray}\label{M7_9c}
   \SU(3)\times S^1 \,\, \supset \,\, \set{(\beta(m_-\theta),e^{in_-\theta})}\cdot H, \,\, \set{(\beta(m_+\theta),e^{in_+\theta})}\cdot H  \,\, \supset \,\, H\\
   \nonumber H_0 = \SU(1)\SU(2)\times 1,\,\, H=H_-\cdot H_+,\,\, K^-\ne K^+, \\
   \nonumber \beta(\theta)= \diag(e^{-i\theta},e^{i\theta},1),\,\, \gcd(n_-,n_+,d)=1\\
   \nonumber \text{ where $d$ is the index of $H\cap K^-_0 \cap K^+_0$ in $K^-_0 \cap K^+_0$}
\end{eqnarray}
\begin{eqnarray}\label{M7_9d}
   \SU(3)\times S^1 \,\, \supset \,\, \set{(\beta(m\theta),e^{i\theta})}\cdot H_0, \,\, \set{(\beta(m\theta),e^{i\theta})}\cdot H_0  \,\, \supset \,\, H_0\cdot \Z_n\\
   \nonumber H_0 = \SU(1)\SU(2)\times 1,\,\, \Z_n\subset \set{(\beta(m\theta),e^{i\theta})},\\
   \nonumber \beta(\theta)= \diag(e^{-i\theta},e^{i\theta},1).
\end{eqnarray}
The first of these families is $N^7_H$ and the second is $Q^7_F$.

\subsection{The remaining cases}

Here we address the rest of the cases from \ref{7dimgrps}.

\subsubsection{Cases 3, 7, 10 and 14}

Notice that in Cases 3, 7 and 10 we have $G=G_{ss}\times T^2$ where $G_{ss}$ is semisimple and $\rank(H)=\rank(G_{ss})$. Therefore by \ref{abfactor}, the resulting actions must all be product actions.

In Case 14 we see from \ref{grplist} that $H_1=\Sp(2)$ is maximal among connected subgroups in $\SU(4)$. Therefore, \ref{abfactor} says $K^-$, $K^+$ and $H$ are connected. Further, we can assume $K^-=H_1\times S^1$ and that $K^+$ is either $H_1\times S^1$ or has the form $K_1\times 1$, for $K_1/H_1\approx S^l$. If $K^+=K_1\times 1$ then by \ref{grplist}, $K_1$ would have to be $\SU(4)$ and in this case we do have $\SU(4)/\Sp(2) \approx S^5$. Therefore we have the two following possibilities, both of which give simply connected manifolds:
\begin{equation}\label{M7_15a}
   \SU(4)\times S^1\,\, \supset \,\, \Sp(2)\times S^1, \,\, \Sp(2)\times S^1 \,\, \supset \,\, \Sp(2)\times 1,
\end{equation}
\begin{equation}\label{M7_15a}
   \SU(4)\times S^1\,\, \supset \,\, \Sp(2)\times S^1, \,\, \SU(4)\times 1 \,\, \supset \,\, \Sp(2)\times 1.
\end{equation}
We notice that the first is a product action and the second is a sum action.

\subsubsection{Cases 9 and 11}

In both cases $G=G_1\times S^3$ and $H_0=H_1\times S^1$ where $H_1$ is maximal among connected subgroups of $G_1$. Then $\proj_1(K^\pm_0)$ are either $H_1$ or $G_1$ and $\proj_2(K^\pm_0)$ are either $S^1$ or $S^3$. It is also clear that if $\proj_2(K^\pm_0)= S^3$ then $K^\pm_0 \supset 1\times S^3$ and so if $\proj_1(K^\pm_0)= G_1$ then $K^\pm_0 \supset G_1\times 1$ as well. Therefore the proper subgroups $K^\pm_0$ must each be either $G_1\times S^1$ or $H_1\times S^3$. Note that $H_1\times S^3/ H_1\times S^1$ is always a sphere. In Case 9, $G_1\times S^1/H_1\times S^1\approx \SU(3)/\U(2)\approx \CP^2$ so this is not a possibility for $K^\pm$ but in Case 11, $G_1\times S^1/H_1\times S^1$ is a sphere. Notice that in all cases $l_\pm>1$ so $H$, $K^-$ and $K^+$ must all be connected by \ref{fundgrp2}. Therefore we have the following possible diagrams:
\begin{equation}\label{M7_15a}
   \SU(3)\times S^3\,\, \supset \,\, \U(2)\times S^3, \,\, \U(2)\times S^3 \,\, \supset \,\, \U(2)\times S^1,
\end{equation}
\begin{equation}\label{M7_15a}
   \Sp(2)\times S^3\,\, \supset \,\, \Sp(1)\Sp(1)\times S^3, \,\, \Sp(1)\Sp(1)\times S^3 \,\, \supset \,\, \Sp(1)\Sp(1)\times S^1,
\end{equation}
\begin{equation}\label{M7_15a}
   \Sp(2)\times S^3\,\, \supset \,\, \Sp(1)\Sp(1)\times S^3, \,\, \Sp(2)\times S^1 \,\, \supset \,\, \Sp(1)\Sp(1)\times S^1,
\end{equation}
\begin{equation}\label{M7_15a}
   \Sp(2)\times S^3\,\, \supset \,\, \Sp(2)\times S^1, \,\, \Sp(2)\times S^1 \,\, \supset \,\, \Sp(1)\Sp(1)\times S^1,
\end{equation}
all of which are simply connected by \ref{Hgens}. The third is a sum action and the remaining three actions are product actions.

\subsubsection{Cases 12, 13 and 15}

In each of these cases, $H_0$ is maximal in $G$ among connected subgroups. Therefore, \ref{Hmax} gives a full description of these types of actions. \ref{topdim} also deals with Case 15 separately.

\subsubsection{Case 1: $G=S^3\times S^3$}

Here $G=S^3\times S^3$ and $H$ is discrete. Since $H$ is discrete it follows that for $K^\pm/H$ to be spheres, $K^\pm_0$ must themselves be covers of spheres. From \ref{grplist} we see that $K^\pm_0$ must be one of the following: $\set{(e^{x_\pm p_\pm\theta},e^{y_\pm q_\pm\theta})}$ for $x_\pm,y_\pm \in \Im(\H)$, $S^3\times 1$, $1\times S^3$ or $\Delta_{g_0} S^3 = \set{(g,g_0gg_0^{-1})}$ for $g_0 \in S^3$. We break this into cases depending on what $K^\pm$ are.

\begin{case}[A] $K^-_0\approx S^3$ and $K^+_0\approx S^3$.

In this case we know from \ref{fundgrp2} that $H$, $K^-$ and $K^+$ must all be connected in this case. Hence $N(H)_0=S^3\times S^3$ and we can conjugate $K^-$ and $K^+$ by anything in $S^3\times S^3$ without changing the manifold, by \ref{normalizer}. In particular if $K^\pm=\Delta_{g_0} S^3$ then we can assume $g_0=1$. Therefore we get the following possible groups diagrams up to automorphism of $G$, all of which are clearly simply connected by \ref{Hgens}:

\begin{equation}\label{M7_6A1}
   S^3\times S^3 \,\, \supset \,\, S^3\times 1, \,\, S^3\times 1  \,\, \supset \,\,  1
\end{equation}
\begin{equation}\label{M7_6A2}
   S^3\times S^3 \,\, \supset \,\, S^3 \times 1, \,\, 1 \times S^3  \,\, \supset \,\,  1
\end{equation}
\begin{equation}\label{M7_6A3}
   S^3\times S^3 \,\, \supset \,\, S^3\times 1, \,\, \Delta S^3  \,\, \supset \,\,  1
\end{equation}
\begin{equation}\label{M7_6A4}
   S^3\times S^3 \,\, \supset \,\, \Delta S^3, \,\, \Delta S^3  \,\, \supset \,\,  1
\end{equation}
The first of these actions is a product action and the second is a sum action. The last two are actions $Q^7_A$ and $Q^7_B$, respectively.
\end{case}  

\begin{case}[B] $K^-_0\approx S^1$ and $K^+_0\approx S^3$.

From \ref{Hgens}, we know that $K^-$ is connected and $H=\Z_n\subset K^-$ such that $H\cap K^+_0 = 1$. After conjugation of $G$ we can assume that $K^- = \set{(e^{ip\theta},e^{iq\theta})}$. If $K^+_0=S^3\times 1$ then the condition $H\cap K^+_0 = 1$ means that $n$ and $q$ are relatively prime. Therefore we have the following family of diagrams:
\begin{eqnarray}\label{M7_6B1}
   S^3\times S^3 \,\, \supset \,\, \set{(e^{ip\theta},e^{iq\theta})}, \,\, S^3\times \Z_n \,\, \supset \,\,  Z_n\\
   \nonumber\text{  where  } (q,n)=1
\end{eqnarray}
which all give simply connected manifolds by \ref{Hgens}. These actions are the actions of type $N^7_C$.

Next suppose that $K^+_0 = \Delta_{g_0} S^3$ for some $g_0\in S^3$. Notice that $N(K^+_0)= \set{(\pm g, g_0gg_0^{-1})}$ and since $L\subset N(L_0)$ for every subgroup $L$ it follows that $K^+$ can have at most two components and hence $H$ can have at most two elements. In particular this means that $H$ is normal in $G$ and hence by \ref{normalizer} we can conjugate $K^+$ by $(1,g_0^{-1})$ without changing the resulting manifold. Lastly, if $n=2$ the condition that $H\cap K^+_0=1$ means that $p$ and $q$ are not both odd and not both even since $(p,q)=1$. Without loss of generality we can assume that $p$ is even and $p$ is odd. Therefore we have the following family of diagrams, all of which are simply connected by \ref{Hgens}:
\begin{eqnarray}\label{M7_6B1}
   S^3\times S^3 \,\, \supset \,\, \set{(e^{ip\theta},e^{iq\theta})}, \,\, \Delta S^3 \cdot \Z_n \,\, \supset \,\,  Z_n\\
   \nonumber\text{  where $n$ is 1 or 2, $p$ even, $(p,q)=1$.}
\end{eqnarray}
This is family $P^7_D$.
\end{case}  

\begin{case}[C] $K^-_0\approx S^1$ and $K^+_0\approx S^1$.

Here we have $K^\pm_0= \set{(e^{x_\pm p_\pm\theta},e^{y_\pm q_\pm\theta})}$. To address this case we will break it up into further cases depending on how big the group generated by $K^-_0$ and $K^-_0$ is.

\begin{case}[C1] $K^-_0$ and $K^+_0$ are both contained in some torus.

After conjugation we can assume that $K^\pm_0 = \set{(e^{ip_\pm\theta},e^{iq_\pm\theta})}$. By \ref{Hgens}, $H = H_-\cdot H_+$ where $H_\pm = \Z_{n_\pm}\subset K^\pm_0$ and conversely by \ref{Hgens}, such groups will always give simply connected manifolds. Therefore we have the following possibilities:
\begin{eqnarray}\label{M7_6C1}
   S^3\times S^3 \,\, \supset \,\, \set{(e^{ip_-\theta},e^{iq_-\theta})}\cdot H_+, \,\, \set{(e^{ip_+\theta},e^{iq_+\theta})} \cdot H_- \,\, \supset \,\,  H_-\cdot H_+\\
   \nonumber H_\pm = \Z_{n_\pm}\subset K^\pm_0
\end{eqnarray}
which make up family $N^7_A$.
\end{case}  

\begin{case}[C2] $K^-_0$ and $K^+_0$ are both contained in $S^3\times 1$.

In this case it follows from \ref{Hgens} that $H$, $K^-$ and $K^+$ are all contained in $S^3\times 1$. It also follows from \ref{Hgens} that $M^7$, given by the diagram $G\supset K^-, K^+\supset H$, will be simply connected if and only if the manifold $N^4$ given by the diagram $S^3\times 1 \supset K^-, K^+\supset H$ is simply connected. Therefore this gives a product action.
\end{case}  

\begin{case}[C3] $K^-_0$ and $K^+_0$ are both contained in $S^3\times S^1$ but not in $T^2$ or $S^3\times 1$.

It follows from \ref{Hgens} that $H$, $K^-$ and $K^+$ must all be contained in $S^3\times S^1$ in this case. Notice further that if both $p_-q_-=0$ and $p_+q_+=0$ then we would be back in one of the previous cases. So after conjugation of $G$ and switching of $-$ and $+$, we can assume that $K^-_0 = \set{(e^{ip_-\theta},e^{iq_-\theta})}$, where $p_-q_-\ne 0$. For $K^+_0$ we can assume that $y_+=i$ and denote $x_+=x$. It also follows that $p_+\ne 0$ and $x\ne\pm i$ since otherwise we would be in a previous case again.

Notice that $N(K^-_0) = \set{(e^{i\theta},e^{i\phi})} \cup \set{(je^{i\theta},je^{i\phi})}$ and $K^-\subset S^3\times S^1$ and hence $K^-\subset \set{(e^{i\theta},e^{i\phi})}$. Similarly if $q_+\ne 0$ then $N(K^+_0) = \set{(e^{x\theta},e^{i\phi})} \cup \set{(we^{x\theta},je^{i\phi})}$ for $w\in x^\bot \cap \Im S^3$. Therefore $K^+ \subset \set{(e^{x\theta},e^{i\phi})}$ in this case as well. However $H$ would then be a subset of the intersection of these two sets, $H\subset \set{(\pm 1, e^{i\phi})}$, and $N(H)_0$ would contain $S^3\times 1$. We would then be able to conjugate $K^+$ into the set $\set{(e^{i\theta},e^{i\phi})}$ without changing the resulting manifold, by \ref{normalizer}. This would put us back into Case (C1), so we can assume that $q_+ = 0$ and $K^+_0 = \set{(e^{x\theta},1)}$.

Therefore $N(K^+_0) = \big( \set{e^{x\theta}} \cup \set{we^{x\theta}} \big) \times S^3$. Again we see that for $N(K^-_0) \cap N(K^+_0) \nsubseteq \set{(\pm 1, e^{i\phi})}$ we need $x\bot i$. So after conjugation we can assume $K^+_0 = \set{(e^{j\theta},1)}$. Then $H\subset \set{\pm 1,\pm i} \times S^1$. By \ref{Hgens}, $H=H_-\cdot H_+$ for $H_\pm = \Z_{n_\pm} \subset K^\pm_0$. We see then that $n_+$ is $1$ or $2$ and the conditions that $H\subset \set{\pm 1,\pm i} \times S^1$ but $H\nsubseteq \set{\pm 1} \times S^1$ mean that $4|n_-$ and $p_- \equiv \pm n_-/4 \mod n_-$. Conversely we see we get the following possible diagrams:
\begin{eqnarray}\label{M7_6C3}
   S^3\times S^3 \,\, \supset \,\, \set{(e^{ip\theta},e^{iq\theta})}\cdot H_+, \,\, \set{(e^{j\theta},1)}\cdot H_- \,\, \supset \,\,  H_-\cdot H_+\\
   \nonumber \text{where $H_\pm = \Z_{n_\pm} \subset K^\pm_0$, $n_+\le 2$, $4|n_-$ and $p_- \equiv \pm \frac{n_-}{4} \mod n_-$},
\end{eqnarray}
all of which give simply connected manifolds by \ref{Hgens}. These actions make up the family $N^7_B$.
\end{case}  

\begin{case}[C4] $K^-_0$ and $K^+_0$ are not both contained in $S^3 \times S^1$ or $S^1 \times S^3$.

As in the previous case, we can assume here that both $p_-q_-\ne0$ and $p_+q_+\ne 0$ and after conjugation $K^-_0 = \set{(e^{ip_-\theta},e^{iq_-\theta})}$ and $K^+_0= \set{(e^{xp_+\theta},e^{yq_+\theta})}$. Then if $u\in x^\bot \cap \Im S^3$ and $w\in y^\bot \cap \Im S^3$ then we have $N(K^-_0) = \set{(e^{i\theta},e^{i\phi})}\cup \set{(je^{i\theta},je^{i\phi})}$ and $N(K^+_0) = \set{(e^{x\theta},e^{y\phi})}\cup \set{(ue^{x\theta},ve^{y\phi})}$ and $H \subset N(K^-_0) \cap N(K^+_0)$.

We now claim that we can assume $x$ and $i$ are perpendicular for suppose they are not perpendicular. Then if we denote the two elements in $i^\bot \cap x^\bot \cap \Im S^3$ by $\pm w$, we would have $H \subset \set{\pm 1, \pm w}\times S^3$. Notice that conjugation by $(e^{w\alpha}, 1)$ fixes $\set{\pm 1, \pm w}\times S^3$ pointwise and hence $(e^{w\alpha}, 1) \in N(H)_0$ for all $\alpha\in \R$. Therefore, by \ref{normalizer} we can conjugate $K^+$ by $(e^{w\alpha}, 1)$ without changing the resulting manifold. Since $w\bot \set{x,i}$, conjugation by $(e^{w\alpha}, 1)$ fixes the $1w$-space and rotates the $ix$-space by $2\alpha$. So for the right choice of $\alpha$ we can rotate $x$ into $i$. Therefore we could assume that $K^+_0= \set{(e^{ip_+\theta},e^{yq_+\theta})}$, bringing us back to an earlier case. Hence we can assume that $x\bot i$ and similarly $y\bot i$. Then after conjugation of $G$ we can take $K^+_0 = \set{(e^{jp_+\theta},e^{jq_+\theta})}$, without affecting $K^-$.

Then the condition $H \subset N(K^+_0)\cap N(K^-_0)$ becomes
\begin{eqnarray}
  H & \subset & \set{\pm1}\times \set{\pm1} \cup \set{\pm i}\times \set{\pm i} \cup\set{\pm j}\times \set{\pm j} \cup\set{\pm k}\times \set{\pm k}\nonumber\\
  & = & \Delta Q \cup \Delta_- Q\nonumber
\end{eqnarray}
where $Q=\set{\pm1,\pm i, \pm j, \pm k}$ and $\Delta_- Q = \set{\pm(1,-1),\pm (i,-i), \pm (j,-j), \pm (k,-k)}$. In particular, if $(h_1,h_2)\in H$ then $h_1= \pm h_2$.

We also know from \ref{Hgens} that $H$ is generated by $H\cap K^-_0 =: H_-$ and $H\cap K^+_0 =: H_+$, where $H_\pm$ are both cyclic subgroups of the circles $K^\pm_0$. Let $h_\pm = (h_1^\pm,h_2^\pm)$ be generators of $H_\pm$, so that $h_-$ and $h_+$ generate $H$. Notice that if both $h_\pm$ have order 1 or 2 then $H$ would be contained in $\set{\pm 1}\times \set{\pm 1}$ and we would be back in a previous case, as before. So assume that $h_-$ has order 4 and after conjugation of $G$ we can assume that $h_-=(i,i)$. The condition that $h_-\in K^-_0$ means that $p_-,q_-\equiv \pm1 \mod 4$, however, after switching the sign of both $p_-$ and $q_-$ we can assume that $p_-,q_-\equiv 1 \mod 4$.

We will now break our study into further cases depending on the order of $h_+$, which is either 1, 2 or 4.

\begin{case}[C4a] $h_+\in \langle (i,i) \rangle$.

Then $H= \langle (i,i) \rangle$. Hence we get the following family of diagrams:
\begin{eqnarray}\label{M7_6C4a}
   S^3\times S^3 \,\, \supset \,\, \set{(e^{ip_-\theta},e^{iq_-\theta})}, \,\, \set{(e^{jp_+\theta},e^{jq_+\theta})}\cdot H \,\, \supset \,\,  \langle (i,i) \rangle\\
   \nonumber \text{where $p_-,q_-\equiv 1 \mod 4$}
\end{eqnarray}
which is the family $P^7_A$.
\end{case} 

\begin{case}[C4b] $\#(h_+)=2$ but $h_+ \notin \langle (i,i) \rangle$.

It follows that $h_+$ must be $(1,-1)$ or $(-1,1)$ and after switching the factors of $G=S^3\times S^3$ we can assume that $h_+ = (1,-1)$. The condition that $h_+\in K^+_0$ means that $p_+$ is even. Therefore we have the following family of possibilities:
\begin{eqnarray}\label{M7_6C4b}
   S^3\times S^3 \,\, \supset \,\, \set{(e^{ip_-\theta},e^{iq_-\theta})} \cdot H, \,\, \set{(e^{jp_+\theta},e^{jq_+\theta})}\cdot H \,\, \supset \,\,  \langle (i,i), (1,-1) \rangle\\
   \nonumber \text{where $p_-,q_-\equiv 1 \mod 4$, $p_+$ even}.
\end{eqnarray}
This is the family $P^7_B$ of the appendix.
\end{case} 

\begin{case}[C4c] $\#(h_+)=4$.

In this last case, $h_+$ must be one of $(j,j)$, $(j,-j)$, $(-j,j)$ or $(-j,-j)$. However, after conjugation of $G$ by $(\pm i, \pm i)$ we can assume that $h_+= (j,j)$. As before, the condition that $h_+\in K^+_0$ means that $p_+,q_+\equiv \pm1 \mod 4$ but we can assume that $p_+,q_+\equiv 1 \mod 4$, after a change of signs on $p_+$ and $q_+$. Then $H= \Delta Q$ and we have the following possibilities:
\begin{eqnarray}\label{M7_6C4c}
   S^3\times S^3 \,\, \supset \,\, \set{(e^{ip_-\theta},e^{iq_-\theta})} \cdot H, \,\, \set{(e^{jp_+\theta},e^{jq_+\theta})}\cdot H \,\, \supset \,\,  \Delta Q\\
   \nonumber \text{where $p_\pm,q_\pm\equiv 1 \mod 4$}.
\end{eqnarray}
This family of actions is family $P^7_C$ of the appendix.

\end{case} 
We see from \ref{Hgens}, that all of the diagrams above do give simply connected manifolds.
\end{case}  
\end{case}  

\subsubsection{Case 4: $G=\SU(3)$}

In this case, $G=\SU(3)$ and $H_0=T^2$. From \ref{grplist}, the proper subgroups $K^\pm_0$ must both be $\U(2)$ up to conjugacy. It then follows from \ref{fundgrp2} that $H$, $K^-$ and $K^+$ are all connected. Now fix $H = \diag(\SU(3)) \approx T^2$. If $K^\pm$ contains this $T^2$ then it must be a conjugate of $\U(2)$ by an element of the Weyl group $W=N(T^2)/T^2$. We see that there are precisely three such conjugates of $\U(2)$ and they are permuted by the elements of $W$. Therefore, there are two possibilities for the pair $K^-,K^+$ up to conjugacy of $G$: $\S(\U(1)\U(2)),\S(\U(1)\U(2))$ or $\S(\U(1)\U(2)),\S(\U(2)\U(1))$. This gives us precisely the following two simply connected diagrams:
\begin{equation}\label{M7_8a}
   \SU(3) \,\, \supset \,\, \S(\U(1)\U(2)), \,\, \S(\U(1)\U(2)) \,\, \supset \,\, T^2,
\end{equation}
\begin{equation}\label{M7_8b}
   \SU(3) \,\, \supset \,\, \S(\U(1)\U(2)), \,\, \S(\U(2)\U(1)) \,\, \supset \,\, T^2.
\end{equation}
The first is action $N^7_G$ and the second is $Q^7_E$.

\subsubsection{Case 5: $G=S^3\times S^3\times S^3$}

Now $G=S^3\times S^3\times S^3$ and $H_0=T^3$. It is clear that if $\proj_1(K^\pm_0)\ne S^1$ then $K^\pm_0 \supset S^3\times 1\times 1$ and similarly for the other factors. Hence each $K^\pm_0$ will be a product of $S^3$ factors and $S^1$ factors. Further, it is clear that for $K^\pm_0/H$ to be a sphere we need $K^\pm_0$ to be one of $S^3\times S^1\times S^1$, $S^1\times S^3\times S^1$ or $S^1\times S^1\times S^3$.  Then by \ref{fundgrp2}, all of $H$, $K^-$ and $K^+$ must be connected. Putting this together we see we have the following possible simply connected diagrams, up to $G$-automorphism:
\begin{equation}\label{M7_9Aa}
   S^3\times S^3\times S^3 \,\, \supset \,\, S^3\times S^1\times S^1, \,\, S^3\times S^1\times S^1 \,\, \supset \,\, S^1\times S^1\times S^1,
\end{equation}
\begin{equation}\label{M7_9Ab}
   S^3\times S^3\times S^3 \,\, \supset \,\, S^3\times S^1\times S^1, \,\, S^1\times S^3\times S^1 \,\, \supset \,\, S^1\times S^1\times S^1.
\end{equation}
It is clear that both of these are product actions.

\subsubsection{Case 8: $G=\Sp(2)$}

Here $G=\Sp(2)$ and $H_0$ is either $\U(2)_{max}=\set{\diag(zg,\bar z g)}$ or $\Sp(1)\SO(2)$. Since $\U(2)_{max}$ is maximal among connected subgroups, and $\Sp(2)/\U(2)_{max}$ is not a sphere, we see this is not a possibility for $H_0$. So assume $H_0=\Sp(1)\SO(2)$. Then from \ref{grplist}, we see the proper subgroups $K^\pm_0$ must be conjugates of $\Sp(1)\Sp(1)$. Since the only conjugate of $\Sp(1)\Sp(1)$ which contains $\Sp(1)\SO(2)$ is the usual $\Sp(1)\Sp(1)$ we see $K^\pm_0= \Sp(1)\Sp(1)$. Then by \ref{fundgrp2}, $H$, $K^-$ and $K^+$ must all be connected. Therefore we get the one possible diagram:
\begin{equation}\label{M7_10}
    \Sp(2)\,\, \supset \,\, \Sp(1)\Sp(1), \,\, \Sp(1)\Sp(1) \,\, \supset \,\, \Sp(1)\SO(2).
\end{equation}
This last action is action $N^7_I$ of the appendix.

%
%
%
%
%
%



\section{Identifying some actions}\label{identify}

Here we will study each action in Tables \ref{5dimlist} to \ref{7dimlist2} and find we can identify many of these actions.

\subsection{Isometric actions on symmetric spaces}\label{symspactions}

In this section we will list all isometric cohomogeneity one actions on compact simply connected symmetric spaces of dimension seven or less. In 1971, Hsiang and Lawson classified cohomogeneity one actions on symmetric spheres in \cite{HL} (see \cite{Straume} for correction) and later Uchida did the same for complex projective spaces in \cite{Uchida}. Then in 2001, Kollross generalized these results to a classification of cohomogeneity one actions on irreducible symmetric spaces of compact type in \cite{Kollross}.

These classifications found that the only maximal isometric cohomogeneity one actions on compact simply connected irreducible symmetric spaces of dimension 7 or less are the following, up to equivalence: the sum actions of $\SO(k_1)\times \SO(k_2)$ on $S^{k_1+k_2-1}\subset \R^{k_1}\times \R^{k_2}$, for $k_i\ge1$; the tensor actions of $\SO(k)\times \SO(2)$ on $S^{2k-1}\subset \R^{k\times 2}$ via $(A,B)\star M = AMB^{-1}$, for $k=3,4$; the exceptional sphere actions of $\SO(3)$ on $S^4$ as in Section \ref{dim4class}, $\SU(3)$ on $S^7\subset \fsu(3)$ via $\Ad$, and $\SO(4)$ on $S^7$ via the isotropy representation of $\G_2/\SO(4)$; the linear actions on complex projective spaces of $\U(n)$ or $\SO(n+1)$ on $\CP^n=\SU(n+1)/\U(n)$, or $\S(\U(2)\times \U(2))$ on $\CP^3=\SU(4)/\U(3)$; and the two remaining symmetric space actions of $\U(2)$ on $\SU(3)/\SO(3)$ and $\SO(4)$ on $\SO(5)/\SO(2)\SO(3)$.

There are several other actions which do not appear in this list since they are equivalent to actions that do appear. For example, the complex tensor action of $\SU(2)\times \U(2)$ on $S^7\subset \C^{2\times 2}$ via $(A,B)\star M = AMB^{-1}$ is equivalent to the real tensor action of $\SO(4)\times \SO(2)$ on $S^{7}\subset \R^{4\times 2}$ mentioned above.

For each action of $G$ on $M$ mentioned above, it is not difficult so see which subgroups of $G$ act on $M$ with the same orbits. Many of these actions are simply sum actions or fixed-point actions as described in Section \ref{special_types}. As we have already fully examined these special cases we will not address them again here. In Table \ref{nonred_irredsymspactions} we list the remaining nonreducible cohomogeneity one actions on irreducible symmetric spaces. 

{\setlength{\tabcolsep}{0.40cm}
\renewcommand{\arraystretch}{1.6}
\stepcounter{equation}
\begin{table}[!h]
\begin{center}
\begin{tabular}{|c|c|}
\hline
$Q^5_B$ (part) & $\SO(3)\times \SO(2)$ on $S^5\subset \R^{3\times 2}$ via $(A,B)\star M = AMB^{-1}$\\
\hline
$P^5$ (part) & $\U(2)$ on $\SU(3)/\SO(3)$\\
\hline
$Q^5_C$ (all) & $S^3\times S^1$ on $S^5\subset \H\times \C$ via $(g,z)\star (p,w)=(gp\bar{z}^n,z^mw)$\\
\hline
$Q^6_A$ & $\SO(4)$ on $\CP^3=\SU(4)/\U(3)$\\
\cline{2-2}
 (all)  & $\SO(4)$ on $\SO(5)/\SO(2)\SO(3)$\\
\hline
$Q^6_C$ (all) & $S^3\times S^3$ on $S^6\subset \H\times \Im(\H)$ via $(g_1,g_2)\star (p,q)=(g_1pg_2^{-1},g_2qg_2^{-1})$\\
\hline
$Q^6_D$ (all) & $\SU(2)\times \SU(2)$ on $\CP^3=\SU(4)/\U(3)$\\
\hline
$Q^7_A$ (all) & $S^3\times S^3$ on $S^7\subset \H\times \H$ via $(g_1,g_2)\star (p,q)=(g_1pg_2^{-1},g_2q)$\\
\hline
$Q^7_D$ (part) & $\SO(4)\times \SO(2)$ on $S^7\subset \R^{4\times 2}$ via $(A,B)\star M = AMB^{-1}$\\
\hline
$Q^7_E$ (all) & $\SU(3)$ on $S^7\subset \fsu(3)$ via $\Ad$\\
\hline
$Q^7_G$ (all) & $\SU(3)\times S^1$ on $S^7\subset \C^3\times \C$ via $(A,z)\star(x,w)=(z^nAx,z^mw)$\\
\hline
$P^7_C$ (part) & $\SO(4)$ on $S^7$ via the isotropy representation of $\G_2/\SO(4)$\\
\hline
\end{tabular}
\end{center}
\caption{Nonreducible isometric cohomogeneity one actions on compact simply connected irreducible symmetric spaces in dimensions 5, 6 and 7 which are not sum actions or fixed-point actions. Also indicated is whether the family of actions listed in the right-hand column makes up all or part of the family listed in the left-hand column.}\label{nonred_irredsymspactions}
\end{table}}

For a complete list of cohomogeneity one actions on compact simply connected symmetric spaces we must only find such actions on product symmetric spaces. By looking at each such product individually, considering its full isometry group, then determining which subgroups of the isometry group can act by cohomogeneity one on the product, we can get a list of all possible actions. Many of these actions will be simple product actions and since we have already addressed those in general in Section \ref{special_types}, we will not consider them again here. The remaining actions which are also nonreducible are listed in Tables \ref{symsp_prod_actions1} and \ref{symsp_prod_actions2}.

{\setlength{\tabcolsep}{0.40cm}
\renewcommand{\arraystretch}{1.6}
\stepcounter{equation}
\begin{table}[!h]
\begin{center}
\begin{tabular}{|c|c|}
\hline
$Q^5_A$ & $S^3\times S^1$ on $S^2\times S^3\subset \Im(\H)\times \H$\\
 (part) &   via $(g,z)\star (p,q) = (gpg^{-1},gqz^{-1})$\\
\hline
$Q^5_A$ & $S^3\times S^1$ on $S^2\times S^3\subset \Im(\H)\times \H$\\
 (all)  &   via $(g,z)\star (p,q) = (z^np\bar z^n,gq\bar z^m)$\\
\hline
$N^6_A$ & $S^3\times T^2$ on $S^3\times S^3\subset \H\times \H$\\
 (part)  &   via $(g,z,w)\star (p,q) = (z^aw^bp\bar z^c\bar w^d, gq\bar z^n\bar w^m)$\\
\hline
$N^6_A$ & $S^3\times T^2$ on $S^3\times S^3\subset \H\times \H$\\
 (part) &   via $(g,z,w)\star (p,q) = (gp\bar z, gq\bar w)$\\
\hline
$Q^6_B$ & $S^3\times S^3$ on $S^3\times S^3\subset \H\times \H$\\
 (all)  &   via $(g_1,g_2)\star (p,q) = (g_1pg_1^{-1}, g_1qg_2^{-1})$\\
\hline
$Q^6_B$ & $\SO(4)$ on $S^3\times S^3$\\
 (all) &   via $A\star(x,y) = (Ax,Ay)$\\
\hline
$N^6_E$ & $S^3\times S^3$ on $S^2\times S^4\subset \Im(\H)\times (\H\times \R)$\\
 (part)  &   via $(g_1,g_2)\star (p,q,t) = (g_1pg_1^{-1}, g_1qg_2^{-1}, t)$\\
\hline
product & $S^3\times S^3$ on $S^3\times S^4\subset \H\times (\H\times \R)$\\
  &   via $(g_1,g_2)\star(p,q,t)=(g_1p,g_1qg_2^{-1},t)$\\
\hline
$N^7_A$ & $S^3\times S^3$ on $S^2\times S^2\times S^3\subset \Im(\H)\times \Im(\H)\times \H$\\
 (part) &  via $(g_1,g_2)\star(p_1,p_2,q)=(g_1p_1g_1^{-1},g_1p_2g_1^{-1},g_1qg_2^{-1})$\\
\hline
$N^7_A$ & $S^3\times S^3$ on $S^2\times S^2\times S^3\subset \Im(\H)\times \Im(\H)\times \H$\\
 (part) &  via $(g_1,g_2)\star(p_1,p_2,q)=(g_1p_1g_1^{-1},g_2p_2g_2^{-1},g_1qg_2^{-1})$\\
\hline
\end{tabular}
\end{center}
\caption{Nonreducible isometric cohomogeneity one actions on compact simply connected products of irreducible symmetric spaces in dimensions 5, 6 and 7 which are not product actions. Also indicated is whether the family of actions listed in the right-hand column makes up all or part of the family listed in the left-hand column. (Table 1 of 2)}\label{symsp_prod_actions1}
\end{table}}
{\setlength{\tabcolsep}{0.40cm}
\renewcommand{\arraystretch}{1.6}
\stepcounter{equation}
\begin{table}[!h]
\begin{center}
\begin{tabular}{|c|c|}
\hline
$P^7_A$ & $S^3\times S^3$ on $S^3\times \CP^2$ via $(g_1,g_2)\star(p,x)=(g_1pg_2^{-1},g_2\star_1x)$\\
 (part) &  where $\star_1$ is the action of $\SO(3)$ on $\CP^2$\\
\hline
$P^7_C$ & $S^3\times S^3$ on $S^3\times S^4\subset \H\times \R^5$ via $(g_1,g_2)\star(p,y)=(g_1pg_2^{-1},g_2\star_1y)$\\
 (part) &  where $\star_1$ is the action of $\SO(3)$ on $S^4$ from Section \ref{dim4class}\\
\hline
$P^7_D$ & $S^3\times S^3$ on $S^3\times \CP^2$ via $(g_1,g_2)\star(p,x)=(g_1pg_2^{-1},g_2\star_1x)$\\
 (part) &  where $\star_1$ is the action of $\SU(2)$ on $\CP^2$\\
\hline
$Q^7_B$ & $S^3\times S^3$ on $S^3\times S^4\subset \H\times (\H\times \R)$\\
   (all) &   via $(g_1,g_2)\star(p,q,t)=(g_1pg_2^{-1},g_2q,t)$\\
\hline
$Q^7_C$ & $S^3\times S^3\times S^1$ on $S^3\times S^4\subset \H\times (\Im(\H)\times \C)$\\
 (all)  &   via $(g_1,g_2,z)\star(p,q,w)=(g_1p\bar z^n,g_2qg_2^{-1},z^mw)$\\
\hline
$Q^7_D$ & $S^3\times S^3\times S^1$ on $S^3\times S^4\subset \H\times (\Im(\H)\times \C)$\\
 (part) &   via $(g_1,g_2,z)\star(p,q,w)=(g_1pg_2^{-1},g_2qg_2^{-1},zw)$\\
\hline
$N^7_F$ & $S^3\times S^3\times S^1$ on $S^2\times S^5\subset \Im(\H)\times (\H\times \C)$\\
 (part) &   via $(g_1,g_2,z)\star(p,q,w)=(g_1pg_1^{-1},g_1qg_2^{-1},zw)$\\
\hline
$Q^7_F$ & $\SU(3)\times S^1$ on $S^2\times S^5\subset \Im(\H)\times \C^3$\\
  (all) &   via $(A,z)\star(p,x)=(z^np\bar z^n,z^mAx)$\\
\hline
\end{tabular}
\end{center}
\caption{Nonreducible isometric cohomogeneity one actions on compact simply connected products of irreducible symmetric spaces in dimensions 5, 6 and 7 which are not product actions. Also indicated is whether the family of actions listed in the right-hand column makes up all or part of the family listed in the left-hand column. (Table 2 of 2)}\label{symsp_prod_actions2}
\end{table}}

In Tables \ref{nonred_irredsymspactions}, \ref{symsp_prod_actions1} and \ref{symsp_prod_actions2}, we also list the families to which the described actions belong. These families are easily determined by computing the group diagrams for each action. The group diagram for the last entry in Table \ref{nonred_irredsymspactions} is described in section 4 of \cite{GWZ}. In one case, the action is equivalent to a product action, even though the action itself is not a product action. In this case the word ``product'' is printed in the left-hand column. The table also indicates whether the given actions make up the entire family from the left-hand column, or just a part of it. The words ``all'' or ``part'' are written to distinguish these two cases.

\subsection{Brieskorn varieties}\label{brieskorn}

There is a well known cohomogeneity one action on the Brieskorn variety
 $$B^{2n-1}_d=\{z\in \C^{n+1} \, | \, z_0^d+z_1^2+z_2^2+\cdots+z_n^2=0, \sum_{i=0}^n  | z_i|^2=1\},$$
by the group $S^1\times\SO(n)$, given by
  $$(w,A)\star (z_0,z_1,z_2,\dots,z_n) = (w^2z_0,w^dA(z_1,z_2,\dots,z_n)^t).$$

These actions were extensively studied in \cite{GVWZ}. In particular they describe the group diagrams for the actions. In dimension 5 the group diagrams are
  $$S^3\times S^1 \,\, \supset \,\, \set{(e^{i\theta},1)}\cdot H, \,\, \set{(e^{jd\theta}, e^{2i\theta})}  \,\, \supset \,\,  \langle(j,-1)\rangle, \text{ for $d$ odd}$$
  $$S^3\times S^1 \,\, \supset \,\, \set{(e^{i\theta},1)}, \,\, \set{(e^{jd\theta}, e^{i\theta})}  \,\, \supset \,\,  1, \text{ for $d$ even}$$
where we have taken a more effective version of the diagram in the second case. The first diagram is precisely the diagram of $Q^5_B$ for $d=p$ and hence $Q^5_B\approx B^5_d$ for $d$ odd. Since $H$ is trivial in the second diagram, \ref{normalizer} says this diagram is equivalent to one of type $N^5$ for certain parameters.

In dimension 7, after lifting the action to $S^3\times S^3\times S^1$, the group diagrams are given by
  $$S^3\times S^3\times S^1 \,\, \supset \,\, \set{(e^{i\phi},e^{i\phi}e^{id\theta},e^{i\theta})}, \,\, \pm\Delta S^3\times \pm 1  \,\, \supset \,\, \pm\Delta S^1\times \pm1 $$
if $d$ is odd, where the $\pm\Delta S^3= \set{(g,\pm g)}$ and where the $\pm$ signs are correlated, and
  $$S^3\times S^3\times S^1 \,\, \supset \,\, \set{(e^{i\phi},e^{i\phi}e^{id\theta},e^{2i\theta})}, \,\, \Delta S^3\times 1  \,\, \supset \,\, \Delta S^1\times 1 $$
if $d$ is even.

This first diagram is exactly diagram $Q^7_D$ in the case that $n=2$ for $d=b$, since if $n=2$, $b$ must be odd for the diagram to be effective. The second diagram above, is exactly $Q^7_D$ in the case $n=1$ since if $d$ is even we can take $d=2b$ for $b$ in diagram $Q^7_D$. So the family $Q^7_D$ exactly corresponds to these actions on the Brieskorn varieties.

\subsection{Important actions in more detail}\label{others}

In this section we will look at the actions appearing in Tables \ref{class:prim} and \ref{class:nonprim} of the introduction one by one and summarize various facts that we have collected about these actions.

\subsubsection{Primitive actions of Table \ref{class:prim}}

\vspace{1ex}
\paragraph{$P^5$:}
One example of this family is the usual $\S(\U(2)\U(1))\subset\SU(3)$ action on $\SU(3)/\SO(3)$. This gives $P^5$ in the case $p=1$. In Section \ref{5dimtop} we will show that $P^5$ is in fact always diffeomorphic to $\SU(3)/\SO(3)$.

\vspace{1ex}
\paragraph{$P^7_A$:} This family is interesting because of its similarity to the families $P^7_B$ and $P^7_C$. One very special case of this family is the action of $S^3\times S^3$ on $S^3\times \CP^2$ given by  $(g_1,g_2)\star(p,x)=(g_1pg_2^{-1},g_2\star_1x)$, where $\star_1$ is the action of $\SO(3)$ on $\CP^2$, which corresponds to the case $p_-=q_-=p_+=q_+=1$.

\vspace{1ex}
\paragraph{$P^7_B$ and $P^7_C$:} These families were studied in \cite{GWZ}, where their homology groups were computed. These two classes contain all the new candidates for compact simply connected cohomogeneity one manifolds with invariant metrics of positive curvature, with one exception, as shown in \cite{GWZ}. The class $P^7_C$ was also shown to contain all $S^3$ principal bundles over $S^4$ in \cite{GZ1}. Two explicit actions of type $P^7_C$ are the isometric actions on $S^7$ and on $S^3\times S^4$ listed in Tables \ref{nonred_irredsymspactions} and \ref{symsp_prod_actions2}, respectively. One example of the family $P^7_B$ is the action of $\SO(3)\times \SO(3)$ on the Aloff-Wallach space $W^7=\SU(3)/\diag(z,z,\bar z^2)$, as described in section 4 of \cite{GWZ}.

\vspace{1ex}
\paragraph{$P^7_D$:} This family contains the cohomogeneity one Eschenburg spaces
  $$E^7_p = \diag(z,z,z^p)\backslash \SU(3)/\diag(1,1,\bar z^{p+2}),$$
where $\SU(2)\times \SU(2)$ acts on $E^7_p$ with the first factor acting on the left and the second on the right, both as the the upper $\SU(2)$ block in $\SU(3)$. These actions correspond to the case $n=2$ and $(p,q)=(p,p+1)$ in the family $P^7_D$. It should be noted that all of these Eschenburg spaces admit invariant metrics of positive sectional curvature, by \cite{Es}. For details, see \cite{GWZ}.

The action on $S^3\times \CP^2$ given in Table \ref{symsp_prod_actions2} is another example of type $P^7_D$, this time with $n=1$ and $p=q=1$. 

\subsubsection{Nonprimitive actions of Table \ref{class:nonprim}}

Recall from Proposition \ref{prop:prim} that for a nonprimitive action of $G$ on $M_G$, with $G\supset L\supset K^-, K^+\supset H$, we have the fiber bundle $M_L\to M_G\to G/L$, where $M_L$ is the cohomogeneity one manifolds given by the diagram $L\supset K^-, K^+\supset H$. Having such a fiber bundle says something about the topology of the manifold $M_G$. Because of this, we will list these bundles below. More details about how we get the specific fiber bundles below can be found in Section \ref{curvsec}.

\vspace{1ex}
\paragraph{$N^5$:}
We saw in Section \ref{brieskorn} that the Brieskorn varieties for $d$ even are all examples of this family. There is one more explicit action that is also of this type. Let $S^3\times S^1$ act on $S^3\times S^2 \subset \H\times \Im \H$ via $(g,\theta)\star (q,v)=(gq\bar ge^{gv\bar g\theta},gv\bar g)$. This gives the diagram
$S^3\times S^1\,\, \supset \,\, \set{(e^{i\theta},1)}, \,\, \set{(e^{i\theta},e^{2i\theta})} \,\, \supset \,\, \set{(\pm1,1)}$
which is a special case of $N^5$. In Section \ref{5dimtop}, we will show every manifold of this type is either $S^3\times S^2$ or the nontrivial $S^3$ bundle over $S^2$.

\vspace{1ex}
\paragraph{$N^6_A$:}
A large subclass of this family of actions is given by the $S^3\times S^1\times S^1$ action on $S^3\times S^3$ given by $(g, z,w)\star (x,y) = (gx\bar z^r\bar w^s, (z^{c_-}\bar w^{b_-},z^{c_+}\bar w^{b_+})\star_1 y)$ where $\star_1$ is the usual torus action on $S^3\subset \C^2$. This action gives the diagram
$$S^3\times T^2  \supset  \set{(z^rw^s,z,w)|z^{c_-}\bar w^{b_-}=1},  \set{(z^rw^s,z,w)|z^{c_+}\bar w^{b_+}=1} \supset$$ $$\nonumber  \set{(z^rw^s,z,w)|z^{c_-}\bar w^{b_-}\!=\!1\!=\!z^{c_+}\bar w^{b_+}}$$
however this family does not constitute every action of type $N^6_A$ since one can always extend $H$ to get a new diagram, not of the form above.

For a general manifold $M$ in the family $N^6_A$, if we take $L=T^3\subset S^3\times S^1\times S^1$ we get the nonprimitivity fiber bundle $S^3\times S^1\to M\to S^2$ or if we take $L=T^2=K^-_0\cdot K^+_0=K^-\cdot K^+$ we get the different bundle $S^3\to M\to S^3$.

%

\vspace{1ex}
\paragraph{$N^6_B$:}Here, if we take $L=T^2$ we get the fiber bundle $S^2\to M \to S^2\times S^2$, for any $M$ in the family $N^6_B$. If $p=0$ in this family, we get the product action on $S^2\times M^4$ where $M^4$ is either $S^2\times S^2$ or $\CP^2\#-\CP^2$, depending on whether $n$ is even or odd, as described in section \ref{dim4class}.

\vspace{1ex}
\paragraph{$N^6_C$:}For each $M$ of this type, we can take $L=S^3\times S^1$ to get the fiber bundle $S^4\to M \to S^2$.

\vspace{1ex}
\paragraph{$N^6_D$:}If we let $L=S^3\times S^1$ in this case as well, we get the fiber bundle $\CP^2\to M \to S^2$, for manifolds $M$ of this type. In the case $p=0$, $N^6_D$ is the product action on $\CP^2\times S^2$.

\vspace{1ex}
\paragraph{$N^6_E$:}
In the case that $p=1$ we get the following $S^3\times S^3$ action on $S^2\times S^4\subset \Im \H\times (\H\times \R)$: $(g_1, g_2)\star (x,(y,t)) = (g_1xg_1^{-1}, (g_1yg_2^{-1},t))$. In this case we get the diagram
$$S^3\times S^3 \,\, \supset \,\, S^3\times S^1, \,\, S^3\times S^1 \,\, \supset \,\, \Delta S^1.$$
Similarly, when $p=0$ this is a product action on $S^2\times S^4$.

For a general $M$ in this family, if we take $L=S^3\times S^1$ we get the nonprimitivity fiber bundle $S^4\to M\to S^2$.

\vspace{1ex}
\paragraph{$N^6_F$:} One special case of this class of actions, is the $\SU(3)$ action on $\CP^3\#-\CP^3$ obtained by gluing two copies of the $\SU(3)$ action on $\CP^3$ along the fixed point. We get
$$\SU(3) \,\, \supset \,\, \S(\U(2)\U(1)), \,\, \S(\U(2)\U(1)) \,\, \supset \,\, \SU(2)\SU(1).$$

In general, for any $M$ in this family, we can take $L=\S(\U(2)\times \U(1))$ and get the fiber bundle $S^2\to M\to \CP^2$.

\vspace{1ex}
\paragraph{$N^7_A$:} One special case of this family is the $S^3\times S^3$ action on $S^2\times S^2\times S^3\subset \Im(\H)\times \Im(\H)\times \H$ via $(g_1,g_2)\star(p_1,p_2,q) = (g_1p_1g_1^{-1},g_1p_2g_1^{-1},g_1qg_2^{-1})$, or the equivalent action of $(g_1,g_2)\star(p_1,p_2,q) = (g_1p_1g_1^{-1},g_2p_2g_2^{-1},g_1qg_2^{-1})$. Also, the case $(p_\pm,q_\pm)=(0,1)$ and $H=\Z_n$ gives the product action on $S^3\times M^4$ where $M^4$ is either $\CP^2\#-\CP^2$ or $S^2\times S^2$ depending on whether $n$ is odd or even.

For this family, the nonprimitivity fiber bundle depends heavily on the parameters of the action. For any $M$ in the family $N^7_A$ corresponding to a diagram with $K^-\ne K^+$, we can take $L=T^2$ to get the fiber bundle $L_m(n)\to M\to S^2\times S^2$, where $L_m(n)$ is some lens space which depends on the parameters of $M$ in $N^7_A$. If $M$ is a member of this family with $K^-=K^+$ then taking $L=T^2$ gives the fiber bundle $S^2\times S^1\to M\to S^2\times S^2$ and taking $L=K^-=K^+$ in this case gives the bundle $S^2\to M \to S^3\times S^2$.

\vspace{1ex}
\paragraph{$N^7_B$:}For each $M$ in this family, we can take $L=S^3\times S^1$ to get the fiber bundle $M_L\to M\to S^2$, where $M_L$ is the cohomogeneity one manifold given by the diagram
  $$S^3\times S^1\supset \set{(e^{ip\theta}, e^{iq\theta})}\cdot H, \set{(e^{j\theta},1)}\cdot H\supset H$$
with the same restrictions on $H$ as in $N^7_B$. The manifolds $M_L$ will depend greatly on the parameters in the diagram. For example the actions $Q^5_B$ and $P^5$, the Brieskorn actions and the family of actions on $\SU(3)/\SO(3)$ respectively, as described in Sections \ref{brieskorn} and \ref{5dimtop-N^5} respectively, are both actions of this type. In fact these are the only cases when $M_L$ will be simply connected, assuming $p\ne0$. Of course if $p=0$ then the original action would be of type $N^7_A$.

In the case $q=0$, action $N^7_B$ becomes the product action on $S^4\times S^3$ or on $\CP^2\times S^3$, depending on $H$.

\vspace{1ex}
\paragraph{$N^7_C$:} For a manifold $M$ of this type, we can take $L=S^3\times S^1$ to get the fiber bundle $S^5/\Z_q\to M\to S^2$, if $q\ne0$, where $S^5/\Z_q$ is the lens space, as described in Section \ref{curvsec}. If $q=0$ then the original action is just a product action on $\CP^2\times S^3$.

\vspace{1ex}
\paragraph{$N^7_D$:} In this case we can take $L=T^2=K^+=K^-$ to get the fiber bundle $S^2\to M\to S^3\times S^2$, for any $M$ in this family. We can also take $L=T^3\subset S^3\times S^3\times S^1$ to get the bundle $S^2\times S^1\to M\to S^2\times S^2$, for any such $M$.

If $q=\lambda=0$ then this action is the product action on $S^2\times Q^5_A$, where $Q^5_A$ is the family of actions on $S^2\times S^3$ described above. Also if $m=0$ this is another product action on $S^3\times S^2\times S^2$, since it is known that $S^3\times S^3/\set{(z^p,z^q)}=S^3\times S^2$ (see \cite{WZ}).

\vspace{1ex}
\paragraph{$N^7_E$:} For each $M$ in this family, let $L=T^3\subset S^3\times S^3\times S^1$. This gives the fiber bundle $L_a(b)\to M\to S^2\times S^2$, where $L_a(b)$ is some lens space which depends on the parameters of $M$ in $N^7_E$ in a complicated way.

As in the previous case, if we take $q=\lambda=0$ we get a product action on $S^2\times N^5$. Later, we will show that the actions of type $N^5$ are always on $S^3\times S^2$ or the nontrivial $S^3$ bundle over $S^2$.

\vspace{1ex}
\paragraph{$N^7_F$:} The action of $S^3\times S^3\times S^1$ on $S^2\times S^5\subset \Im(\H)\times (\H\times \C)$ given by $(g_1,g_2,z)\star(p,q,w)=(g_1pg_1^{-1},g_1qg_2^{-1},zw)$ is one example of this family. Also, in the case $p=0$, we get the product action on $S^2\times Q^5_C$, where $Q^5_C$ is an action on $S^5$.

For a general $M$ in this family, taking $L=S^3\times S^1\times S^1$ gives the fiber bundle $S^5\to M\to S^2$.

\vspace{1ex}
\paragraph{$N^7_G$:}For this manifold, we have the bundle $S^3\to N^7_G\to \CP^2$, after taking $L=U(2)=K^\pm$.

\vspace{1ex}
\paragraph{$N^7_H$:} In this case, if we take $L=\U(2)\times S^1\subset \SU(3)\times S^1$ we get the fiber bundle $L_a(b)\to M\to \CP^2$, for each $M$ in this family, where $L_a(b)$ is some lens space which depends on the parameters of $M$ in the class $N^7_H$.

\vspace{1ex}
\paragraph{$N^7_I$:}Finally, we have the fiber bundle $S^3\to N^7_I\to S^4$, for this manifold after taking $L=\Sp(1)\Sp(1)=K^\pm$.

\section{Curvature properties}\label{curvsec}

Here we will prove Theorem \ref{curvature}. First, we have shown, through the classification, that every nonreducible cohomogeneity one action on a simply connected manifold in dimension 5, 6 or 7 must be a product action, a sum action, a fixed-point action or one of the actions in Tables \ref{5dimlist} through \ref{7dimlist2}. We know from Section \ref{special_types} that sum actions and fixed-point actions are isometric actions on symmetric spaces, and hence they admit invariant metrics of nonnegative curvature. For product actions, suppose $M=N\times L/J$ where $N$ is a lower dimensional cohomogeneity one manifold, $L/J$ is is a homogeneous space, and $G=G_1\times L$ acts as a product. It is clear that the action of $G_1$ on $N$ is nonreducible if and only if the original action of $G$ on $M$ is. From the classification of lower dimensional cohomogeneity one manifolds, we see that all the possibilities for $N$ admit $G_1$-invariant metrics of nonnegative sectional curvature. Obviously $L/J$ admits an $L$ invariant metric of nonnegative sectional curvature. Therefore, all nonreducible product actions admit $G$ invariant metrics of nonnegative sectional curvature as well.

Therefore we must only check that the actions listed in Tables \ref{5dimlist} through \ref{7dimlist2} admit invariant metrics of nonnegative sectional curvature, except for the two families $P^7_D$ and $Q^7_D$ listed in Theorem \ref{curvature} as exceptions. Section \ref{symspactions} shows that many of these actions are isometric actions on symmetric spaces, and hence admit invariant metrics of nonnegative curvature. Many of the actions also have codimension two singular orbits. Therefore, by the main result in \cite{GZ1}, these also admit invariant nonnegative curvature. After these two considerations, we are only left with a handful of actions that we still need to check: $N^6_C$, $N^6_D$, $N^6_E$, $N^7_C$, $N^7_F$, $N^7_G$ and $N^7_I$.

Notice that these seven actions are all nonprimitive. We will use Proposition \ref{prop:prim} to write each manifold $M$ of these types in the form $G\times_L M_L$ and in each case we will see that the $L$ action on $M_L$ admits an invariant metric of nonnegative sectional curvature. This will show that $M\approx G\times_L M_L$ admits a $G$ invariant metric of nonnegative sectional curvature, since we can take the metric mentioned above on $M_L$ and a biinvariant metric on $G$ to induce a submersed metric on $M$. This metric will still be nonnegatively curved by O'Niell's formula (see \cite{Petersen}). We will proceed to do this case by case.

In the case of $N^6_C$, the subdiagram corresponding to $M_L$ is given by
  $$S^3\times S^1 \,\, \supset \,\, T^2, \,\, S^3\times \Z_n \,\, \supset \,\, S^1\times \Z_n.$$
This action is ineffective, with the effective version given by taking $n=1$. We then recognize this effective version as an isometric sum action on $S^4$. Therefore $M_L$ admits an $L$ invariant metric of nonnegative sectional curvature.

We can do a similar thing for actions $N^6_D$ and $N^6_E$. The primitive subdiagram for $N^6_D$ is
  $$S^3\times S^1 \,\, \supset \,\, T^2, \,\, S^3\times S^1 \,\, \supset \,\, \set{(e^{ip\theta}, e^{i\theta})}.$$
This is the action of $\SU(2)\times S^1$ on $\CP^2$ given by $(A,w)\star [z_1,z_2,z_3] = [w^pz_1,A(z_2,z_3)^t]$. Since this is an isometric action for the usual metric on $M_L=\CP^2$, this gives an invariant metric on $M\approx G\times_L M_L$. Similarly the primitive subdiagram for $N^6_E$ is
  $$S^3\times S^1 \,\, \supset \,\, S^3\times S^1, \,\, S^3\times S^1 \,\, \supset \,\, \set{(e^{ip\theta}, e^{i\theta})}.$$
This is the action of $S^3\times S^1$ on $S^4\subset \H \times \R$ by $(g,z)\star (p,t)= (gp\bar z^p, t)$. As above this also gives $N^6_E$ an invariant metric of nonnegative sectional curvature.

For $N^7_C$, the primitive subdiagram is given by
  $$S^3\times S^1 \,\, \supset \,\, \set{(e^{ip\theta},e^{iq\theta})}, \,\, S^3\times \Z_n  \,\, \supset \,\, \Z_n$$
where $\gcd(p,q)=1$ and $\gcd(q,n)=1$. We claim this is an isometric action on the lens space $S^5/\Z_q$. It is easy to check that the special case of this action, when $q=1$, is the modified sum action of $\SU(2)\times S^1$ on $S^5\subset \C^2\times \C$ by $(A,w)\star (x,z)= (w^pAx,w^nz)$. Then consider $S^5/\Z_q$ where $\Z_q$ acts as $\Z_q\subset 1\times \Z_q$ with this same action, $\star$. Then $\SU(2)\times S^1$ still acts on the quotient $S^5/\Z_q$ and does so isometrically in the induced metric. If the original group diagram is taken along the geodesic $c$ in $S^5$ then we can take the group diagram of the induced action on $S^5/\Z_q$ along the image of $c$. When we do this we see we get exactly the diagram shown above. Hence this is an isometric action on $S^5/\Z_q$ in the usual, positively curved, metric. As above, this induces an invariant metric on nonnegative sectional curvature on $N^7_C$.

The last three cases are slightly easier to handle. For $N^7_F$ the primitive subdiagram is given by
  $$S^3\times S^1\times S^1 \,\, \supset \,\, \set{(e^{ip\phi}e^{ia\theta},e^{i\phi},e^{i\theta})}, \,\, S^3\times S^1\times \Z_n  \,\, \supset \,\, \set{(e^{ip\phi},e^{i\phi},1)}\cdot \Z_n.$$
We see this is the modified sum action of $S^3\times S^1\times S^1$ on $S^5\subset \H\times \C$ given by $(g,z_1,z_2)\star (y,w)= (gy\bar z_1^p\bar z_2^a, z_2^nw)$. Similarly the primitive subdiagram for $N^7_G$ is given by the action of $\U(2)$ on $\SU(2)\approx S^3$ by conjugation and the subdiagram for $N^7_I$ corresponds to the action of $\Sp(1)$ on $\Sp(1)$ by conjugation. Both of these are isometries in the positively curved biinvariant metric on $\Sp(1)\approx\SU(2)\approx S^3$. Therefore these seven remaining cases do admit invariant metrics of nonnegative sectional curvature. This proves Theorem \ref{curvature}.

\section{Topology of the 5-dimensional manifolds}\label{5dimtop}

In this section we will determine the diffeomorphism type of the five dimensional manifolds appearing in the classification, and prove \tref{5class}. By the results of Smale and Barden (\cite{Barden}) the diffeomorphism type of a closed simply connected 5-manifold is determined by the second homology group and the second Stiefel-Whitney class. As we will see, we can compute the homology of our manifolds relatively easily. To compute the second Stiefel-Whitney class, we will use the topology of the frame bundle. Recall that the second Stiefel-Whitney class of a simply connected manifold is zero if and only if the manifold is a $\Spin$-manifold, i.e. the orthonormal frame bundle lifts to a $\Spin$-bundle (see \cite{Petersen}). With this motivation, we will now look at the frame bundle in more detail.

Suppose $M^n$ is an oriented cohomogeneity one manifold with the group diagram $G\supset K^-,K^+\supset H$ as usual. Assume further that $G$ is connected so that the $G$ action preserves the orientation of $M$. Then let
  $$FM=\set{f=(f_1,\dots,f_n) | f_1,\dots,f_n \text{ is an oriented o.n. frame at $p\in M$} }$$
denote the orthonormal oriented-frame bundle of $M$. Recall that $\SO(n)$ acts on $FM$ from the left as
\begin{equation}\label{SO(n)actn}
(a_{ij})_{ij} \star (f_1,\dots,f_n) = (\sum_j a_{1j}f_j, \cdots, \sum_j a_{nj} f_j).
\end{equation}
This action makes $FM$ into an $\SO(n)$-principal bundle over $M$. We can put a metric on $FM$ by choosing a biinvariant metric on $\SO(n)$, keeping the original metric on $M$ and specifying a horizontal distribution. To describe this distribution, fix a point $p_0\in M$ and a frame $f_{p_0}$ at $p_0$. For each $p$ in a normal neighborhood of $p_0$, let $f_p$ be the frame gotten by parallel translating $f_{p_0}$ to $p$ along the minimal geodesic from $p_0$ to $p$. This gives a local orthonormal frame field, and we define the horizontal space at $f_{p_0}\in FM$ to be the tangent space of this frame field. Since parallel transport is an isometry, the action of $\SO(n)$ preserves this horizontal distribution.

Recall that $G$ acts on $M$ by isometry and hence takes orthonormal frames to orthonormal frames, while preserving orientation. Therefore we have an induced action of $G$ on $FM$ given by $g\star(p,f)= (gp,dgf) = (gp,(dgf_1,\dots,dgf_n))$. This action is isometric since it takes the horizontal space to the horizontal space and acts by isometry on both the vertical and horizontal spaces. This $G$-action also commutes with the action by $\SO(n)$ and so we have an action by $G\times \SO(n)$ on $FM$. Furthermore, this $G\times \SO(n)$ action on $FM$ is clearly cohomogeneity one since $\SO(n)$ acts transitively on the fibers of $FM$. If $c$ denotes a minimal geodesic in $M$ between nonprincipal orbits, then choose $f(t)$ to be a parallel orthonormal frame along $c$. Then $f$ is a horizontal curve in $FM$ and therefore a geodesic. $f(t)$ is clearly perpendicular to the $\SO(n)$ orbits and it is perpendicular to the $G$ orbits since $c(t)$ is in $M$. Therefore $f(t)$ is a minimal geodesic in $FM$ between nonprincipal orbits.

Our next goal is to determine the isotropy groups of $G\times \SO(n)$ along $f(t)$. We see $(g,A)\star (p,f)= (p,f)$ if and only if $g\in G_p$ and $A\star dg f= f$, where $\star$ is from \ref{SO(n)actn}. To understand this second equality we rewrite it as
  $$(dg^{-1}f_1,\dots,dg^{-1}f_n) = dg^{-1} f= A\star f = (\sum_j a_{1j}f_j, \cdots, \sum_j a_{nj} f_j).$$
This precisely means $A^t=[dg^{-1}]_f$ or $A= [dg]_f$ where $[dg]_f$ denotes the representation of the linear operator $dg: T_p M\to T_p M$ as a matrix in the basis $f_1,\dots,f_n$. In conclusion, the isotropy group of $G\times \SO(n)$ at $(p,f)$ is $\set{(g,[dg]_f)|g\in G_p}$. We have proven the following proposition.

\begin{prop}\label{frame_bundle}
Let $M^n$ be an oriented cohomogeneity one manifold with group diagram $G\supset K^-,K^+\supset H$ for the normal geodesic $c$, and assume $G$ is connected. The orthonormal oriented frame bundle $FM$ of $M$ admits a natural cohomogeneity one action with group diagram
\begin{equation}\label{frame_diagram}
\begin{split}
\xymatrix{
&   G\times \SO(n)   \ar@{-}[dr] \ar@{-}[dl] &\\
\set{(k,[dk]_{f(-1)})|k\in K^-}  \ar@{-}[dr] & &   \set{(k,[dk]_{f(1)})|k\in K^+}   \ar@{-}[dl] \\
 & \set{(h,[dh]_{f(0)})|h\in H} &}
\end{split}
\end{equation}
where $f(t)$ is a parallel frame along $c(t)$.
\end{prop}

This proposition, together with \ref{Hgens}, gives the following corollary.

\begin{cor}\label{frame_bundle_cor}
Let $M$ be a cohomogeneity one manifold as in \ref{frame_bundle} and assume that $H$ is discrete. Let $\alpha_\pm:[0,1]\to K^\pm$ be paths, based at the identity, which generate $\pi_1(K^\pm/H)$. If $M$ is simply connected then $FM$ is simply connected if and only if there is some curve $\gamma=\alpha_-^n\cdot\alpha_+^m$ which gives a contractible loop in $G$ and where $[d\alpha_-^n]_{f(-1)}\cdot[d\alpha_+^m]_{f(1)}$ generates $\pi_1(\SO(n))$.
\end{cor}
\begin{proof}
Notice that the maps $k\mapsto (k,[dk]_{f(\pm1)})$ give isomorphisms of $K^\pm$ with $\widehat K^\pm:=\set{(k,[dk]_{f(\pm1)})|k\in K^\pm}$, the nonprincipal isotropy subgroups of the $\widehat G:=G\times \SO(n)$ action on $FM$. Furthermore, this map takes $H$ to $\widehat H:=\set{(h,[dh]_{f(0)})|h\in H}$, the principal isotropy group of this action. Therefore we see that $\widehat H$ is generated by $\widehat H\cap \widehat K^-_0$ and $\widehat H\cap \widehat K^+_0$ by \ref{Hgens}, since $M$ is already assumed to be simply connected. Then by \ref{Hgens}, $FM$ is simply connected if and only if the curves $\hat \alpha_\pm(t)=(\alpha_\pm(t),[d\alpha_\pm(t)]_{f(\pm1)})$ generate $\pi_1(\widehat G/\widehat H_0)$. Further, since $H$ is discrete, $\pi_1(\widehat G/\widehat H_0)=\pi_1(\widehat G)\approx \pi_1(G)\times \pi_1(\SO(n))$.

Therefore, if $\pi_1(FM)=0$ then the curve $\gamma$ from the corollary must exist. Conversely, suppose such a curve $\gamma$ exists. We already know from \ref{Hgens}, that $\alpha_-$ and $\alpha_+$ generate $\pi_1(G)$, since $M$ is simply connected. Then it is clear that $\hat \alpha_-$, $\hat \alpha_+$ and $\gamma$ would generate all of $\pi_1(\widehat G)$, proving $\pi_1(FM)=0$.
\end{proof}

\subsection{The family $P^5$}\label{5dimtop-P^5}
We will now compute the homology of the manifolds $P^5$ using the Mayer-Vietoris sequence. Because of the commutative diagram \ref{orbitcomdiagram}, we have the Mayer-Vietoris sequence for the spaces $M$, $G/K^-$, $G/K^+$ and $G/H$:
\begin{equation}\label{M-V_seq}
\begin{CD}
\cdots\to H_n(G/H)@>(\rho_*^-,\rho_*^+)>> H_n(G/K^-)\oplus H_n(G/K^+) @>i_*^--i_*^+>>H_n(M)\to\cdots\\
\end{CD}
\end{equation}
To compute $H_n(P^5)$ in our case, first notice that $G/K^+= S^3\times S^1/\set{(e^{jp\theta},e^{i\theta})}\approx S^3$, since $S^3$ acts transitively on this space with trivial isotropy group.

Next we claim that $G/H=S^3\times S^1/ \langle(j,i)\rangle$ is $S^3\times S^1$. For this, denote $\alpha: S^3\times S^1\to S^3\times S^1: (g,z)\mapsto (gj,zi)$, so that $G/H=G/\langle \alpha \rangle$. Then define the map $\phi: S^3\times S^1\to S^3\times S^1: (g,e^{i\theta})\mapsto (ge^{-j\theta},e^{i\theta})$, a diffeomorphism of manifolds. We notice that $\beta:=\phi\alpha\phi^{-1}:(g,z)\mapsto (g,zi)$. Therefore $G/\langle \alpha \rangle$ is diffeomorphic to $G/\langle \beta \rangle\approx S^3\times S^1$.

Finally we study $G/K^-= S^3\times S^1/\set{(e^{i\theta},1)}\cdot \langle(j,i)\rangle$. Since $K^-_0$ is normal in $K^-$ we have $G/K^-\approx (G/K^-_0)/(K^-/K^-_0)$. We see $G/K^-_0 = S^3\times S^1/\set{(e^{i\theta},1)}\approx S^2\times S^1 = \Im(S^3) \times S^1$ via $(gS^1,z)\mapsto (gig^{-1},z)$. Further, it is easy to see that $(j,i)$ acts on $S^2\times S^1$ as $(-\Id, i)$, via this correspondence. Hence $G/K^-\approx S^2\times S^1/\langle (-\Id, i)\rangle$. We can identify this space with $S^2\times [0,1] / \sim$ where $(x,0)\sim (-x,1)$. Using Mayer-Vietoris for this space we can easily compute that
  $$H_i(G/K^-)= \left\{ \begin{array}{cc}
     \Z & \text{ if } i=1\\
     \Z_2 & \text{ if } i=2\\
     0  & \text{ otherwise}
  \end{array} \right.$$

We are now ready to use the Mayer-Vietoris sequence for $P^5$. \ref{M-V_seq} becomes:
  $$\cdots\to 0 \to  \Z_2 \oplus 0 \to H_2(P^5) \to \Z \to \Z \oplus 0 \to 0 \to \cdots$$
since we know $H_1(P^5)=0$. Since the map $\Z \to \Z \oplus 0$ is onto, it must have trivial kernel and hence the map from $H_2(P^5)$ must be trivial. Therefore $\Z_2 \oplus 0 \to H_2(P^5)$ must be an isomorphism. That is
  $$H_2(P^5) = \Z_2.$$
Poincare duality then determines the rest of the homology groups.

To determine the second Steifel-Whitney class, we look at the fundamental group of the frame bundle $F(P^5)$. In the notation of \ref{frame_bundle_cor}, we can take $\alpha_-(\theta)=(e^{i\theta},1)$, in this case, since this is a curve in $K^-$ which generates $\pi_1(K^-/H)$. We need to determine how $d(\alpha_-(\theta))$ acts on $T_{c(-1)}M\approx T_{K^-}(G/K^-)\oplus D_-$. On $T_{K^-}(G/K^-)$, $d(\alpha_-(\theta))$ has the form $\diag(R(2\theta),1)$ in the basis $\set{(j,0),(k,0),(0,1)}$ and since $d(\alpha_-(\theta))$ is an isometry of $T_{K^-}(G/K^-)$ there must be an orthonormal basis of $T_{K^-}(G/K^-)$ in which $d(\alpha_-(\theta))$ still has this form. On $D_-$, $d(\alpha_-(\theta))$ acts isometrically as $R(\theta)$. Therefore there is an oriented orthonormal basis $f_-$ of $T_{c(-1)}M$ for which $[d(\alpha_-(\theta))]_{f_-}=\diag(R(2\theta),1,R(\theta))$. Since this generates $\pi_1(\SO(5))$ and since $\alpha_-$ is contractible in $G$, it follows from \ref{frame_bundle_cor} that $F(P^5)$ is simply connected, independent of $p$.

Therefore $P^5$ is not $\Spin$, and hence has nontrivial second Steifel-Whitney class for each $p$. By the results of Smale and Barden mentioned above this proves that the diffeomorphism type of $P^5$ is independent of $p$. In Section \ref{others}, we showed that $\SU(3)/\SO(3)$ is one example in this family. Hence $P^5$ is diffeomorphic to $\SU(3)/\SO(3)$ for all $p$.

\subsection{The family $N^5$}\label{5dimtop-N^5}
We will now compute the topology of the manifolds $N^5$, this time, using the nonprimitive fiber bundle. Notice first that these manifolds are not primitive. In fact if we take $L=T^2$ then $K^-,K^+,H\subset L$. Therefore, by \ref{prop:prim}, $N^5$ is fibered over $G/L\approx S^2$ with fiber $M_L$ where $M_L$ is the cohomogeneity one manifold given by
  $$S^1\times S^1 \,\, \supset \,\, \set{(e^{ip_-\theta},e^{iq_-\theta})}\cdot H, \,\, \set{(e^{ip_+\theta},e^{iq_+\theta})}\cdot H  \,\, \supset \,\,  H_-\cdot H_+.$$
Since $H$ is normal in $T^2$ it follows that the effective version of the $L$ action on $M_L$ is given by
\begin{equation}\label{N^5_fiber}
S^1\times S^1 \,\, \supset \,\, \set{(e^{i\theta},1)}, \,\, \set{(e^{ip\theta},e^{iq\theta})}  \,\, \supset \,\,  1
\end{equation}
after taking an automorphism of $T^2$, where $q\ne0$ since $K^+\ne K^-$ in the original diagram. To identify this action first recall that $T^2$ acts by cohomogeneity one on $S^3\subset \C^2$, by multiplication on each factor. If we take $S^3/\langle(\xi_q,\xi_q^p)\rangle$, where $\xi_q$ is a $q$th root of unity, this gives the lens space $L_q(p)$. Since the $T^2$ action on $S^3$ commutes with this subaction by $\Z_q$, we get an induced action on $L_q(p)$. It is not difficult to see that the effective version of this action is precisely the action given by \ref{N^5_fiber}. Therefore $M_L$ is $L_q(p)$ and we have the fibration
  $$L_q(p)\to N^5 \to S^2.$$

Given that $N^5$ is simply connected, the long exact sequence of homotopy groups induced from this fibration contains the following short exact sequence
  $$0\to \pi_2(N^5)\to\pi_2(S^2)\to \pi_1(L_q(p))\to 0.$$
Since the middle group is $\Z$ and the last group is $\Z_q$, for $q\ne0$, it follows that $\pi_2(N^5)\approx \Z$ and hence $H_2(N^5)\approx \Z$, by the Hurewicz theorem.

We claim here that the frame bundle $F(N^5)$ can either be simply connected or not, depending on the parameters of the diagram. In Section \ref{others}, we saw one example of an action in this family on $S^3\times S^2$. This shows that some of these actions are on spin manifolds. To see that some of these manifolds are not spin we take a simple example. The manifold $M_1$ with group diagram
\begin{equation}\label{nonspineg}
S^3\times S^1 \,\, \supset \,\, S^1\times 1, \,\, 1\times S^1  \,\, \supset \,\,  1
\end{equation}
is an example of type $N^5$. If we let $\alpha_-(\theta)=(e^{i\theta},1)$, then $\alpha_-$ generates $\pi_1(K^-/H)$. By precisely the same argument as in the case of $P^5$, we see that $F(M_1)$ is simply connected. Therefore the family $N^5$ contains both spin and nonspin manifolds, but always with the homology of $S^3\times S^2$. Using \cite{Barden} again, this proves \tref{5class}.

%
%
%
%
%
%

\section{Appendix}

In this appendix we give several tables for the convenience of the reader. The first four tables list all the nonreducible cohomogeneity one actions from our classification which are not sum, product or fixed point actions. The last table summarizes the relation between sum, product and fixed point actions and their group diagrams. These special types of actions are described in more detail in Section \ref{special_types}.

{\setlength{\tabcolsep}{0.40cm}
\renewcommand{\arraystretch}{1.6}
\stepcounter{equation}
\begin{table}[!h]
\begin{center}
\begin{tabular}{|c|c|}
\hline
$Q^5_A$ & $S^3\times S^1\,\, \supset \,\, \set{(e^{ip\theta},e^{i\theta})}, \,\, \set{(e^{ip\theta},e^{i\theta})} \,\, \supset \,\, \Z_n$\\
\hline
$N^5$ & $S^3\times S^1 \,\, \supset \,\, \set{(e^{ip_-\theta},e^{iq_-\theta})}\cdot H, \,\, \set{(e^{ip_+\theta},e^{iq_+\theta})}\cdot H  \,\, \supset \,\,  H_-\cdot H_+$\\
 &  $K^-\ne K^+$, $(q_-,q_+)\ne \mathbf{0}$,  $\gcd(q_-,q_+,d)=1$\\
 &  where $d= \#(K^-_0\cap K^+_0)/\#(H\cap K^-_0\cap K^+_0)$\\
\hline
$Q^5_B$ & $S^3\times S^1 \,\, \supset \,\, \set{(e^{i\theta},1)}\cdot H, \,\, \set{(e^{jp\theta}, e^{2i\theta})}  \,\, \supset \,\,  \langle(j,-1)\rangle$\\
 &  where $p>0$ is odd.\\
\hline
$P^5$ & $S^3\times S^1 \,\, \supset \,\, \set{(e^{i\theta},1)}\cdot H, \,\, \set{(e^{jp\theta}, e^{i\theta})}  \,\, \supset \,\,  \langle(j,i)\rangle$\\
 &  where $p\equiv 1 \mod 4$\\
\hline
$Q^5_C$ & $S^3\times S^1 \,\, \supset \,\, \set{(e^{ip\theta},e^{i\theta})}, \,\, S^3\times\Z_n  \,\, \supset \,\,  \Z_n$\\
\hline
\end{tabular}
\end{center}
\vspace{0.1cm}
\caption{Nonreducible 5-dimensional diagrams which are not products, sums or fixed point actions.}\label{5dimlist}
\end{table}}

{\setlength{\tabcolsep}{0.40cm}
\renewcommand{\arraystretch}{1.6}
\stepcounter{equation}
\begin{table}[!h]
\begin{center}
\begin{tabular}{|c|c|}
\hline
$N^6_A$ & $S^3\times T^2 \,\, \supset \,\, \set{(e^{i a_-\theta},e^{ib_-\theta},e^{ic_-\theta})}\cdot H, \,\, \set{(e^{i a_+\theta},e^{ib_+\theta},e^{ic_+\theta})}\cdot H \,\, \supset \,\, H$\\
 & where $K^-\ne K^+$, $H=H_-\cdot H_+$, $\gcd(b_\pm,c_\pm)=1$,\\
 & $a_\pm=rb_\pm+sc_\pm$, and $K^-_0\cap K^+_0\subset H$\\
\hline
$N^6_B$ & $S^3\times S^3 \,\, \supset \,\, \set{(e^{i\theta},e^{i\phi})}, \,\, \set{(e^{i\theta},e^{i\phi})} \,\, \supset \,\, \set{(e^{ip\theta},e^{iq\theta})} \cdot \Z_n$\\
\hline
$N^6_C$ & $S^3\times S^3 \,\, \supset \,\, T^2, \,\, S^3\times \Z_n \,\, \supset \,\, S^1 \times \Z_n$\\
\hline
$Q^6_A$ & $S^3\times S^3 \,\, \supset \,\, T^2, \,\, \Delta S^3\cdot \Z_n \,\, \supset \,\, \Delta S^1 \cdot \Z_n$\\
 & where $n = 1$ or 2\\
\hline
$N^6_D$ & $S^3\times S^3 \,\, \supset \,\, T^2, \,\, S^3\times S^1 \,\, \supset \,\, \set{(e^{ip\theta},e^{i\theta})}$\\
\hline
$Q^6_B$ & $S^3\times S^3 \,\, \supset \,\, \Delta S^3, \,\, \Delta S^3 \,\, \supset \,\, \Delta S^1$\\
\hline
$Q^6_C$ & $S^3\times S^3 \,\, \supset \,\, \Delta S^3, \,\, S^3\times S^1 \,\, \supset \,\, \Delta S^1$\\
\hline
$N^6_E$ & $S^3\times S^3 \,\, \supset \,\, S^3\times S^1, \,\, S^3\times S^1 \,\, \supset \,\, \set{(e^{ip\theta},e^{i\theta})}$\\
\hline
$Q^6_D$ & $S^3\times S^3 \,\, \supset \,\, S^3\times S^1, \,\, S^1\times S^3 \,\, \supset \,\, \Delta S^1$\\
\hline
$N^6_F$ & $\SU(3) \,\, \supset \,\, \S(\U(2)\U(1)), \,\, \S(\U(2)\U(1)) \,\, \supset \,\, \SU(2)\SU(1)\cdot \Z_n$\\
\hline
\end{tabular}
\end{center}
\vspace{0.1cm}
\caption{Nonreducible 6-dimensional diagrams which are not products, sums or fixed point actions.}\label{6dimlist}
\end{table}}

{\setlength{\tabcolsep}{0.20cm}
\renewcommand{\arraystretch}{1.6}
\stepcounter{equation}
\begin{table}[!h]
\begin{center}
\begin{tabular}{|c|c|}
\hline
$N^7_A$ & $S^3\times S^3 \,\, \supset \,\, \set{(e^{ip_-\theta},e^{iq_-\theta})}\cdot H_+, \,\, \set{(e^{ip_+\theta},e^{iq_+\theta})} \cdot H_- \,\, \supset \,\,  H_-\cdot H_+$\\
 & $H_\pm = \Z_{n_\pm}\subset K^\pm_0$\\
\hline
$N^7_B$ & $S^3\times S^3 \,\, \supset \,\, \set{(e^{ip\theta},e^{iq\theta})}\cdot H_+, \,\, \set{(e^{j\theta},1)}\cdot H_- \,\, \supset \,\,  H_-\cdot H_+$\\
 & where $H_\pm = \Z_{n_\pm} \subset K^\pm_0$, $n_+\le 2$, $4|n_-$ and $p_- \equiv \pm \frac{n_-}{4} \mod n_-$\\
\hline
$P^7_A$ & $S^3\times S^3 \,\, \supset \,\, \set{(e^{ip_-\theta},e^{iq_-\theta})}, \,\, \set{(e^{jp_+\theta},e^{jq_+\theta})}\cdot H \,\, \supset \,\,  \langle (i,i) \rangle$\\
 & where $p_-,q_-\equiv 1 \mod 4$\\
\hline
$P^7_B$ & $S^3\times S^3 \,\, \supset \,\, \set{(e^{ip_-\theta},e^{iq_-\theta})} \cdot H, \,\, \set{(e^{jp_+\theta},e^{jq_+\theta})}\cdot H \,\, \supset \,\,  \langle (i,i), (1,-1) \rangle$\\
 & where $p_-,q_-\equiv 1 \mod 4$, $p_+$ even\\
\hline
$P^7_C$ & $S^3\times S^3 \,\, \supset \,\, \set{(e^{ip_-\theta},e^{iq_-\theta})} \cdot H, \,\, \set{(e^{jp_+\theta},e^{jq_+\theta})}\cdot H \,\, \supset \,\,  \Delta Q$\\
 & where $p_\pm,q_\pm\equiv 1 \mod 4$\\
\hline
$N^7_C$ & $S^3\times S^3 \,\, \supset \,\, \set{(e^{ip\theta},e^{iq\theta})}, \,\, S^3\times \Z_n \,\, \supset \,\,  \Z_n$\\
 & where $(q,n)=1$\\
\hline
$P^7_D$ & $S^3\times S^3 \,\, \supset \,\, \set{(e^{ip\theta},e^{iq\theta})}, \,\, \Delta S^3 \cdot \Z_n \,\, \supset \,\,  \Z_n$\\
 & $n=2$ and $p$ even; or $n=1$ and $p$ arbitrary\\
\hline
$Q^7_A$ & $S^3\times S^3 \,\, \supset \,\, S^3\times 1, \,\, \Delta S^3  \,\, \supset \,\,  1$\\
\hline
$Q^7_B$ & $S^3\times S^3 \,\, \supset \,\, \Delta S^3, \,\, \Delta S^3  \,\, \supset \,\,  1$\\
\hline
$N^7_D$ & $S^3\!\!\times\! S^3\!\!\times\! S^1 \supset  \set{(z^pw^{\lambda m}\!,z^qw^{\mu m}\!,w)},   \set{(z^pw^{\lambda m}\!,z^qw^{\mu m}\!,w)}   \supset  H_0\cdot\Z_n $\\
 & where $H_0=\set{(z^p,z^q,1)}$, $p\mu-q\lambda=1$ and $\Z_n\subset \set{(w^{\lambda m}\!,w^{\mu m}\!,w)}$\\
\hline
$N^7_E$ & $S^3\!\!\times\! S^3\!\!\times\! S^1\! \supset\!  \set{\!(z^pw^{\lambda m_-}\!,z^qw^{\mu m_-}\!,w^{n_-})\!}\! H,   \set{\!(z^pw^{\lambda m_+}\!,z^qw^{\mu m_+}\!,w^{n_+})\!}\! H \! \supset \! H$\\
 & where $H=H_-\cdot H_+$, $H_0=\set{(z^p,z^q,1)}$, $K^-\ne K^+$, $p\mu-q\lambda=1$,\\
 & $\gcd(n_-,n_+,d)=1$ where $d$ is the index of $H\cap K^-_0 \cap K^+_0$ in $K^-_0 \cap K^+_0$\\
\hline
\end{tabular}
\end{center}
\vspace{0.1cm}
\caption{Nonreducible 7-dimensional diagrams which are not products, sums or fixed point actions. (1 of 2).}\label{7dimlist1}
\end{table}}
{\setlength{\tabcolsep}{0.2cm}
\renewcommand{\arraystretch}{1.6}
\stepcounter{equation}
\begin{table}[!h]
\begin{center}
\begin{tabular}{|c|c|}
\hline
$Q^7_C$ & $S^3\times S^3\times S^1 \,\, \supset \,\, \set{(e^{i\phi},e^{ib\theta},e^{i\theta})}, \,\, S^3\times 1\times 1\cdot \Z_n  \,\, \supset \,\, S^1\times 1\times 1 \cdot \Z_n$\\
 & $\Z_n\subset \set{(1,e^{ib\theta},e^{i\theta})}$\\
\hline
$Q^7_D$ & $S^3\times S^3\times S^1 \,\, \supset \,\, \set{(e^{i\phi},e^{i\phi}e^{ib\theta},e^{i\theta})}, \,\, \Delta S^3\times 1\cdot \Z_n  \,\, \supset \,\, \Delta S^1\times 1\cdot \Z_n$\\
 & $\Z_n\subset \set{(1,e^{ib\theta},e^{i\theta})}$ where $n$ is 1 or 2\\
\hline
$N^7_F$ & $S^3 \! \times \! S^3 \! \times \! S^1 \,\, \supset \,\, \set{(e^{ip\phi}e^{ia\theta},e^{i\phi},e^{i\theta})}, \,\, S^3\! \times \! S^1\! \times\! \Z_n  \,\, \supset \,\, \set{(e^{ip\phi},e^{i\phi},1)}\cdot \Z_n$\\
 & $\Z_n\subset \set{(e^{ia\theta},1,e^{i\theta})}$\\
\hline
$N^7_G$ & $\SU(3) \,\, \supset \,\, \S(\U(1)\U(2)), \,\, \S(\U(1)\U(2)) \,\, \supset \,\, T^2$\\
\hline
$Q^7_E$ & $\SU(3) \,\, \supset \,\, \S(\U(1)\U(2)), \,\, \S(\U(2)\U(1)) \,\, \supset \,\, T^2$\\
\hline
$Q^7_F$ & $\SU(3)\times S^1 \,\, \supset \,\, \set{(\beta(m\theta),e^{i\theta})}\cdot H_0, \,\, \set{(\beta(m\theta),e^{i\theta})}\cdot H_0  \,\, \supset \,\, H_0\cdot \Z_n$\\
 & $H_0 = \SU(1)\SU(2)\times 1$, $\Z_n\subset \set{(\beta(m\theta),e^{i\theta})}$,\\
 & $\beta(\theta)= \diag(e^{-i\theta},e^{i\theta},1)$\\
\hline
$N^7_H$ & $\SU(3)\times S^1 \,\, \supset \,\, \set{(\beta(m_-\theta),e^{in_-\theta})}\cdot H, \,\, \set{(\beta(m_+\theta),e^{in_+\theta})}\cdot H  \,\, \supset \,\, H$\\
 & $H_0 = \SU(1)\SU(2)\times 1$, $H=H_-\cdot H_+$, $K^-\ne K^+$,\\
 & $\beta(\theta)= \diag(e^{-i\theta},e^{i\theta},1)$, $\gcd(n_-,n_+,d)=1$\\
 & where $d$ is the index of $H\cap K^-_0 \cap K^+_0$ in $K^-_0 \cap K^+_0$\\
\hline
$Q^7_G$ & $\SU(3)\times S^1 \,\, \supset \,\, \set{(\beta(m\theta),e^{i\theta})}\cdot H_0, \,\, \SU(3)\times \Z_n  \,\, \supset \,\, H_0\cdot \Z_n$\\
 & $H_0 = \SU(1)\SU(2)\times 1,\,\, \Z_n\subset \set{(\beta(m\theta),e^{i\theta})}$,\\
 & $\beta(\theta)= \diag(e^{-i\theta},e^{i\theta},1),\,\, \gcd(m,n)=1$\\
\hline
$N^7_I$ & $\Sp(2)\,\, \supset \,\, \Sp(1)\Sp(1), \,\, \Sp(1)\Sp(1) \,\, \supset \,\, \Sp(1)\SO(2)$\\
\hline
\end{tabular}
\end{center}
\vspace{0.1cm}
\caption{Nonreducible 7-dimensional diagrams which are not products, sums or fixed point actions (2 of 2).}\label{7dimlist2}
\end{table}}

{\setlength{\tabcolsep}{0.40cm}
\renewcommand{\arraystretch}{1.6}
\stepcounter{equation}
\begin{table}[!h]
\begin{center}
\begin{tabular}{|c|c|}
\hline
\multicolumn{2}{|c|}{Sum actions}\\
\hline
$S^n$  & $G_i/H_i\approx S^{m_i}\subset \R^{m_i+1}$, $i=1,2$.\\
       & $G_1\times G_2$ on $S^{m_1+m_2+1}\subset \R^{m_1+1}\times \R^{m_2+1}$:\\
       & $(g_1,g_2)\star (x,y)=(g_1\cdot x,g_2\cdot y)$\\
\cline{2-2}
     & $G_1\times G_2 \,\, \supset \,\, G_1\times H_2, \,\, H_1\times G_2  \,\, \supset \,\, H_1\times H_2$\\
\hline
\hline
\multicolumn{2}{|c|}{Product actions}\\
\hline
$M\times L/J$  & $M$ cohomogeneity one for $G\supset K^-,K^+ \supset H$,\\
               & $L/J$ homogeneous.\\
     &  $G\times L$ on $M\times L/J$:\\
     &  $(g,l)\star (p, \ell J)=(g\cdot p, l\ell J)$\\
\cline{2-2}
     & $G\times L \,\, \supset \,\, K^-\times J, \,\, K^+\times J \,\, \supset \,\, H\times J$\\
\hline
\hline
\multicolumn{2}{|c|}{Fixed-point actions}\\
\hline
CROSS  &  See Section \ref{special_types}\\
\cline{2-2}
     & $G \,\, \supset \,\, G, \,\, K \,\, \supset \,\, H$\\
\hline
\end{tabular}
\end{center}
\vspace{0.1cm}
\caption{Sum, product and fixed point actions.}\label{SaP}
\end{table}}

%
%

%
%
%
%
%
%

\end{document}